\documentclass[11pt]{article}

\usepackage{graphicx}
\usepackage[dvipsnames]{xcolor}
\usepackage{caption}
\usepackage{bm}
\usepackage{comment}
\usepackage{psfrag}
\usepackage{latexsym,amsmath,amsfonts,amscd, amsthm}
\usepackage{changebar}
\usepackage{color}
\usepackage{bm}
\usepackage{tikz}
\usepackage{multirow}
\usepackage{subfigure}

\usetikzlibrary{arrows,backgrounds,snakes}
\usepackage[title]{appendix}

\usepackage{chngcntr,colonequals}
\usepackage{soul}
\newtheoremstyle{plainNoItalics}{}{}{\normalfont}{}{\bfseries}{.}{ }{}
\usepackage{amssymb}
\usepackage{bbm}

\theoremstyle{plain}
\newtheorem{thm}{Theorem}[section]

\theoremstyle{plainNoItalics}

\newtheorem{defn}[thm]{Definition}
\newtheorem{rem}[thm]{Remark}
\newtheorem{prop}[thm]{Proposition}

\newcommand{\beq}{\begin{equation}}
\newcommand{\eeq}{\end{equation}}
\newcommand{\beqa}{\begin{eqnarray}}
\newcommand{\eeqa}{\end{eqnarray}}
\newcommand{\bit}{\begin{itemize}}
\newcommand{\eit}{\end{itemize}}
\newcommand{\bedef}{\begin{defn}}
\newcommand{\edefn}{\end{defn}}
\newcommand{\bpro}{\begin{prop}}
\newcommand{\epro}{\end{prop}}

\newcommand{\mE}{\mathcal E}

\newcommand{\mL}{\mathcal{L}}

\newcommand\pzc[1]{\textcolor{black}{#1}}

\newcommand{\vareps}{\varepsilon}

\newcommand{\lgl}{\langle}
\newcommand{\rgl}{\rangle}

\newcommand{\half}{\frac{1}{2}}
\newcommand{\BX}{\mathbf{X}}
\newcommand{\BOmega}{{\boldsymbol{\Omega}}}
\newcommand{\BF}{\mathbf{f}}
\newcommand{\Bx}{\mathbf{x}}
\newcommand{\BD}{\mathbf{D}}
\newcommand{\BRHO}{{\bm{\rho}} }

\newcommand{\Bg}{{\bf{g}} }

\newcommand{\BSIGMA}{{\bf{\Sigma}} }

\newcommand{\BC}{{\bf{c}} }
\newcommand{\BFF}{{\bf{F}} }
\newcommand{\BH}{{\bf{H}} }
\newcommand{\BU}{{\bf{U}} }
\newcommand{\WBOmega}{\widetilde{{\boldsymbol{\Omega}}}}

\newcommand{\wtheta}{\widetilde{\theta}}

\newcommand{\mP}{\mathcal{P}}
\newcommand{\mU}{\mathcal{U}}

\newcommand{\mT}{\mathcal{T}}
\newcommand{\BMU}{\boldsymbol{\mathcal{U}}}

\newcommand{\BMV}{\boldsymbol{\mathcal{V}}}
\newcommand{\BLambda}{\boldsymbol{\Lambda}}

\newcommand{\WOmega}{\widetilde{\Omega}}

\newcommand{\wa}{\widetilde{a}}
\newcommand{\wb}{\widetilde{b}}

\usepackage{algcompatible}
\usepackage[ruled,vlined]{algorithm2e}

\newcommand\fli[1]{\textcolor{black}{#1}}

\counterwithin*{equation}{section}
\counterwithout{figure}{section}

\begin{document}
\baselineskip=1.2pc
\title{A reduced basis method for radiative transfer equation\footnote{
This material is based upon work supported by the National Science Foundation under Grant No. DMS-1439786 and by the Simons Foundation Grant No. 50736 while the authors were in residence at the Institute for Computational and Experimental Research in Mathematics in Providence, RI, during the ``Model and dimension reduction in uncertain and dynamic systems" program.
}}

\author{
Zhichao Peng\thanks{Department of Mathematics, Michigan State University, East Lansing, MI 48824 U.S.A. Email: {\tt pengzhic@msu.edu}.}
\and
Yanlai Chen\thanks{Department of Mathematics, University of Massachusetts Dartmouth, 285 Old Westport Road, North Dartmouth, MA 02747, USA. Email: {\tt{yanlai.chen@umassd.edu}}. Research is partially supported by National Science Foundation grant DMS-1719698, and by the UMass Dartmouth Marine and UnderSea Technology (MUST) Research Program made possible via an Office of Naval Research grant N00014-20-1-2849.}
\and
Yingda Cheng \thanks{Department of Mathematics, Department of  Computational Mathematics, Science and Engineering, Michigan State University,
East Lansing, MI 48824 U.S.A.
Email: {\tt ycheng@msu.edu}. Research is supported by NSF grants    DMS-2011838 and AST-2008004.}
\and
Fengyan Li\thanks{Department of Mathematical Sciences, Rensselaer Polytechnic Institute, Troy, NY 12180, U.S.A. Email: 
{\tt lif@rpi.edu}. Research is supported by NSF grants DMS-1719942 and DMS-1913072.}}
\maketitle

\abstract{
Linear kinetic transport equations play a critical role in optical tomography, radiative transfer and neutron transport.  The fundamental difficulty hampering their efficient and accurate numerical resolution lies in the high dimensionality of the physical and velocity/angular variables and the fact that the problem is multiscale in nature. Leveraging the existence of a hidden low-rank structure hinted by the diffusive limit, in this work, we design and test the  angular-space reduced order model for the  linear radiative transfer equation, the first such effort based on the celebrated reduced basis method (RBM). 

Our method is built upon a high-fidelity solver employing the discrete ordinates method in the angular space, an asymptotic preserving  upwind discontinuous Galerkin method for the physical space, and an efficient synthetic accelerated source iteration for the resulting linear system. Addressing the challenge of the parameter values (or angular directions) being coupled through an integration operator, the first novel ingredient of our method is an iterative procedure  where the macroscopic density is constructed from the RBM snapshots, treated explicitly and allowing a transport sweep, and then updated afterwards. A greedy algorithm can then proceed to adaptively select the representative samples in the angular space and  form a surrogate solution space. The second novelty is a least squares density reconstruction strategy, at each of the relevant physical locations, enabling the robust and accurate integration over an arbitrarily unstructured set of angular samples toward the macroscopic density. Numerical experiments indicate that our method is \pzc{effective} for computational cost reduction in a variety of regimes. 
}

\section{Introduction}
\label{sec:introduction}
Linear kinetic transport equations model particles propagating through, and interacting with, background media.  They provide prototype models for optical tomography \cite{arridge2009optical}, radiative transfer \cite{mihalas2013foundations,pomraning2005equations} and neutron transport \cite{lewis1984computational}. In this work, we consider the following steady-state {linear radiative transfer equation}
\begin{subequations}
\label{eq:kinetic}
\begin{align}
&\BOmega\cdot\nabla f= \sigma_s\lgl f\rgl -\sigma_t f + G, \quad\forall \;\Bx\in\BX,\;\BOmega \in \mathbb{S}^{d-1},\\
\intertext{whose solution $f=f(\Bx,\BOmega)$ delineates particle distribution at location $\Bx$ in a physical domain $\BX\subset\mathbb{R}^d$, and $\BOmega\in \mathbb{S}^{d-1}$ denotes the angular variable. We enforce a Dirichlet boundary condition on the inflow boundary}
&f(\Bx,\BOmega) = f_\textrm{inflow}(\Bx,\BOmega),\; \Bx\in \partial \BX\;\;\text{and}\;\;\BOmega\cdot{\bf{n}}(\Bx)<0.
\end{align}
\end{subequations}
{Here ${\bf{n}}(\Bx)$ stands for the unit outward normal on $\partial\BX$. In addition, $\sigma_s(\Bx)\geq 0$ is the scattering cross section, $\sigma_t(\Bx)=\sigma_s(\Bx)+\sigma_a(\Bx)$ is the total cross section, with $\sigma_a(\Bx)\ge 0$ being the absorption cross section, 
and  $G(\Bx)$ is the source. $\lgl\cdot\rgl$ encodes a normalized integration in the angular space, namely,
 \begin{align}
\label{eq:integral_operator}
 \lgl f\rgl = \frac{1}{|\mathbb{S}^{d-1}|}\int_{\mathbb{S}^{d-1}}f d\BOmega.
\end{align}
This gives the macroscopic density, $\rho(\Bx)=\lgl f\rgl$, defined on $\BX$.
}
A fundamental difficulty in numerically resolving \eqref{eq:kinetic} originates from the high dimension ($2 d - 1$) of $(\Bx, \BOmega)$ space.  

Given  the non-dimensional mean free path length $\vareps>0$, one can further define  the rescaled cross sections as 
${\sigma}_s=\frac{\widehat{\sigma}_s}{\vareps}$, ${\sigma}_a=\vareps\widehat{\sigma}_a$, and $G=\vareps\widehat{G}$.
As proved in \cite{bardos1984diffusion}, when $\vareps\rightarrow0$, $f(\Bx,\BOmega)\rightarrow \rho(\Bx)$, and  $\rho(\Bx)$ solves  the diffusion equation
\pzc{\begin{align}
\label{eq:diffusion_limit}
\nabla\cdot\left( \widehat{\sigma}_s^{-1} D\nabla\rho\right)=\widehat{\sigma}_a\rho+\widehat{G}.
\end{align}
Here $D = \textrm{diag}\left(\lgl \Omega_1^2\rgl,\dots,\lgl\Omega_d^2\rgl\right)$, with $\Omega_i$ being the $i$-th component of $\BOmega\in\mathbb{S}^{d-1}$.
} As a result, when the problem is in its diffusive regime, {the solution is close to a rank-1 manifold} in the angular space (i.e. the diffusion limit). Hence, {in or near such  regime,} 
it is possible to design a reduced order model (ROM) to capture the low-rank structure of the solution, and this   provides an opportunity to {mitigate} the curse of dimensionality. On the other hand, when  $\vareps$ varies  and possibly spans through several magnitude in the physical domain, the problem is multiscale in its nature. This poses challenges not only for traditional numerical schemes but also for any attempt to design an effective ROM. 
{It is well known that a standard numerical method may fail to capture the diffusion limit and lead to unphysical solutions when spatial meshes are under-resolved with respect to $\vareps$, i.e. when the mesh size $h$ satisfies $h\gg \vareps$ \cite{larsen1987asymptotic,naldi1998numerical}.
}
 This issue should be addressed before designing any ROM. One solution is to apply the asymptotic preserving (AP) schemes \cite{jin2010asymptotic}, which preserve the asymptotic limit on the discrete level and work well for both the kinetic regime and the diffusive regime. In other words,  AP methods are able to capture the diffusion limit even on  under-resolved meshes. 
In this paper, we apply the discrete ordinates ($S_N$) method \cite{pomraning1973equations} in the angular space and the upwind discontinuous Galerkin (DG) method in physical space.   
{Upwind DG method was investigated 
for \eqref{eq:kinetic} in \cite{larsen1987asymptotic,larsen1989asymptotic, adams2001discontinuous}, and it is proved to be AP \cite{guermond2010asymptotic} when the approximation space  contains  continuous functions 
that are at least linear on each mesh element.} {To solve the resulting linear system, we apply iterative solvers based on source iteration \cite{lewis1984computational}. Standard source iteration converges slowly when the problem is scattering dominant \cite{adams2002fast}, we here apply the synthetic accelerated source iteration (SASI) \cite{adams2002fast} to achieve good efficiency for various regimes.}

As discussed above, the existence of a low-rank {structure of the solution manifold in the diffusive regime}  
gives hope to the success of ROM techniques for kinetic transport equations \eqref{eq:kinetic}, the focus of our paper. This is by no means the first such attempt.  In fact,  ROM has become an increasingly popular technique in kinetic simulations in the last few years. 
The proper orthogonal decomposition (POD) type techniques were applied to \eqref{eq:kinetic} and its time transient counterpart in \cite{buchan2015pod,behne2019model,coale2019reduced}. Other related works include space-time POD \cite{choi2020space}, AP random singular value decomposition (RSVD) \cite{chen2020random}, dynamic mode decomposition (DMD) \cite{mcclarren2019calculating}, 
proper generalized decomposition (PGD) \cite{dominesey2018reduced,dominesey2019reduced,prince2019separated,alberti2020reduced}, and dynamical low rank approximations (DLRA) \cite{ding2019error, einkemmer2020asymptotic,peng2020low}. 
These existing approaches either lack the hallmark greedy algorithm or is not a method of snapshots. A POD method with the greedy algorithm to adaptively select angular samples is proposed in \cite{tencer2016reduced}, but \cite{tencer2016reduced} only considers problems without the scattering effect.
In this work, taking the scattering effect into account, we design and test a model reduction technique for the linear radiative transfer equation \eqref{eq:kinetic} that features both ingredients, namely a celebrated greedy algorithm adaptively selecting the representative samples in the angular space and a resulting surrogate solution space spanned by the corresponding snapshots. Indeed, it is under the framework of  the reduced basis (RB) method \cite{quarteroni2011certified,quarteroni2015reduced,hesthaven2016certified,haasdonk2017reduced} and to the best of our knowledge, the very first such attempt, which takes the scattering effect into account.

RB method (RBM) has become the go-to option for efficiently simulating parametric partial differential equation (PDE). Its hallmark feature is a greedy algorithm embedded in an offline-online decomposition procedure. 
The offline (i.e. training) stage is devoted to a judicious exploration of the parameter-induced solution manifold. It adaptively selects a number of representative parameter values via a mathematically rigorous greedy algorithm \cite{BinevCohenDahmenDevorePetrovaWojtaszczyk}. Solution snapshots for these parameter values are then obtained through a user-specified (potentially expensive) full order (i.e. accurate) solver.
 The buildup of the surrogate solution space spanned by these snapshots is done step-by-step. Each iteration of the greedy algorithm adds the parameter value which is the maximizer of a mathematically rigorous {\em a posteriori} error estimator or an effective error indicator were the current surrogate space used as a reduced solver space via e.g. a Galerkin or Petrov-Galerkin projection. 
For parametric systems bearing a small Kolmogorov N-width \cite{pinkus_n-widths_1985}, the dimension of the surrogate space is orders of magnitude smaller than the total degrees of freedom for the full model in order to reach a high degree of accuracy.  
This difference in size leads to a dramatic decrease in computation time for the online simulations when a {\em reduced} solution is sought in the terminal surrogate space for each parameter value. Moreover, unlike other ROM techniques (e.g. POD-based approaches), the number of full order inquiries RBM takes offline is minimum i.e. equal to the dimension of the surrogate space.

Leveraging the low-rank structure induced by the angular space, it is natural to treat the angular variable $\BOmega$ as our parameter. 
What prevents a direct application of RBM is the integral operator \eqref{eq:integral_operator} which poses a two-fold extra challenge. 
First, unlike the standard setting when RBM applies, the solutions for different parameter values $\{f(\cdot,\BOmega): \BOmega \in {\mathbb S}^{d-1}\}$ are all coupled through \eqref{eq:integral_operator}. A milder version of this type of coupling is addressed in the stochastic PDE setting \cite{LiuChenChenShu2019}. 
Second, from a practical viewpoint, the robustness and the efficiency of the SASI iterative linear solver is highly sensitive to the quality of  the density approximation $\rho=\lgl f\rgl$.  In particular, the unstructured nature of the selected RB parameter samples {in the angular variable} prevents a robust and accurate numerical integration which usually requires a structured set of quadrature points.

To address
these challenges, we effectively decouple the solutions for different angular samples by designing {a greedy} iterative procedure where an approximation of the macroscopic density $\rho$ is constructed from the RB snapshots at the beginning of each iteration. It is frozen during the iteration and reconstructed before the next. 
To resolve the lack of structure for a robust numerical integration, we develop a least squares density reconstruction strategy capable of integrating over an arbitrary set of selected angular samples {robustly}. Our numerical experiments show that the proposed method is effective for the radiative transfer problems in the scattering dominant and intermediate regimes, as well as for the multiscale problems with large scattering dominant subregions. \pzc{The proposed RB method (including both the offline and online stages) can be seen as  a surrogate model for \fli{the full order upwind DG solver}. The online stage of the proposed method can be further applied to predict solutions at \fli{``unseen'' angular samples.} The proposed method can also be utilized as a building block to construct ROMs for problems with essential physical parameters such as the magnitude of the scattering cross section whose ``multi-query'' nature will lead to more pronounced saving for our approach.}

The rest of the paper is organized as follows. We introduce the full order numerical scheme and  the SASI iterative solver in Section \ref{sec:full_order}. Section \ref{sec:RB} is devoted to the least squares density reconstruction and  the RB algorithm. We demonstrate the performance of the proposed method through a series of  numerical experiments in Section \ref{sec:numerical}. Finally, conclusions are made in Section \ref{sec:conclusion}.

\section{Full order numerical method}\label{sec:full_order}
In this section, we describe our full order numerical scheme for solving \eqref{eq:kinetic} focusing on the 1D slab geometry and the 2D case with $\BOmega\in \mathbb{S}^1$. 
{The equation \eqref{eq:kinetic} on the 1D slab geometry is given as 
\begin{align}
v\partial_x f = \frac{\sigma_s}{2}\int_{-1}^1 f dv-\sigma_t f+G,\quad v\in[-1,1],\quad x\in\BX,
\end{align}
where we write $\Omega$ as $v$ here following convention.}
The 2D equation with $\BOmega\in \mathbb{S}^1$ can be written as
\begin{align}
\cos(\theta)\partial_x f+\sin(\theta)\partial_y f= \sigma_s\lgl f\rgl -\sigma_t f + G, \quad\theta\in[0,2\pi],\quad \Bx\in\BX.
\end{align}
In the next three subsections, we detail  the discrete ordinates method \cite{pomraning1973equations} for the angular discretization, the upwind DG method  for the spatial discretization, and finally the  SASI iterative solver.

\subsection{Angular discretization}
In the angular space, we apply the discrete ordinates method \cite{pomraning1973equations} by sampling $f$ at quadrature points  $\{\BOmega_j\}_{j=1}^{N_\BOmega}$. Letting $\{\omega_j\}_{j=1}^{N_\BOmega}$ be the corresponding normalized quadrature weights, we can discretize the integral operator as
\begin{align*}
\lgl f \rgl \approx \lgl f\rgl_h=\sum_{j=1}^{N_{\BOmega}} \omega_j f(\cdot,\BOmega_j),\quad \sum_{j=1}^{N_\BOmega}\omega_j=1,
\end{align*}
and, as a consequence, equation \eqref{eq:kinetic} can be discretized as
\begin{align}
\label{eq:sn}
\BOmega_j\cdot\nabla f_j= \sigma_s\sum_{i=1}^{N_\BOmega} \omega_i f_i -\sigma_t f_j + G,\quad j=1,\dots,N_\BOmega,
\end{align}
where \fli{$f_j(\Bx) \approx f(\Bx,\BOmega_j)$. }
This method is also referred to as the $S_N$ method, if one use $N_\Omega=N$ quadrature points  for the 1D slab geometry or $N_{\BOmega}=2N$ quadrature points for the 2D case. Particularly, we use the following quadrature rules.
\begin{enumerate}
\item [] {\bf 1D}:\;\;  $\{\BOmega_j\}_{j=1}^N$ is the collection of  $N$-point Gauss-Legendre quadrature points on $[-1,1]$.
\item []{\bf 2D}:\;\; 
$\BOmega_j = \left(\cos(\theta_j),\sin(\theta_j)\right)^T,\;\text{with}\; \theta_j = \frac{(2j-1)\pi}{2N},\;\omega_j=\frac{1}{2N},\; j=1,\dots,N_\BOmega.$
\end{enumerate}

\subsection{Spatial discretization}
It is well known that the  asymptotic preserving (AP) schemes \cite{jin2010asymptotic, naldi1998numerical} can capture the correct diffusion limit without a highly refined mesh resolving the small $\vareps$-scale. 
We adopt the upwind DG method \cite{larsen1987asymptotic,larsen1989asymptotic,guermond2010asymptotic}, which has been proven to be AP {if the approximation space contains continuous functions that are at least linear on each mesh element
\cite{adams2001discontinuous,guermond2010asymptotic}.  Without loss of generality, we assume $\BX=[x_L, x_R]$ in 1D and $\BX=[x_L, x_R]\times[y_L, y_R]$  in 2D. 
Let $\BX_h = \{T_i,i=1,\dots,N_x\}$ be a partition of $\BX$, with each $T_i$ being an interval in 1D or a rectangle in 2D.} We introduce a discrete space
\begin{align}
U_h^K=\{u(\Bx): u(\Bx)|_{T_i}\in Q^K(T_i),i=1,\dots,N_x\},\; K\geq 1,
\label{DG:space}
\end{align}
where $Q^K(T_i)$ denotes the $d$-variate polynomial on $T_i$ with degree up to $K$ for each variable. We are now ready to state the upwind DG spatial discretization:  we seek $\big\{f_h(\Bx,\BOmega_j)\in U_h^K: \,\, j=1,\dots,N_{\BOmega}\big\}$, such that 
\begin{align}
-\int_{T_i} \Big(\BOmega_j\cdot\nabla\phi_h(\Bx)\Big) f_h(\Bx,\BOmega_j) d\Bx&+ \int_{\partial T_i} \widehat{\BH}(\BOmega_j,f_h, {\bf{n}}_i)\phi_h(\Bx) ds+\int_{T_i}\sigma_t f_h(\Bx,\BOmega_j)\phi_h(\Bx) d\Bx
\notag
\\
\fli{= \sum_{k=1}^{N_\BOmega}\omega_k\int_{T_i}\sigma_s f_h(\Bx,\BOmega_k)\phi_h(\Bx) d\Bx}
& + \int_{T_i}G(\Bx)\phi_h(\Bx) d\Bx, \quad \forall\phi_h\in U_h^K.
\label{eq:DG}
\end{align}
Here $\widehat{\BH}(\BOmega_j, f_h,{\bf{n}}_i)$ is the upwind numerical flux along $\partial T_i$, that is defined, for an element $T_i=T^-$ with the neighboring element $T^+$, as
\begin{align}
\widehat{\BH}(\BOmega_j, f_h,{\bf{n}}_i) = \frac{\BOmega_j\cdot{\bf{n}}_i}{2}\Big(f_h^+(\Bx,\BOmega_j)+f_h^-(\Bx,\BOmega_j)\Big)+\frac{|\BOmega_j\cdot {\bf{n}}_i|}{2}\Big(f_h^-(\Bx,\BOmega_j)-f_h^+(\Bx,\BOmega_j)\Big).
\label{eq:upwind}
\end{align}
We use $f_h^{\pm}$ to denote the restriction of $f_h$ to $T^{\pm}$, while ${\bf{n}}_i$ is  the unit outward normal on $\partial T_i$. 
Based on $f_h$, the density $\rho$ is further approximated by 
\begin{align}
\rho_h(\Bx)=\lgl f_h\rgl_h = \sum_{j=1}^{N_\BOmega}\omega_j f_h(\Bx,\BOmega_j).
\label{eq:quadrature_relation}
\end{align}

Next we will rewrite the DG scheme \eqref{eq:DG} into its matrix-vector form. To this end, 
 we assume that $\{\phi_k\}_{k=1}^{ N_{\textrm{dof}} }$ is a basis for $U_h^K$, 
 and  $f_h(\Bx,\BOmega_j)$, $\rho_h(\Bx)$ are then expanded as
\[
f_h(\cdot,\BOmega_j) = \sum_{k=1}^{N_\textrm{dof}}\alpha^f_k(\BOmega_j)\phi_k, \quad \rho_h = \sum_{k=1}^{N_\textrm{dof}}\alpha^\rho_k\phi_k, \mbox{ with } \alpha^\rho_k = \sum_{j=1}^{N_\BOmega}\omega_j\alpha_k^f(\BOmega_j). 
\]
The last equality is due to the numerical quadrature \eqref{eq:quadrature_relation}. We further define 
\begin{align*}
&\BF_{\BOmega_j}=\left(\alpha^f_1(\BOmega_j),\dots,\alpha^f_{N_\textrm{dof}}(\BOmega_j)\right)^T\in \mathbb{R}^{N_\textrm{dof}},\;\;\BFF=\left(\BF_{\BOmega_1},\dots,\BF_{\BOmega_{N_\Omega}}\right)\in \mathbb{R}^{N_\textrm{dof}\times N_\BOmega}, 
\\ 
&\BRHO=\left(\alpha^\rho_1,\dots,\alpha^\rho_{N_\textrm{dof}}\right)^T\in \mathbb{R}^{N_\textrm{dof}},
\end{align*}
and $\BU_{\BOmega_j},j=1,\dots,N_{\BOmega}$, $\BSIGMA_t,\BSIGMA_s\in \mathbb{R}^{N_\textrm{dof}\times N_\textrm{dof}}$ and $\Bg\in \mathbb{R}^{N_\textrm{dof}}$ as
\begin{align*}
&(\BU_{\BOmega_j})_{kl} = \sum_{i=1}^{N_x}\left(-\int_{T_i} (\BOmega_j\cdot\nabla\phi_k(\Bx)) \phi_l(\Bx) d\Bx+ \int_{\partial T_i} \widehat{\BH}\left(\BOmega_j, \phi_l,{\bf{n}}_i\right)\phi_k(\Bx) ds\right),
\\
&(\BSIGMA_t)_{kl} = \sum_{i=1}^{N_x}\left(\int_{T_i} \sigma_t\phi_k(\Bx) \phi_l(\Bx) d\Bx\right),
\quad (\BSIGMA_s)_{kl} = \sum_{i=1}^{N_x}\left(\int_{T_i} \sigma_s\phi_k(\Bx) \phi_l(\Bx) d\Bx\right),
\\
&(\Bg)_k= \sum_{i=1}^{N_x}\left(\int_{T_i} G(\Bx)\phi_k(\Bx) d\Bx\right).
\end{align*}
\fli{We here adopt the commonly used basis functions $\{\phi_k\}_{k=1}^{ N_{\textrm{dof}} }$, with each being nonzero only on one mesh element as a scaled Legendre polynomial or its tensor version. With such a choice,} 
$\BSIGMA_t$ and $\BSIGMA_s$ are block-diagonal, symmetric and semi-positive definite. If mesh elements in space are suitably recorded, each $\BU_{\BOmega_j}$ can be block lower triangular.
With the notation above, the DG scheme  \eqref{eq:DG} can be rewritten into  its matrix-vector form:
\begin{align}
\BU_{\BOmega_j} \BF_{\BOmega_j} + \BSIGMA_t \BF_{\BOmega_j} - \BSIGMA_s \BRHO = \Bg,\quad \forall j=1,\dots,N_\Omega.
\label{eq:matrix-vector}
\end{align}
We end this subsection by noting that the AP upwind DG scheme \eqref{eq:DG} also exists for unstructured meshes and general geometries.

\subsection{Synthetic accelerated source iteration}
\label{sec:SASI}
Due to the high (i.e. $2d-1$) dimensional nature of the problem,  
iterative methods must be adopted when solving \eqref{eq:matrix-vector}. However, when the problem is scattering dominant, iterative solvers \fli{such as the standard source iterations} 
may converge slowly \cite{adams2002fast}.
To efficiently solve \eqref{eq:matrix-vector}, 
we apply the synthetic accelerated source iteration (SASI) scheme \cite{adams2002fast}. Each iteration of a typical SASI scheme consists of two main steps. The first step is a transport sweep based on the known {$\BRHO^k$} from the previous iteration. More specifically,  using the given data 
{$\BRHO^k$} and $\Bg$, 
we invert $\BU_{\BOmega_j}+\BSIGMA_t$ in \eqref{eq:matrix-vector} 
and obtain $\BF_{\BOmega_j}^{k+1}$ for each $\BOmega_j$. After that, we numerically integrate $f_{\BOmega_j}^{k+1}$ in the angular space to obtain an initial update of the density $\BRHO^{k,*}$. The second step 
is to compute a correction, $\BRHO^{k,c}$, for the density by a computationally less expensive procedure. {One can then  update $\BRHO^{k+1}=\BRHO^{k,*}+\BRHO^{k,c}$ and proceed to the next iteration.} 

In this work, we \pzc{mainly focus on}
the  $S_2$ synthetic acceleration (S2SA) method following \cite{lorence1989s,adams2002fast}. {To elaborate the detail of the second step,}  we assume $f^{k+1}$ is the solution to
\begin{align}
\label{eq:kinetic_f_k}
\BOmega\cdot\nabla f^{k+1}+\sigma_t f^{k+1}= \sigma_s\rho^k  + G,
\end{align} 
and define  $\rho^{k,*}=\lgl f^{k+1}\rgl$. Let $\delta f^{k}=f-f^k$. By subtracting \eqref{eq:kinetic_f_k} from \eqref{eq:kinetic}, we obtain the equation for the correction 
$\delta f^{k+1}$, namely
\begin{align}
\label{eq:kinetic_correction}
\BOmega\cdot\nabla (\delta f^{k+1})+\sigma_t \delta f^{k+1}-\sigma_s\lgl \delta f^{k+1}\rgl =\sigma_s(\rho^{k,*}-\rho^k).
\end{align} 
The correction  of the density can then be calculated as $\rho^{k,c}=\lgl \delta f^{k+1}\rgl$.
The main idea of the S2SA is to apply the $S_2$ approximation in the angular space when solving the correction equation \eqref{eq:kinetic_correction}. This means that we only work with $N_d$ quadrature points \fli{in the angular space}, with $N_d=2$ in 1D and $N_d=4$ in 2D. As a result, the direct solver in the correction step can be implemented very efficiently. 
With the DG spatial discretization, the discretized  linear system for the correction step is
\begin{align}
\BU_{\BOmega_j}^d (\delta\BF_{\BOmega_j}) + \BSIGMA_t (\delta\BF_{\BOmega_j}) - \BSIGMA_s \lgl\delta\BRHO\rgl_h = \BSIGMA_s(\BRHO^{k,*}-\BRHO^k), \;\forall j=1,\dots,N_{d},
\label{eq:matrix-vector-preconditioner}
\end{align}
where $\BU_{\BOmega_j}^d$ is the ``upwind" matrix corresponding to the angular samples of the $S_2$ approximation.  Details of the algorithm are presented in Algorithm \ref{alg:source_iteration}.
\begin{rem}
\pzc{An alternative synthetic acceleration strategy is the diffusion synthetic acceleration (DSA)  \cite{adams1992diffusion,wareing1993new,adams2002fast}, which approximates the correction equation \eqref{eq:kinetic_correction} \fli{through a diffusion model.} It is well known that a so-called ``consistent" discretization must be applied to the diffusion \fli{model}, otherwise the source iteration with the DSA may converge slowly or even diverge 
 \cite{adams1992diffusion,wareing1993new,adams2002fast}. With the S2SA, one can reuse the kinetic solver with fewer angular samples, but it does have more degrees of freedom 
  compared with the DSA. In our reduced order algorithm, we observe that the S2SA is more robust \fli{for different regimes} and slightly more accurate than the DSA, though the full order solvers with  the DSA and S2SA are comparable with respect to the robustness and accuracy. More details of the DSA and the comparison between the RB method with both acceleration strategies can be found in Appendix \ref{sec:DSA}.}
\end{rem}
\begin{algorithm}[H]
\caption{Synthetic accelerated source iteration to solve \eqref{eq:matrix-vector}\label{alg:source_iteration} }
\label{alg:source}
\begin{algorithmic}
\STATE{Given initial guess: $\BRHO^0$ }
\STATE{\textbf{Transport sweep}: solve $(\BU_{\BOmega_j}+\BSIGMA_t) \BF_{\BOmega_j}^1 = \Bg + \BSIGMA_s \BRHO^0$, 
$j=1,\dots,N_\BOmega$.}
\STATE{Reconstruct $\BRHO^{0,*}$ through numerical integration based on $\BF_{\BOmega_j}^1$, $j=1,\dots,N_\Omega$.}
\STATE{\textbf{Correction}: {obtain the correction $\BRHO^{0,c}$}. Let $\BRHO^1=\BRHO^{0,*}+\BRHO^{0,c}$. \fli{Set $k=1$.} }
\WHILE{ $ ||\rho^k-\rho^{k-1}||_\infty>\textrm{error}_{\rm tol}$ and $k\leq \textrm{iter}_{\rm tol}$}
\STATE{\textbf{Transport sweep}: solve $(\BU_{\BOmega_j}+\BSIGMA_t) \BF_{\BOmega_j}^{k+1} = \Bg + \BSIGMA_s \BRHO^k$, $j=1,\dots,N_\BOmega$}
\STATE{Reconstruct $\BRHO^{k,*}$ through numerical integration based on $\BF_{\BOmega_j}^{k+1}$, $j=1,\dots,N_\Omega$.}
\STATE{\textbf{Correction}: {obtain the correction $\BRHO^{k,c}$ by solving \eqref{eq:matrix-vector-preconditioner}.} Let $\BRHO^{k+1}=\BRHO^{k,*}+\BRHO^{k,c}$. }
\STATE{Set $k:= k+1$.}
\ENDWHILE
\end{algorithmic}
\end{algorithm}

\section{Reduced basis method in the angular space}\label{sec:RB}
In this section, we briefly review the basics of RBM and describe the main challenges for designing a RBM for \eqref{eq:kinetic} in Section \ref{sec:rb_bg}. We then present all elements of our algorithm in Section \ref{sec:rb_alg}.

\subsection{Background}
\label{sec:rb_bg}
RBM is a popular approach for obtaining reduced order models for a parametric differential equation
\begin{align}
r(\eta;\mu) = \mathcal{L}(\eta;\mu) -g_{\mu}=0, \; \mu\in \mathcal{X}_\mP,\; \eta(\cdot;\mu)\in \Gamma_\mu,
\label{eq:parametric_pde}
\end{align}
where $\mu\in\mathcal{X}_\mP$ stands for a parameter, $\mL(\cdot;\mu)$ encodes a 
steady-state or time-dependent
parametric differential operator. $\eta(\cdot;\mu) \in \Gamma_{\mu}$ denotes the solution corresponding to $\mu$ and is often referred to as a snapshot. The method assumes a full order (potentially expensive) solver of high accuracy for \eqref{eq:parametric_pde} which, for simplicity, we write in a strong form
\begin{align}
r_h(\eta_h;\mu) = 0.
\label{eq:parametric_full}
\end{align}

Its hallmark feature is a greedy algorithm embedded in an offline-online decomposition procedure. 
Given a (sufficiently fine) training set $\mT\subset\mathcal{X}_\mP$. The \pzc{judicious exploration 
offline}
aims to build a low-dimensional surrogate for the $\mT$-induced solution manifold $\{\eta_h(\mu): \mu \in \mT\}$. This surrogate space is iteratively constructed via a hierarchical series of reduced basis $\{\eta_j\}_{j=1}^r$ by a greedy algorithm. At each iteration the snapshot corresponding to the most under-represented parameter value (were the current reduced space to be adopted), as identified by an error indicator or {\em a posteriori} error estimator, is added to the current set of bases. 
These snapshots are obtained through the full order solver \eqref{eq:parametric_full}. 
A defining feature of RBM is that the number of full model solves is minimum, i.e. the same as the surrogate space dimension. 
For parametric systems bearing a small Kolmogorov N-width \cite{pinkus_n-widths_1985}, the dimension of the surrogate space is orders of magnitude smaller than the total degrees of freedom for the full model in order to reach a high degree of accuracy. This difference in size leads to a dramatic decrease in computation time for the online simulations when a {\em reduced} solution is sought in the terminal surrogate space for each parameter value as the Galerkin or Petrov-Galerkin projection into the reduced space $\Gamma^r=\textrm{span}\{\eta_j,j=1,\dots,r\}$ constructed offline.

Leveraging the low-rank structure induced by the angular space for our problem, it is natural to treat the angular variable $\BOmega$ as our parameter. 
What prevents a direct application of RBM is the integral operator \eqref{eq:integral_operator} which poses a two-fold extra challenge. 
First, unlike the standard setting when RBM applies, the solutions for different parameter values $\{f(\cdot,\BOmega): \BOmega \in {\mathbb S}^{d-1}\}$  are all coupled through \eqref{eq:integral_operator} or its discrete counterpart $\lgl f \rgl \approx \sum_{j=1}^{N_\BOmega} \omega_j f(\cdot,\BOmega_j)$. 
Second, the robustness and the efficiency of the SASI iterative solver relies on the high quality of  the density approximation $\rho=\lgl f\rgl$.  In particular, the unstructured nature of the selected RB parameter samples {in the angular variable} prevents a robust and accurate numerical integration which usually requires a structured set of quadrature points. 

\subsection{The \fli{proposed} algorithm}
\label{sec:rb_alg}

We propose to effectively decouple the solutions for different angular samples by designing an iterative procedure where an approximation of the macroscopic density $\rho$ is constructed from the RB snapshots  and gradually refined as the RB space is built and RB solutions get more accurate. This iterative procedure manifests the first novel ingredient of our method in that a quantity  that is indirectly dependent on the parameter 
(i.e. the macroscopic density defined as an integral over the parameter domain) 
is fixed during one greedy iteration and updated only at the end of such iteration.  It allows the greedy algorithm to proceed efficiently which in turn enriches the surrogate solution space rendering the RBM solutions and the dependent macroscopic quantity more accurate. We emphasize that once the greedy algorithm converges, we consider the macroscopic density well-resolved and will adopt its terminal value online for any new parameter value. This is reasonable since the macroscopic density is not directly dependent on a particular parameter value. 
To resolve the lack of structure for a robust numerical integration, we develop a least squares density reconstruction strategy capable of integrating over an arbitrary set of selected angular samples.  In the remaining part of this subsection, we first describe our online solver with any given RB space which is repeatedly called offline to construct the terminal RB space from scratch. Next, the least squares reconstruction algorithm is presented and an $L^1$-based residual-free error indicator is reviewed. Finally, we finish by detailing our algorithm and making a few relevant remarks.

\subsubsection*{Online stage}\label{sec:online}
Assuming that 
the (discrete) reduced basis space corresponds to the column space of $\mU_{RB}$ and 
 the currently  reconstructed density is $\BRHO_{RB}$, 
the online solver amounts to seeking a Galerkin projection of $\BF_{\BOmega}$ into  $\mU_{RB}$ that satisfies the weak formulation for any given $\BOmega$. That is, we assume
\begin{align}
\label{eq:rbsol}
\BF_{\BOmega} \approx \mU_{RB}\BC_{RB}(\BOmega),
\end{align}
and compute $\BC_{RB}(\BOmega)$ by solving the reduced formulation of \eqref{eq:matrix-vector}, namely
\begin{align}
(\BMU_{RB})^T \BU_{\BOmega_i} \BMU_{RB} \BC_{RB}(\BOmega) + (\BMU_{RB})^T \BSIGMA_t \BMU_{RB} \BC_{RB}(\BOmega) = (\BMU_{RB})^T \BSIGMA_s\BRHO_{RB} + (\BMU_{RB})^T \Bg.
\label{eq:online_computation}
\end{align}

\begin{rem}
Note that the online solver \eqref{eq:online_computation} is repeatedly called during an iterative procedure offline to build up $\mU_{RB}$. By assuming that the approximation for the macroscopic density $\BRHO_{RB}$ is given and performing a reduced transport sweep (i.e. obtaining $\BC_{RB}(\BOmega)$ for all $\BOmega$) with a fixed known density, we effectively decouple the angular dependence of the system. Without this technique, all $f_\BOmega$ are coupled
\pzc{and this} leads to a larger system being inverted. The resulting increase in online (and thus offline) time is not amenable. An implication of this strategy is that we would need an {\em initial guess} for $\BRHO_{RB}$. 
We noticed that a poor initial guess may lead to inaccurate reduced approximations. In this paper, we sample at a  group of fixed but small number of quadrature points  $\mP_0=\{ \widehat{\BOmega}_k\}_{k=1}^{N_0}$ to perform a fully coupled full order solve to obtain this initial guess. 

\end{rem}

\subsubsection*{Least squares density reconstruction}\label{sec:reconstruct_rho}

The existence of a low rank structure for $\BF_\BOmega$ in its angular dependence,  {at least in the diffusive regime,} 
 indicates that the density can be accurately captured by the reduced space constructed offline. 
However, one main challenge in the offline stage is that the selected angular samples may not automatically give us a robust numerical quadrature formula {in the angular space}. Nevertheless, careful design of a density reconstruction algorithm approximating $\lgl\cdot\rgl$ is essential for the robustness and efficiency of the SASI solver. 
In this paper, we propose a reconstruction algorithm based on a least squares procedure. 
Indeed, we first fix a group of robust high order quadrature points $\{(\bar{\BOmega}_j, \bar{\omega}_j)\}_{j=1}^{\bar{N}_\BOmega}$.
Given an arbitrary group of angular samples $\mP_{RB}$ and corresponding distribution functions $\{f_\BOmega\}$, we   construct a least squares approximation of $f$, denoted by $f_{\rm ls}$ in the angular space.
We then calculate the macroscopic density $\rho$ based on $f_{\rm ls}$. 
\begin{align}
\rho = \lgl f \rgl \approx \sum_{j=1}^{\bar{N}_\BOmega} \bar\omega_j f_{\rm ls}(\bar{\BOmega}_j).
\label{eq:rho_reconstruct}
\end{align}
{A natural set of space to consider in 1D are polynomial function space. In 2D, to suit  the periodic structure, we use trigonometric function space.} The details of the $f_{\rm ls}$ reconstruction are as follows.  Here we denote the vector of degrees of freedom in space for $\BF_{\BOmega}$ as $\left(\BF_\BOmega(1),\dots,\BF_\BOmega(N_{dof})\right)^T$.
\begin{enumerate}
\item [{\bf 1D:}] With $\BOmega=v$, we find $\BF_{\rm ls}(v)=(p_1(v),\dots,p_{N_{dof}}(v))^T$ such that
\begin{align}
p_i(v) = \arg\min_{p\in P^s([-1,1])} \sum_{v\in\mP_{RB}}(p(v)-\BF_v(i) )^2,\; i=1,\dots, N_{dof}, 
\label{eq:1d_least_square}
\end{align}
where $P^s([-1,1])$ denotes the set of polynomials  on $[-1,1]$ of degree at most $s$. 
\item [{\bf 2D:}] With $\BOmega=(\cos(\theta),\sin(\theta))$, we find $\BF_{\rm ls}(\BOmega)=\BF_{\rm ls}(\theta)=(t_1(\theta),\dots,t_{N_{dof}}(\theta))^T$ such that,
\begin{align}
&t_i(\theta) = a_0 + \sum_{k=1}^s a_k \cos( k\theta) + \sum_{k=1}^{s-1} b_k\sin(k\theta),\;\text{ and }(a_0,\dots,a_s,b_1,\dots,b_{s-1})\;\text{solves} \notag\\
&{\displaystyle {{\arg\min}_{\wa_0,\dots,\wa_s,\wb_1,\dots,\wb_{s-1}}}} \sum_{\theta\in\mP_{RB}} \left( \wa_0 + \sum_{k=1}^s \wa_k \cos( k\theta) + \sum_{k=1}^{s-1} \wb_k\sin(k\theta)-\BF_{\theta}(i)\right)^2.
\label{eq:2d_least_square}
\end{align}
We use the trigonometric least squares approximation \eqref{eq:2d_least_square}, as $\theta=0$ and $\theta=2\pi$ represent the same point on the unit circle. {We note that the resulting $\BF_{\rm ls}$ is a $2\pi$-periodic function with respect to $\theta$, preserving the property of the original distribution function.}
\end{enumerate}
The choice of $s$ will be specified for our numerical experiments and discussed in Section \ref{sec:numerical}.

\subsubsection*{$L^1$ residual-free error indicator} 

A critical piece for the RB greedy algorithm is an {\em a posteriori} error estimator which guides the surrogate space construction and certifies the accuracy of the RB solution. It is often residual-based and can be derived by mimicking the {\em a posteriori} error analysis of the underlying full order scheme \cite{quarteroni2015reduced, hesthaven2016certified}, with the RB solution taking the place of finite element solution which plays the role of the exact solution. 
{\em A posteriori} error analysis of the streamline upwind finite element method for the kinetic equation \eqref{eq:kinetic} is considered in \cite{fuhrer1997posteriori}. With the need of solving its dual problem \cite{hartmann2003adaptive}, its extension to the RB setting is not computationally appealing. 
As a result, we turn to the highly efficient and provably reliable $L^1$ residual-free error indicator proposed in \cite{chen2019robust}. 

Indeed, using the notation from Section \ref{sec:rb_bg}, we assume the $r$-dimensional RB space is given by $\Gamma^r=\textrm{span}\{\eta_j,j=1,\dots,r\}$ where $\eta_j$ is the solution to \eqref{eq:parametric_pde} when the parameter takes value $\mu^j$. For a new parameter value $\mu$ whose corresponding RB solution is identified as
\[
\eta(\mu) = \sum_{j=1}^r c_j(\mu) \eta_j.
\]
The greedy choice informed by the $L^1$ error indicator then proceeds as follows
\begin{equation}
\mu^{r+1} \colonequals \arg\max_{\mu} \Delta_r(\mu), \mbox{ where } \Delta_r(\mu) = \sum_{j=1}^r \big | c_j(\mu) \big |.
\footnote{In actual computation, the basis $\{\eta_j\}$ will be processed e.g. subject to a Gram-Schmidt orthonormalization for numerical stability. The subsequent $L^1$ error indicator will have to incorporate the transformation matrix. 
}
\label{eq:L1greedy}
\end{equation}

As shown in \cite{chen2019robust}, the { $\{c_j(\mu)\}_{j=1}^r$} is the Lagrange interpolation basis in the parameter space. Taking the maximizer of  $\Delta_r(\mu) $ then amounts to controlling the growth of the Lebesgue constant. 
This lead to its effectiveness for selecting the RB snapshots \cite{chen2019robust} including for nonlinear steady-state or time-dependent problems \cite{ChenJiNarayanXu2021, ChenSigalJiMaday2021}. Finally, it is imperative to note that this indicator is straightforward to implement and essentially free to compute.

\subsubsection*{Greedy algorithm}
We are now ready to describe the greedy algorithm for iteratively constructing the reduced basis in the offline stage. It starts with a small set of quadrature points $\mP_0=\{ \widehat{\BOmega}_k\}_{k=1}^{N_0}$ with the initial reduced space constructed as the span of the resulting {snapshots} and the initial density $\BRHO^0$ computed accordingly. In each iteration, we calculate the projection of solution for an unselected parameter into the current reduced space and greedily expand the reduced space according to \eqref{eq:L1greedy}. 
With matrix-vector formulation, details of this algorithm is presented in Algorithm \ref{alg:greedy}. 
\begin{algorithm}[hbp]
\caption{RBM greedy algorithm for radiative transfer equation.\label{alg:greedy} }
\begin{algorithmic}[1]
\STATE{{\bf{Input:}} iteration number tolerance $M_{\rm tol}$, spectral ratio tolerance $r_{\rm tol}$, training set $\mP = \{ \BOmega_1,\dots, \BOmega_{N_\BOmega} \}$, and an initial quadrature rule $\mP_0=\{ \widehat{\BOmega}_k\}_{k=1}^{N_0}.$}
\STATE{Solve (globally coupled) \eqref{eq:matrix-vector} with a direct solver for $\BOmega \in \mP_0$ to obtain $\BF^0_{\widehat{\BOmega}_k}$ and $\BRHO^0$.}

\STATE {{\bf{Iteration:}}}
\STATE {Set $m = 0, r^m = 2 r_{\rm tol}$, $\BFF^{0}_{RB} = (\BF^0_{\widehat{\BOmega}_1},\dots, \BF^0_{\widehat{\BOmega}_{N_0} })$ and $\mP_{RB} =\{\widehat{\BOmega}_k,k=1,\dots, N_0\} $.}
\WHILE{$m\leq M_{\rm tol}$ and $r^m>r_{\rm tol}$}
		\STATE{Perform an SVD $\BFF_{RB}^{m}=\BMU^m_{RB}\BLambda_{RB}^m(\BMV^m_{RB})^T$.}
		\STATE{Calculate the spectral ratio $r^m = \frac{\lambda_{\min}^m}{tr(\BLambda^m_{RB})}$ where $\lambda_{\min}^m$ is the minimal singular value.}
	\FOR{$i=1:N_\BOmega$}	
			\IF{ $\BOmega_i\not\in \mP_{RB}$,} 			
			\STATE{\textbf{{Compute the RBM solution for $\BOmega_i$,}} $\BF^m_{\BOmega_i} = \BMU^m_{RB}\fli{ \BC^m_{RB}(\BOmega_i)}$ with \fli{$\BC^m_{RB}(\BOmega_i)$} solving 
			\begin{align*}
			(\BMU^m_{RB})^T \BU_{\BOmega_i} \BMU^m_{RB} \BC^m_{RB}(\BOmega_i) + (\BMU_{RB}^m)^T \BSIGMA_t \BMU^m_{RB} \BC^m_{RB}(\BOmega_i) = (\BMU^m_{RB})^T \BSIGMA_s\BRHO^m + (\BMU^m_{RB})^T \Bg.
			\end{align*}
			}
			\STATE{\textbf{{Calculate the \fli{$L^1$} error indicator \fli{$\mathcal{E}^m_{\BOmega_i} = \lVert \BMV_{RB}^n\left(\Lambda_{RB}^n\right)^{-1}\BC^m_{RB}(\BOmega_i) \rVert_{\ell^1}$.}
			}}}
			\ENDIF
			
	\ENDFOR
		\STATE{Set $ i_{\rm new} = \arg \max_i\{\mathcal{E}^m_{\BOmega_i}\}$, $m:=m+1$, and $\mP_{RB} = \mP_{RB}\bigcup \{\BOmega_{i_{\rm new}},\WBOmega_{i_{\rm new}}\}$.}
		\STATE{Execute the SASI method to solve \eqref{eq:matrix-vector}  and update $\BF_{\BOmega_{j_k}}^m$ for $\BOmega_{j_k}\in \mP_{RB}$ to assemble 
		\begin{align*}
		\BFF_{RB}^{m}=\left(\BF^{m}_{\BOmega_{j_1}},\dots,\BF^{m}_{\BOmega_{j_{N_m}}}\right), \; N_m = 2m+N_0\quad\text{and}\quad\mP_{RB}=\{\BOmega_{j_k}\}_{k=1}^{N_m}.
		\end{align*}}
		\STATE{Reconstruct $\BRHO^m$ based on $\BF_{\BOmega}^m$ for $\BOmega\in\mP_{RB}$ via \eqref{eq:rho_reconstruct} and \eqref{eq:1d_least_square} (or \eqref{eq:2d_least_square}).}
\ENDWHILE
\STATE{{\bf{Output:}}  the reduced basis $\mU_{RB}=\mU_{RB}^m$ and $\BRHO_{RB}=\BRHO^m$.}
\end{algorithmic}
\end{algorithm}

\begin{rem}
We emphasize two features of the greedy algorithm. The first feature is  the {\em symmetry-enhancing greedy addition}. 
 When we augment the RB space, in addition to the maximizer of the error indicator $\BOmega_{i_{new}}$ as determined by \eqref{eq:L1greedy}, we include its symmetric counterpart  $\WBOmega_{i_{new}}$.  Assuming that the training set is $\mP=\{\BOmega_j\}_{j=1}^{N_\BOmega}$. For each $\BOmega$, in 1D, we define $\WOmega=-\Omega$. In 2D, given $\BOmega=\left(\cos\theta,\sin\theta\right)$ with $\theta \in (0, 2 \pi)$, we define $\WBOmega=(\cos(\wtheta),\sin(\wtheta))$ with $\wtheta = {\rm mod}(\theta + \pi, 2 \pi)$. By construction, if $\BOmega \in \mP$ then $\WBOmega \in \mP$. In our numerical simulations, we observe that this method is {more robust} than adding  $\BOmega_{i_{new}}$ into $\mP_{RB}$ alone.  Indeed, for the 2D examples in  Section \ref{sec:numericaL^2d}, if one angular sample is added per iteration, the SASI iterative solver may fail to converge.

The second feature is a {\em spectral ratio stopping criteria}. The purpose is to mitigate the fact that our $L^1$ residual-free error indicator, albeit highly effective in identifying the next representative parameter value, is not an error estimator. Inspired by the POD
method, we monitor a spectral ratio as {an additional} stopping criteria. We define the spectral ratio  for the $m$-th iteration $r^m$ as
\begin{align}  
r^m=\frac{\lambda^m_{\min}}{Tr(\BLambda_{RB}^m)},
\label{eq:spectral_ratio}
\end{align}
where $\lambda^m_{\min}$ is the smallest diagonal element of $\BLambda_{RB}^m$, and $Tr(\cdot)$ is the trace operator.
\end{rem}

\section{Numerical results}\label{sec:numerical}
In this section, we present 
a series of one- and two-dimensional numerical examples to showcase the performance of the proposed RB method. {For the underlying DG spatial discretization, the discrete space in \eqref{DG:space} with $K=1$ is used.}  
{With the  consideration for the efficiency and robustness of the SASI method},  {we set the degree parameter $s$ in the density reconstruction 
as $s=m+1$ for 1D slab geometry and $s=\min(5,m+1)$ in 2D} during the $m$-th iteration of the greedy algorithm. 
Throughout the experiments, we measure the following absolute and relative $L^2$ errors. 
\begin{align*}
\mE_f = \max_{j} \lVert f_{\rm{F}}(\BOmega_j)-f_{\rm{R}}(\BOmega_j)\rVert , \quad \quad & R_f = \max_{j} \frac{\lVert  f_{\rm{F}}(\BOmega_j)-f_{\rm{R}}(\BOmega_j)\rVert }{\lVert f_{\rm{F}}(\BOmega_j)\rVert },\\
\mE_\rho = \lVert \rho_{\rm{F}}-\rho_{\rm{R}}\rVert , \quad \quad &R_\rho = \frac{\lVert \rho_{\rm{F}}-\rho_{\rm{R}}\rVert }{\lVert \rho_{\rm{F}}\rVert }.
\end{align*}
Here $f_{\rm{F}}$, $\rho_{\rm{F}}$ denote the full order numerical solutions. The RB solutions are $f_{\rm{R}}$, $\rho_{\rm{R}}$ and $\lVert\cdot\rVert$ is the standard $L^2$ norm of $L^2(\BX)$.

\pzc{In the current setting, there are no essential physical parameters (e.g. scattering cross section, boundary conditions) whose ``multi-query'' nature will make the (one-time) offline investment more worthwhile. However, we still compare the proposed RB method (including both the offline and online stages) against (a single query of) the full order DG solver using the training set for angular space discretization. We note that the online stage can be utilized to predict solutions at angular samples outside of the training set, a feat out of reach by the full order DG scheme. When the proposed method is utilized as a surrogate for the full order solve with respect to the training set, we concern both the offline and online efficiency. When its online stage is applied to predict solutions at \fli{``unseen'' angular samples,}
 we only take the online efficiency into account. We call the error associated with the training set ``training error" and the error associated with the test set ``testing error".} 
\fli{Here and below, the training set refers to the set of parameter values of the angular variable used during the offline  stage to build the surrogate space, while the test set refers to that used during the online stage to test the performance of the RB method. }

\subsection{One-dimensional examples }
\label{sec:numericaL^1d}
We perform the one-dimensional experiments
on a slab geometry domain $\BX=[x_L, x_R]$, which is discretized by a uniform mesh with $\Delta x=0.125$. The initial guess {for our algorithm} is obtained by $2$ Gauss-Legendre points, i.e. $N_0 = 2$. {We consider 5 examples and conduct two tests for each of them.} \pzc{The training set consists of 24 Gauss-Legendre points.}   We consider two different  
$r_{\rm tol}$ 
as 
$10^{-4}$ and $10^{-6}$. 
The algorithm stops once $r_{\rm tol}$  is reached leading to different RB dimensions for different examples. \fli{The testing errors associated with a test set will be reported.}
For the second test, we ask the algorithm to generate $M_{\rm tol}=12$ reduced bases for all problems and record the training error between the RB solutions and the full order  solutions. \fli{In the end, we demonstrate the robustness of the algorithm by varying the strength of the scattering cross section $\sigma_s$.} \pzc{These tests aim at  
showing the capability of our method to predict solutions at angular samples outside of 
the training set and the effectiveness of the $L^1$ residual-free error indicator.}

\bigskip
\noindent \textbf{Example 1 (scattering dominant):}
\begin{align*}
\BX=[0,10],\; G = 0.01, \; \sigma_t = 100,\;\sigma_s = 100,\; f(0,v) = 0\;\text{with } v>0, \;f(10,v) = 0 \;\text{with } v\leq 0.
\end{align*}
\textbf{Example 2 (spatially varying scattering coefficient):}
\begin{align*}
& \BX=[0,10],\;G = 0.01, \;\sigma_t = 100(1+x),\;\sigma_s = 100(1+x),\notag\\
&f(0,v) = 0\;\text{with } v>0, \;f(10,v) = 0 \;\text{with } v\leq 0.
\end{align*}
\textbf{Example 3 (two-material problem 1):}
\begin{align*}
&\BX=[0,20], \;
G = \begin{cases}
		5,\; 0<x<10,\\
		0,\;10<x<20,
	\end{cases}
\sigma_t = 100,\;
\sigma_s = \begin{cases}
		90,\; 0<x<10,\\
		100,\;10<x<20,
		\end{cases}\notag\\
&f(0,v) = 0\;\text{with } v>0, \;f(20,v) = 0 \;\text{with } v\leq 0. 
\end{align*}
\textbf{Example 4 (two-material problem 2):}
\begin{align*}
&\BX=[0,11] ,\;G = 0, \;
\sigma_t = \begin{cases}
		100,\; 1<x<11,\\
		2,\;0<x<1,
		\end{cases}\;
\sigma_s = \begin{cases}
		100,\; 1<x<11,\\
		0,\;0<x<1,
		\end{cases}\notag\\
&f(0,v) = 5\;\text{with } v>0, \;f(11,v) = 0 \;\text{with } v\leq 0.  
\end{align*}
\textbf{Example 5 (transport dominant):}
\begin{align*}
\fli{\BX=[0,10],} \;G = 0.01, \;\sigma_t = 1.2,\;\sigma_s = 1,\;f(0,v) = 0\;\text{with } v>0, \;f(10,v) = 0 \;\text{with } v\leq 0.
\end{align*}

\begin{figure}
\begin{center}
\subfigure[Example 1, $4$ reduced basis\label{fig:1d_test1}]{
 \includegraphics[width=0.31\textwidth]{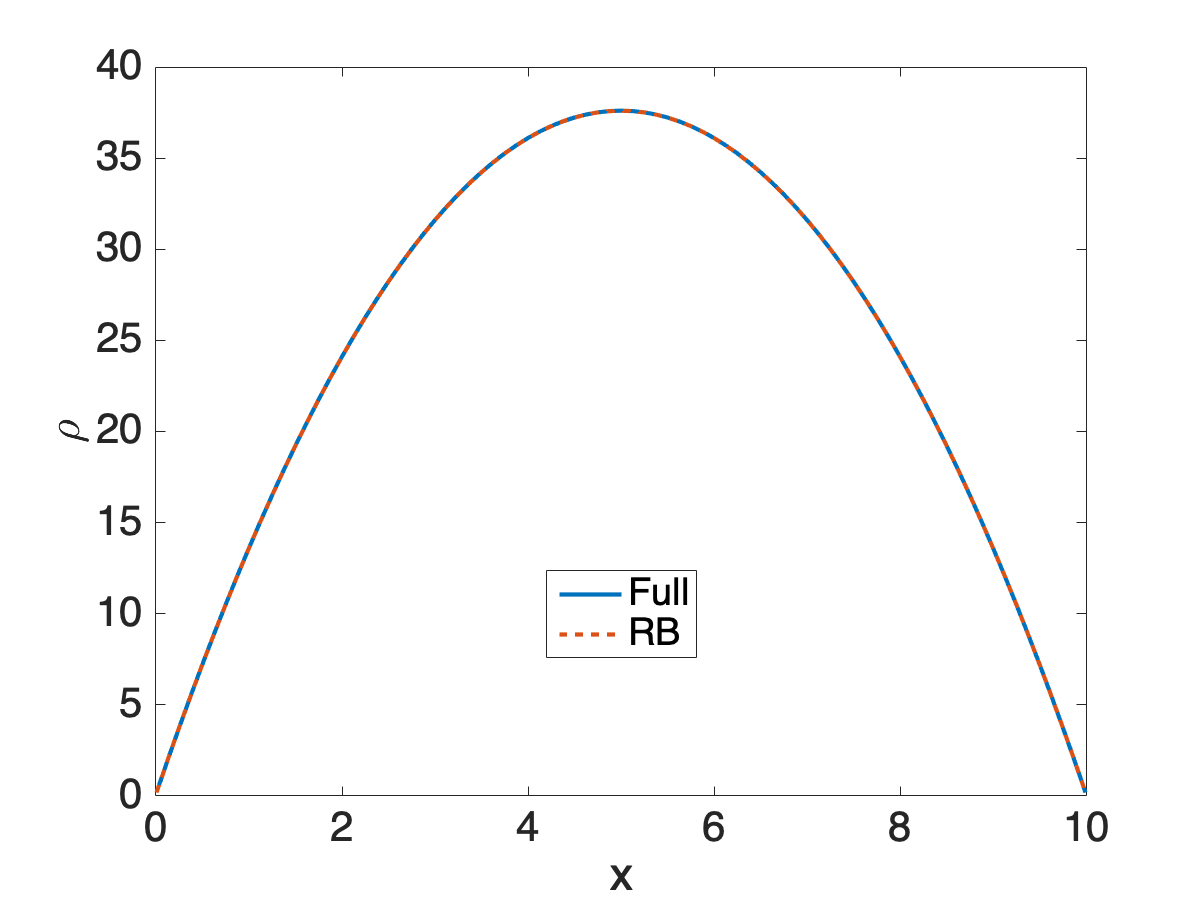}
}
\subfigure[Example 2, $4$ reduced basis\label{fig:1d_test2}]{
 \includegraphics[width=0.31\textwidth]{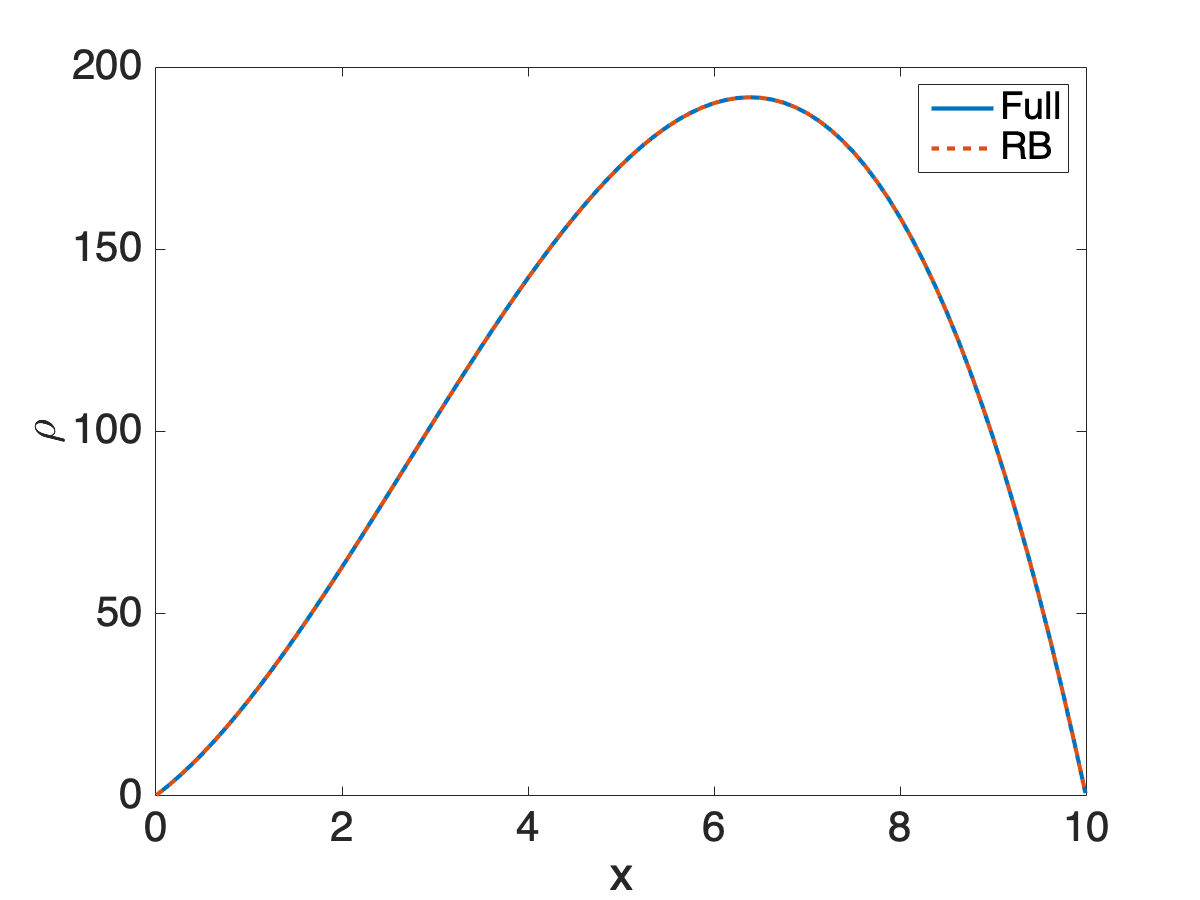}
}
\subfigure[Example 3, $4$ reduced basis\label{fig:1d_test5}]{
 \includegraphics[width=0.31\textwidth]{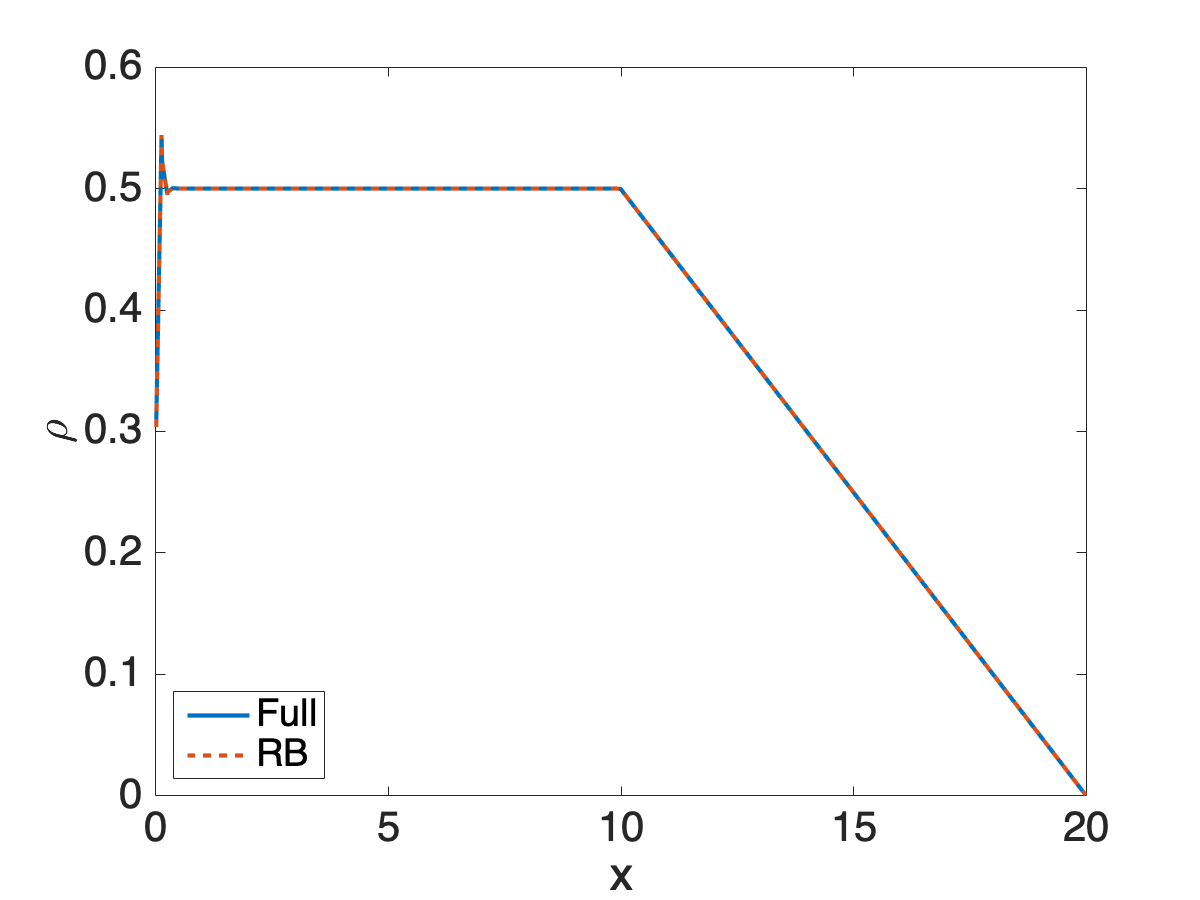}
}\\
\subfigure[Example 4, $8$ reduced basis\label{fig:1d_test4}]{
 \includegraphics[width=0.31\textwidth]{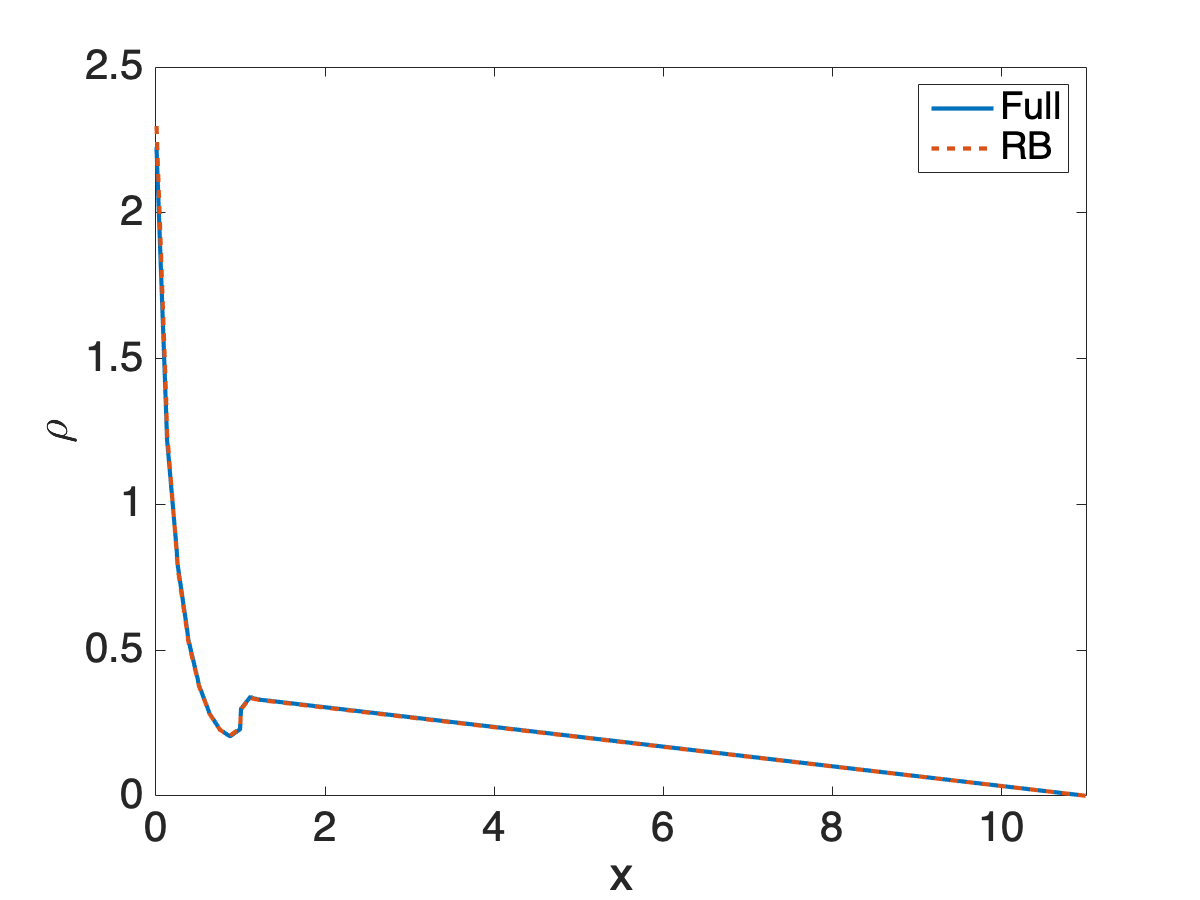}
}
\subfigure[Example 5, $10$ reduced basis\label{fig:1d_test5}]{
 \includegraphics[width=0.31\textwidth]{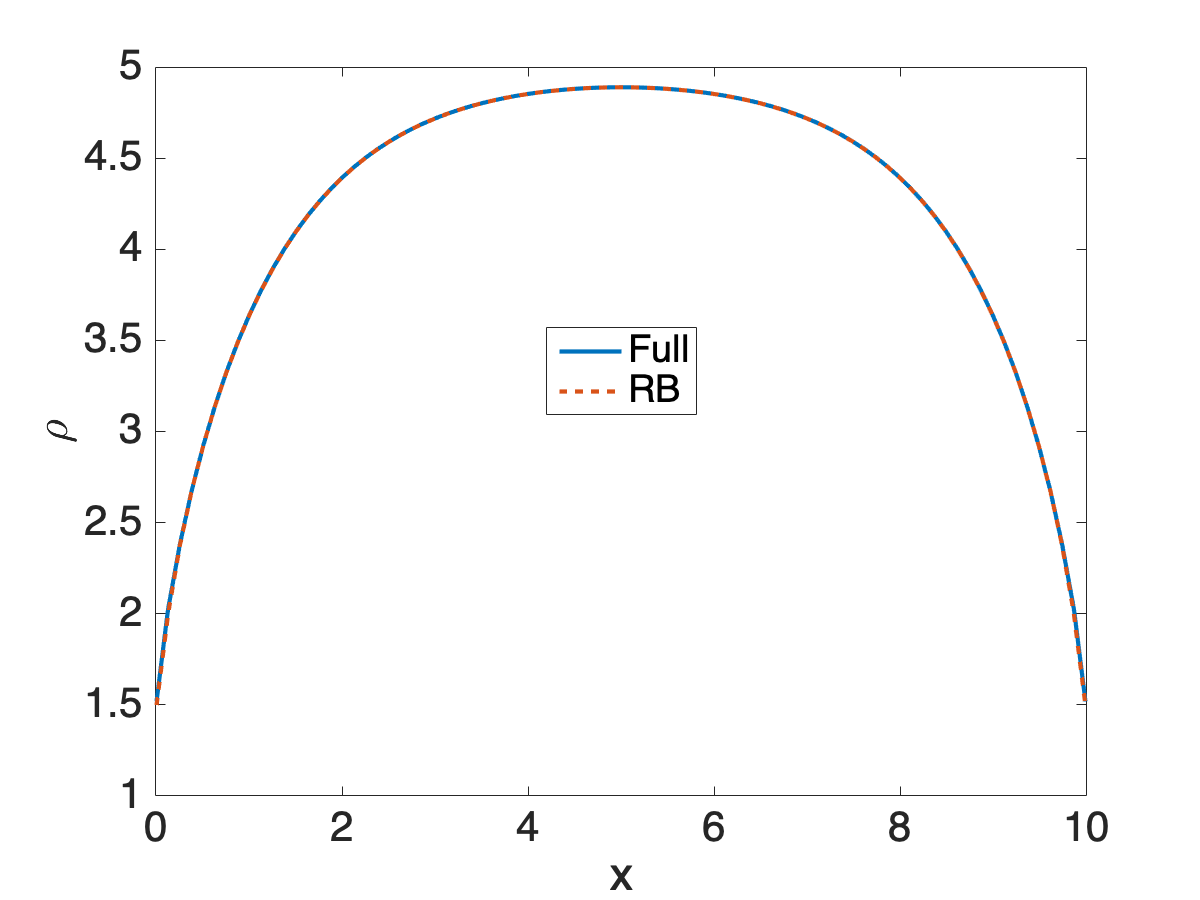}
}
\end{center}
\caption{The true (i.e. full order) and RB densities for the one-dimensional test examples 1-5 of 1D slab geometry. The number of reduced basis is determined by $r_{\rm tol}=10^{-4}$.\label{fig:1d_test_rho1}}
\end{figure}

\medskip
\noindent {\bf \pzc{Online prediction accuracy:}} \pzc{We first present the prediction accuracy for a test set with $32$ quadrature points} in the angular variable, \fli{that are different from the angular samples in the training set}. With $r_{\rm tol}=10^{-4}$, the comparison between  the density obtained by the RB solutions and that by the full order solutions is presented in Figure \ref{fig:1d_test_rho1}. 
For all five examples, the RB solutions match the full order solutions well. We present \fli{the testing errors in $f$ and training errors in $\rho$} 
for these examples in Tables \ref{tab:1d_test1_1e-4} and \ref{tab:1d_test1_1e-6}. One can see that the RB method is most effective for scattering dominant problems (i.e. Examples 1-3). It {achieves} 
very good \fli{degree} of accuracy (3 to 4 digits) for both $f$ and $\rho$ with just 4 to 6 reduced basis functions. The solution of the two-material problem (Example 4), {that involves a purely absorbing subregion without scattering}, has more complicated structure. 
Consequently, more basis functions are needed, and the accuracy  of $f$ is worse than the other examples.  Example 5 is transport dominant leading to distribution functions being far {from the macroscopic density. 
} As a result, it demands the most number of RB functions.  It is worth noting that, for all five examples, the method commits less than $1\%$ relative $L^2$ error for $\rho$. To reach this accuracy, the full model needs $160$ {degrees of freedom} for each angular sample, while our RBM  uses at most $14$ global reduced basis \fli{functions. } 
\begin{table}[htbp]
  \centering
  \medskip
    \begin{tabular}{|l|c|c|c|c|c|c|}
    \hline
 		&RB dimension &  	$\mE_f$&   $R_f$ & $\mE_\rho$  & $R_\rho$\\ \hline
 Example 1&	4	&	8.04e-3&  9.26e-3\% & 8.00e-3 & 9.21e-3\%\\ \hline	
 Example 2&	4	&	8.51e-3&  2.05e-3\% & 8.49e-3 & 2.04e-3\%\\ \hline	
 Example 3&	4	&	9.35e-4&  5.13e-2\% & 8.25e-4 & 4.53e-2\%\\ \hline	
 Example 4&	8	&	9.99e-2&  1.28e+1\% & 2.59e-3 &2.54e-1\%\\ \hline	
 Example 5&	10	&	7.50e-2&  5.39e-1\% & 4.37e-3 & 3.21e-2\%\\ \hline	
 \end{tabular}
   \caption{\pzc{Testing error for $f$ and training error for $\rho$} with $r_{\rm tol}=10^{-4}$, 1D slab geometry. \label{tab:1d_test1_1e-4}  }
  \centering
  \medskip
    \begin{tabular}{|l|c|c|c|c|c|c|}
    \hline
 		&RB dimension &  	$\mE_f$&   $R_f$ & $\mE_\rho$  & $R_\rho$\\ \hline
 Example 1&	6	&	3.34e-3&  3.84e-3\% & 3.31e-3 & 3.80e-3\%\\ \hline	
 Example 2&	4	&	8.51e-3&  2.05e-3\% & 8.49e-3 & 2.05e-3\%\\ \hline	
 Example 3&	6	&	4.31e-4&  4.83e-2\% & 3.77e-4 & 2.07e-2\%\\ \hline	
 Example 4&	10	&	4.73e-2&  6.09e-0\% & 1.91e-3 & 1.87e-1\%\\ \hline	
 Example 5&	14	&	2.78e-2&  2.00e-1\% & 1.93e-3 & 1.42e-2\%\\ \hline	
 \end{tabular}
      \caption{\pzc{Testing error for $f$ and training error for $\rho$} with $r_{\rm tol}=10^{-6}$, 1D slab geometry.\label{tab:1d_test1_1e-6} }
\end{table}

\begin{figure}[htbp]
\subfigure[Error history of $f$\label{fig:1d_test_error_history}]{
 \includegraphics[width=0.5\textwidth]{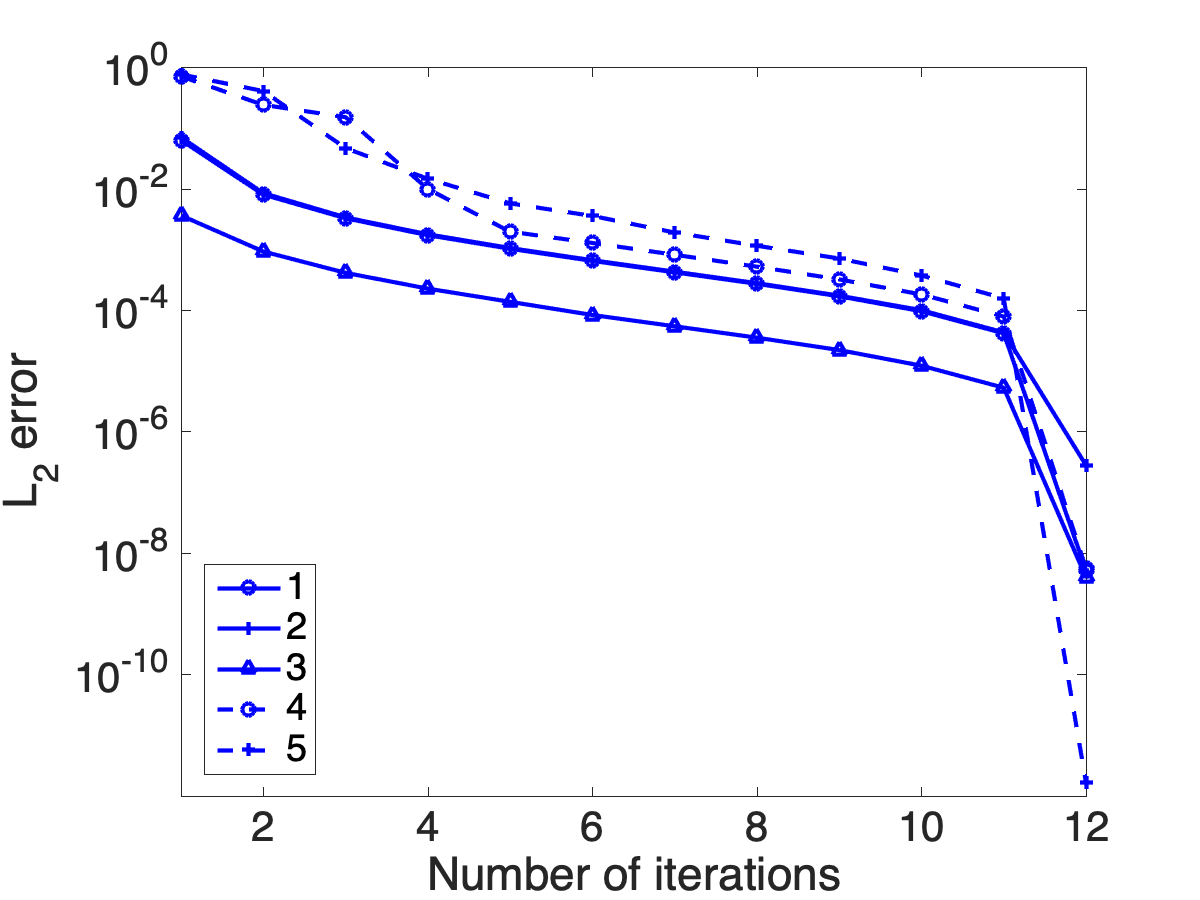}
 }
 \subfigure[History of spectral ratio 
 \label{fig:1d_test_svd_history}]{
 \includegraphics[width=0.5\textwidth]{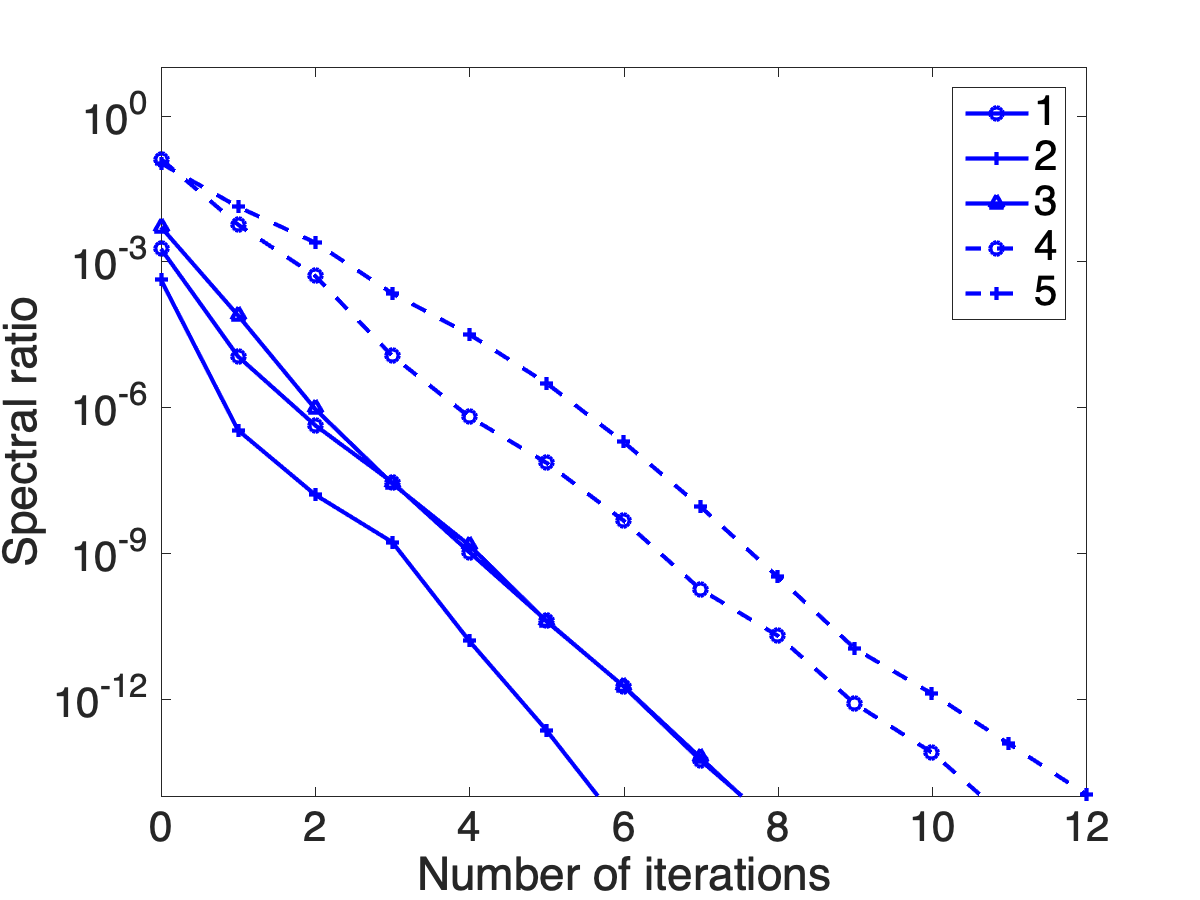}
 }
 \caption{Histories of convergence 
 of the \fli{$L^2$ error  of $f$ } and the monitored spectral ratio during training for 1D Examples 1-5.\label{fig:1d_test_history}}
\end{figure}

\medskip
\noindent  {\bf \pzc{Effectiveness of the error indicator}:} In Figure \ref{fig:1d_test_history}, we present the histories of the convergence of the $L^2$ error of  $f$ and the spectral ratio $r^m$ \fli{as the number of iteration grows}. We observe that the error decays as the dimension of the reduced space increases with a terminal error smaller than $10^{-6}$ for the worst case. We also note that the spectral ratio $r^m$ decays exponentially, albeit with different speed for different examples. The pattern, across different problems, of the spectral ratio decay leads to the different terminal RB dimensions for different problems as shown in Tables \ref{tab:1d_test1_1e-4}  and \ref{tab:1d_test1_1e-6}. 

\begin{figure}
\subfigure[Dimension of reduced space for $r_{\rm tol}=10^{-8}$]{
 \includegraphics[width=0.5\textwidth]{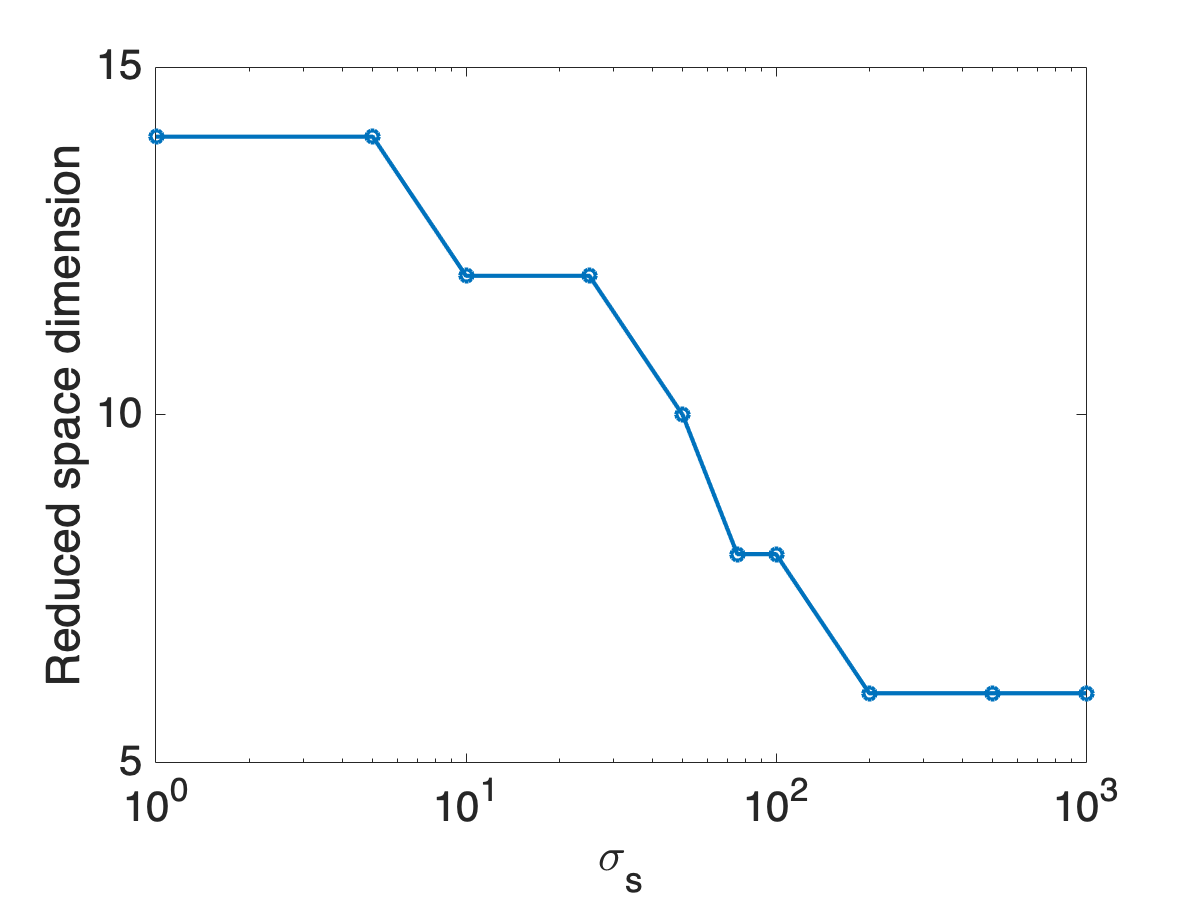}
 }
 \subfigure[\fli{Relative testing error}  of $f$]{
 \includegraphics[width=0.5\textwidth]{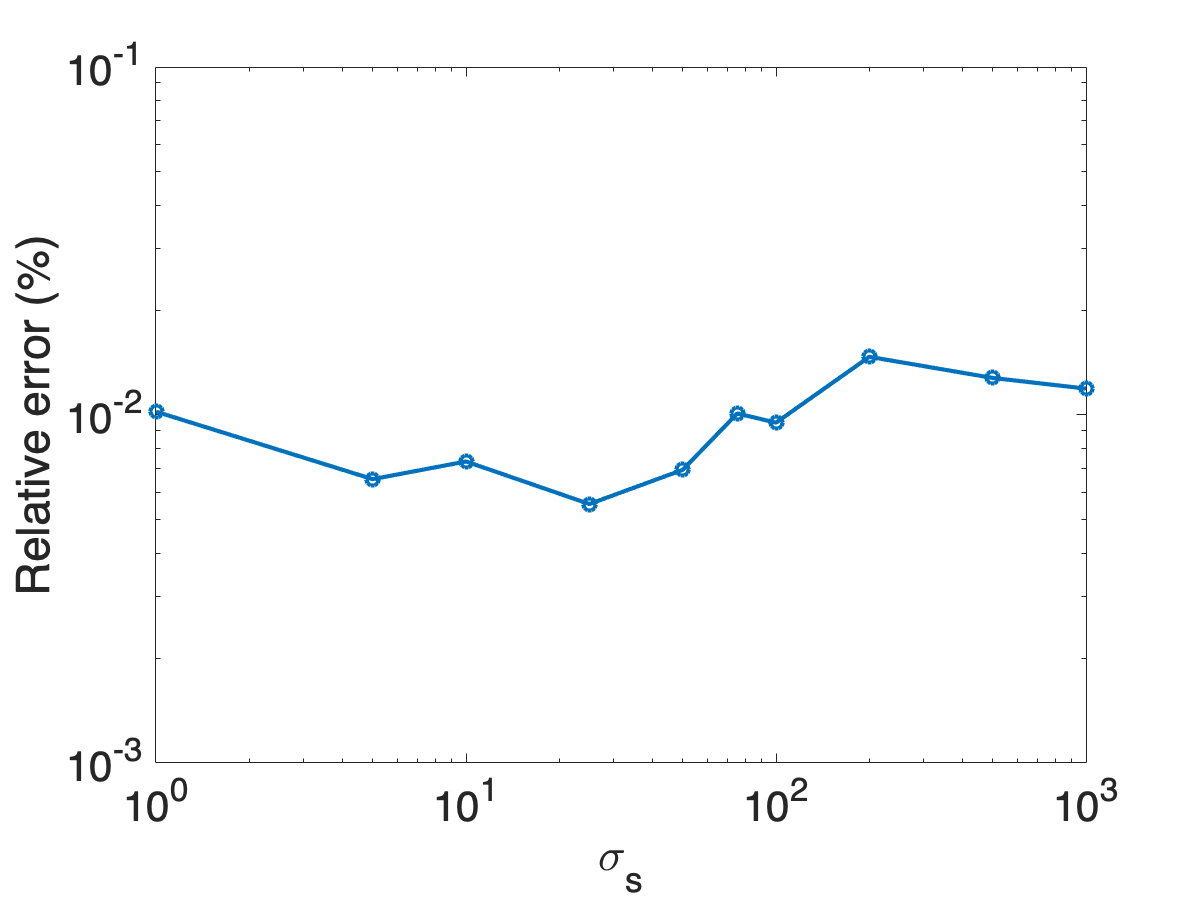}
 }
 \caption{RBM robustness test: Dimension of \fli{the reduced space 
 and the relative testing error} versus different scattering cross section. \label{fig:1d_sigma_s_vs_dim}}
\end{figure}

\medskip
\noindent {\bf Robustness test {with respect to $\sigma_s$:}} 
To showcase the robustness of our method, we consider the following example {by} 
varying strength of the scattering cross section $\sigma_s$,
\begin{align*}
\fli{\BX=[0,10],\;} G = 0.01, \;\sigma_t \equiv C+0.5,\;\sigma_s \equiv C,\;f(0,v) = 0\;\text{with } v>0, \;f(10,v) = 0 \;\text{with } v\leq 0,
\end{align*}
where $C=1,5,10,25,50,75,100,200,500,1000$. We use $40$ {Gauss-Legendre}
points as the training set, $32$ 
{Gauss-Legendre}
points as the test set, and take $r_{\rm tol}=10^{-8}$. The dimensions of the resulting RB spaces and the corresponding relative \pzc{testing} errors are reported in Figure \ref{fig:1d_sigma_s_vs_dim}. We clearly observe that the more scattering dominant the problem is, the fewer reduced basis \fli{functions are }
needed. The fact that the relative {testing} error 
 is on the same level for different $\sigma_s$, even though we only monitor the spectral ratio, attests to the robustness of our method and the reliability of our error indicators.

\subsection{Two-dimensional examples}\label{sec:numericaL^2d}
In this section, we consider the following four examples, all with 
the source term 
$G(x,y) = \exp\left( -100( (x-5)^2+(y-5)^2)\right)$ and the zero inflow boundary condition on the computational domain $\BX=[0,10]^2$. 

\begin{figure}
\centering
{
\includegraphics[width=0.5\textwidth]{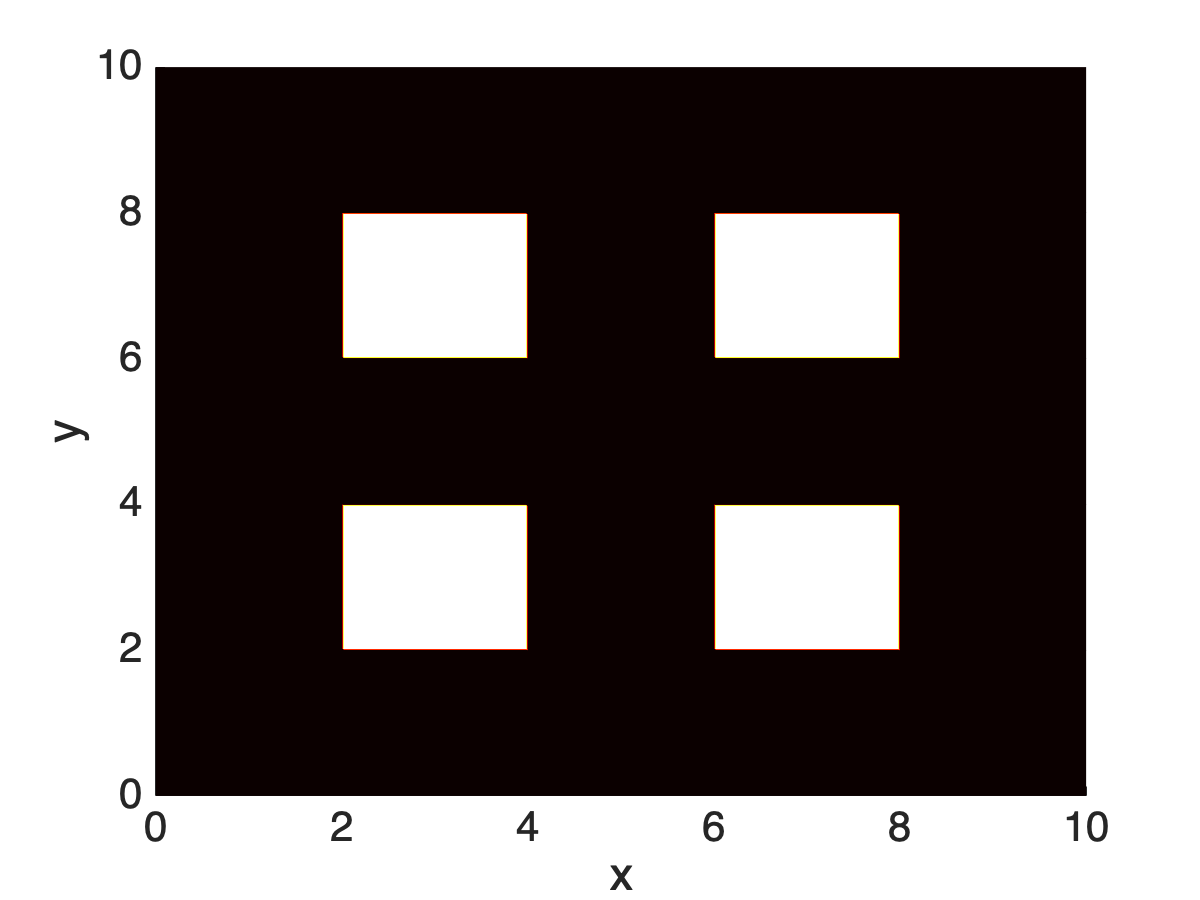}
 }
 \caption{Geometry of the checkerboard (2D Example 1). Black region: scattering dominant;  white region: transport dominant}
  \label{fig:checkerboard_geometry}
\end{figure}

\medskip
\noindent \textbf{Example 1 (checkerboard):} This is a multiscale problem, with part of the domain being scattering dominant and the rest being transport dominant. The checkerboard geometry is shown in Figure \ref{fig:checkerboard_geometry}. The white region is defined as
$\bigcup_{i,j=1}^{2}\{\max( |x-x_i|,|y-y_j|)<1\},\quad\text{with}\quad x_1=y_1=3,\;x_2=y_2=7, $
where we set $\sigma_s(x,y)=1$ and $\sigma_t(x,y)=\sigma_s(x,y)+1$. In the black region, we have $\sigma_s(x,y)=\sigma_t(x,y)=100$.

\smallskip
\noindent\textbf{Example 2 (scattering dominant):} $\sigma_s(x,y)=\sigma_t(x,y) = 100.$

\smallskip
\noindent\textbf{Example 3 (intermediate regime):} $\sigma_s(x,y)=\sigma_t(x,y) = 10.$

\smallskip
\noindent\textbf{Example 4 (transport dominant):} $\sigma_s(x,y)=\sigma_t(x,y) = 1.$
 
 \medskip
In our experiments, a uniform rectangular mesh of $40 \times 40$ is used for the DG scheme. 
To generate the initial guess, 
we start from $\{\theta_j\}_{j=1}^{N_0}$, with $\theta_j=\frac{2(j-1)\pi}{N_0}$, and use $N_0=4$ for Examples 1-3 and $N_0=8$ for Example 4.

\fli{We start with applying the proposed RB algorithm  with a training set of $N_\BOmega=32$ quadrature points in the angular variable, and obtain a surrogate model for the full order upwind DG solver.  To show the efficiency and accuracy of this surrogate solver, we first present the relative computational time and the accuracy with respect to a full order solve based on the same $N_\BOmega=32$  quadrature points, and then use the RB method to predict $f$ at angular samples that are not included in the training set.}

\begin{figure}
\centering
\subfigure[\pzc{Relative total computational time for  Examples 1-3} ]{
 \includegraphics[width=0.45\textwidth]{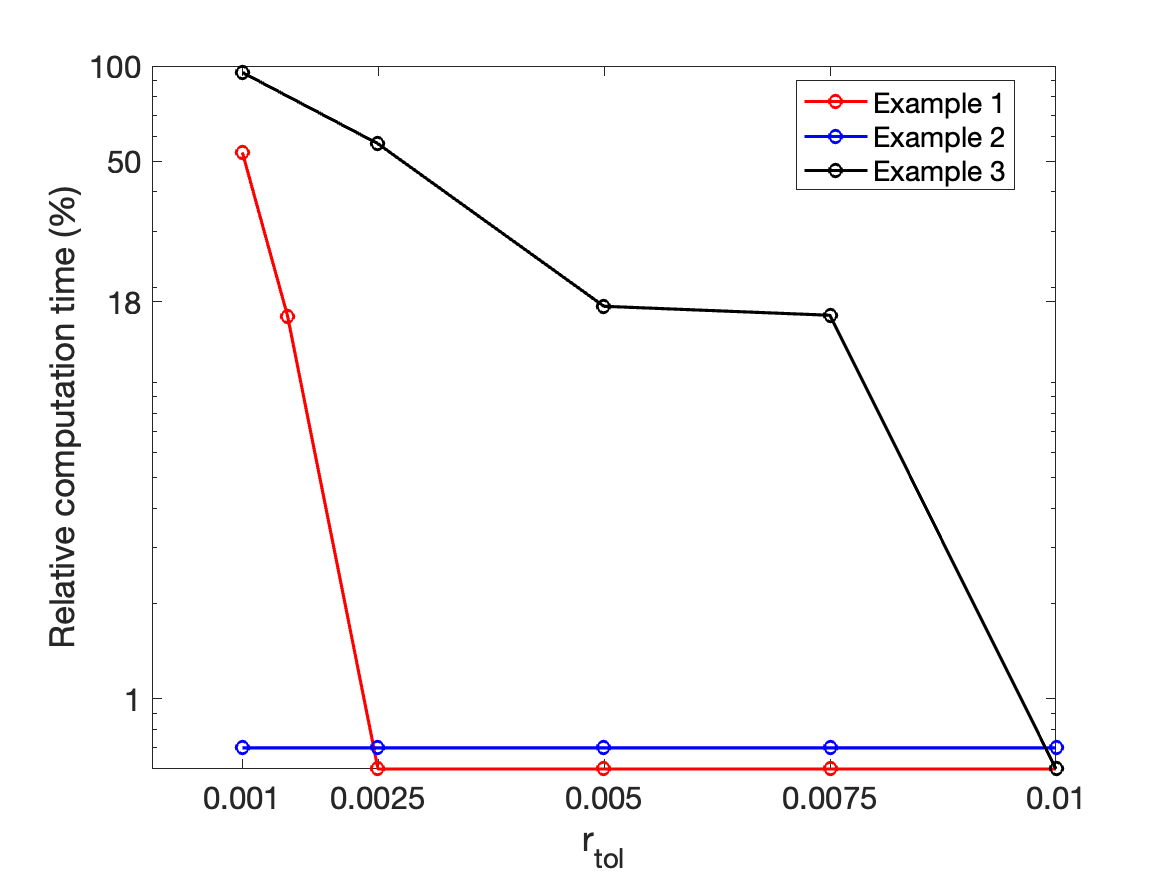}
}
\subfigure[\pzc{Relative total computational time for Example 4}]{
 \includegraphics[width=0.45\textwidth]{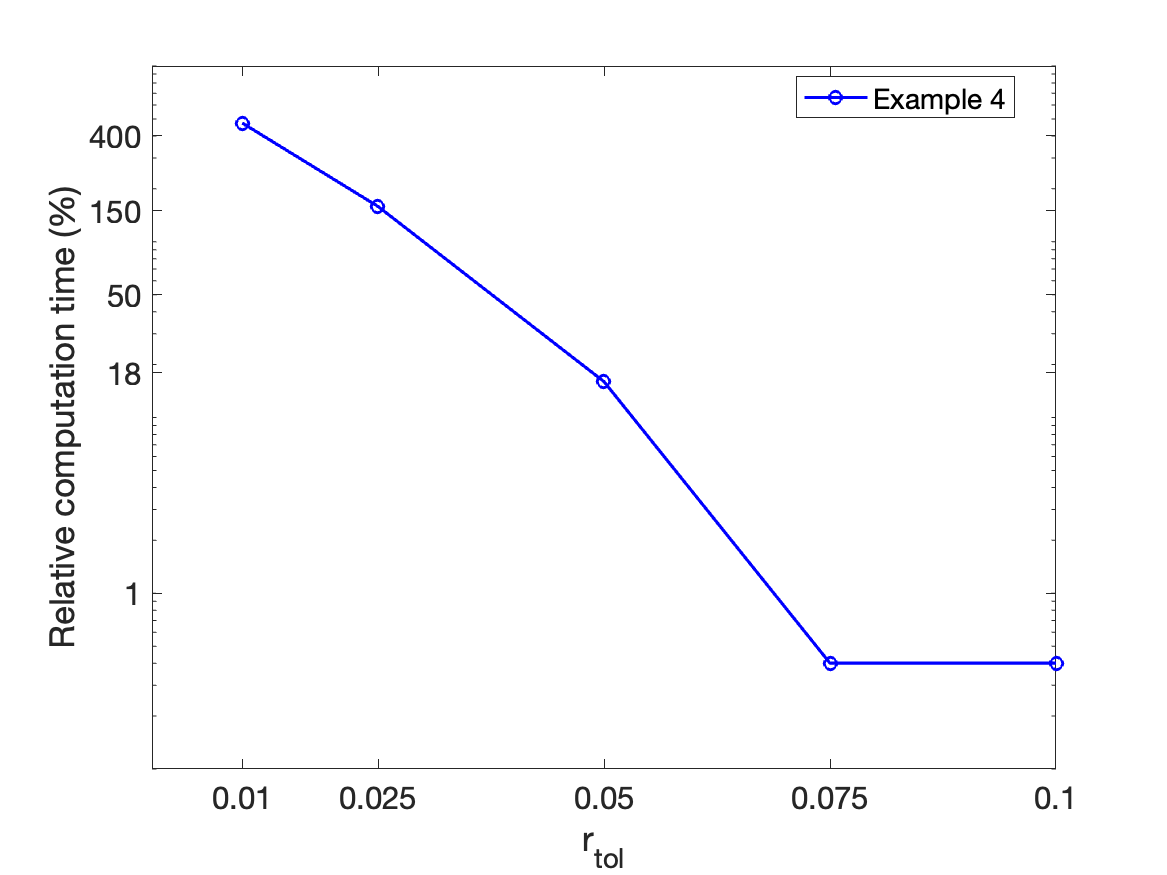}
}
\subfigure[\pzc{Relative errors for  Examples 1-3}]{
 \includegraphics[width=0.45\textwidth]{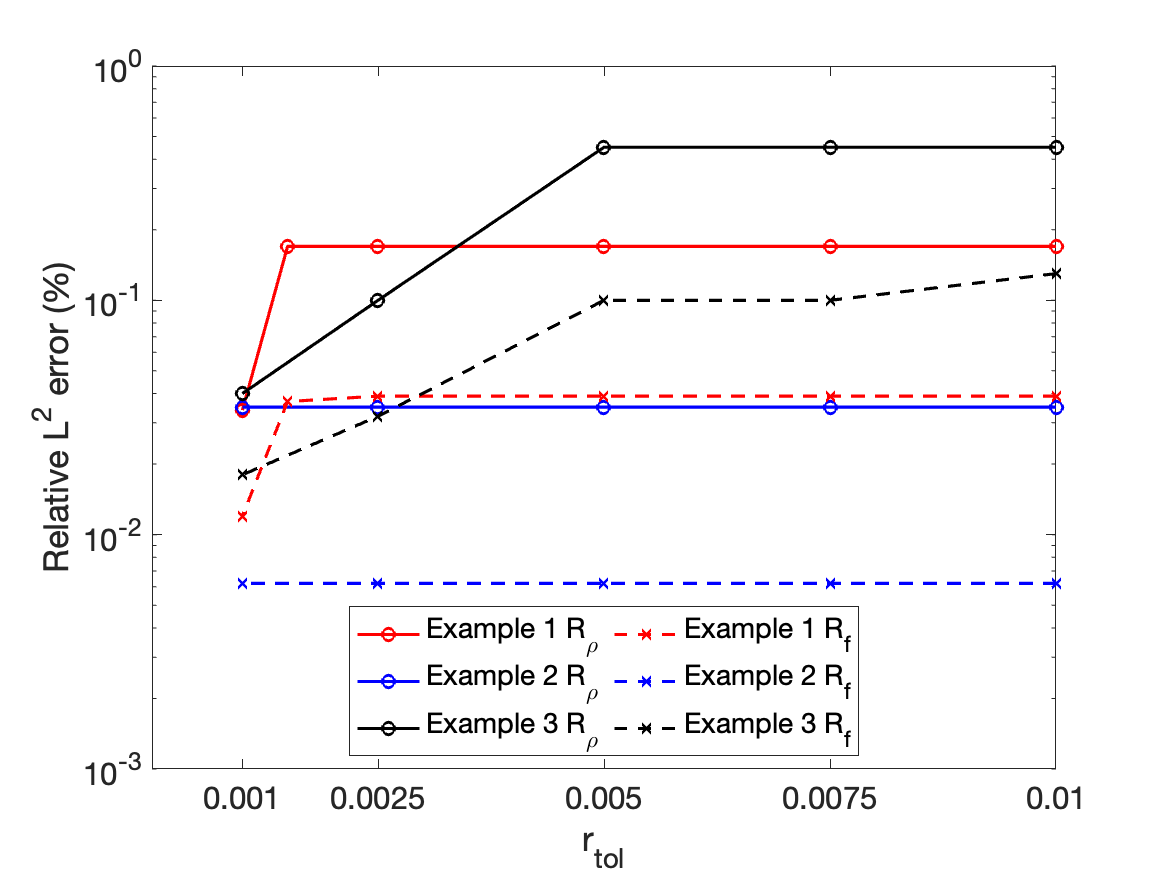}
}
\subfigure[\pzc{Relative errors for Example 4}]{
 \includegraphics[width=0.45\textwidth]{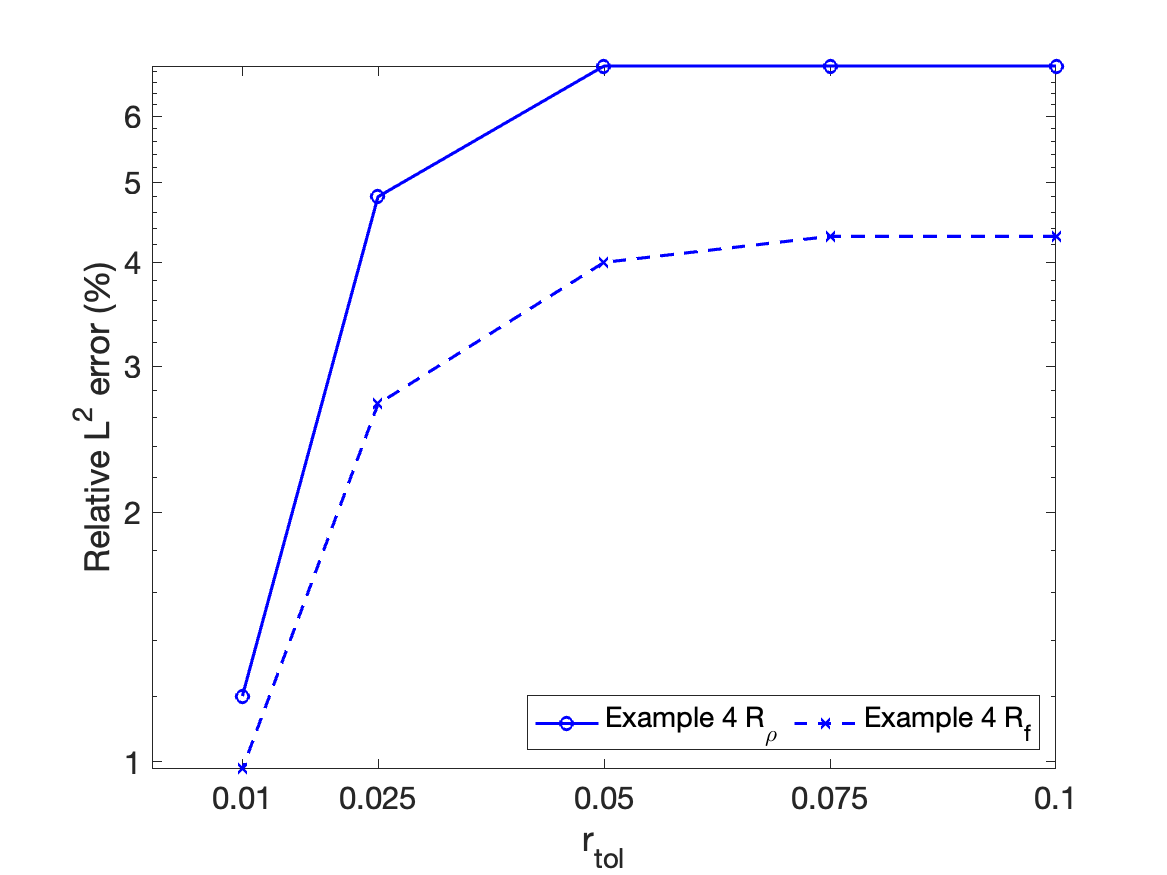}
}
 \caption{\pzc{Relative total computational time and relative $L^2$ training errors  of 2D examples with $N_\BOmega=32$ and different $r_\textrm{tol}$.}\label{fig:error_and_time_vs_tol}}
\end{figure}

\medskip
\noindent {\bf \pzc{Efficiency and the training accuracy as a surrogate solver:}}
\fli{In Figure \ref{fig:error_and_time_vs_tol}, the relative total computational time (including both the online and offline time) with respect to the full order solver based on $N_\BOmega=32$ quadrature points as well as  the relative $L^2$ error  
are presented as a function of $r_{\rm tol}$. As expected, the error decays as $r_{\rm tol}$ becomes smaller. Examples 1-3 are in the diffusive regime, intermediate regime and as a multiscale problem with large scattering dominant subregions, respectively. For these problems,  with a suitable $r_{\rm tol}$, one can achieve less than $1\%$ error with less than $18\%$ relative computational time.  Example 4 is transport dominant. To achieve less than $2\%$ error for this example, we need to use $r_{\rm tol}$$=10^{-2}$ and the offline computational time is longer than that of a full order solve. To understand the poor efficiency observed  for Example 4, one can recall that at the end of each greedy iteration in the offline stage, a full order method will be applied based on the angular sample set $\mathcal{P}_{RB}$ selected so far. Though each problem is of small size,  the transport-dominant nature of the model can require relatively more greedy iterations hence more solves of such small problems, due to the known  slow decay of the Kolmogorov $N$-width for transport dominant problems \cite{greif2019decay,ohlberger2015reduced}.  To improve the offline efficiency \fli{for such problems}, nonlinear reduced order models may be needed.}

\pzc{In Table \ref{tab:2d_test1_1e-3}, the training errors and the dimensions of the reduced order model with $N_\BOmega=32$ are presented. \fli{For Examples 1-3 with $r_{\textrm{tol}}=10^{-3}$, we achieve less than $0.1\%$ relative $L^2$ errors, and for Example 4 with  $r_{\textrm{tol}}=10^{-2}$, we achieve less than $2\%$ relative $L^2$ errors.}
 The scattering dominant problem, the multiscale checkerboard problem, and the intermediate regime problem all  need relatively  \fli{smaller RB spaces.} The transport dominant example requires the most \fli{reduced basis functions,}  
 and has relatively lower accuracy.  Nevertheless, the full order solve takes $6400$ degrees of freedom for each sample of the angular variable. In comparison, our RB algorithm only needs  fewer than $0.5\%$ degrees of freedom online. } 

\begin{figure}
\centering
\subfigure[\pzc{Relative \fli{testing errors}} \label{fig:2d_relative_error}]{
 \includegraphics[width=0.45\textwidth]{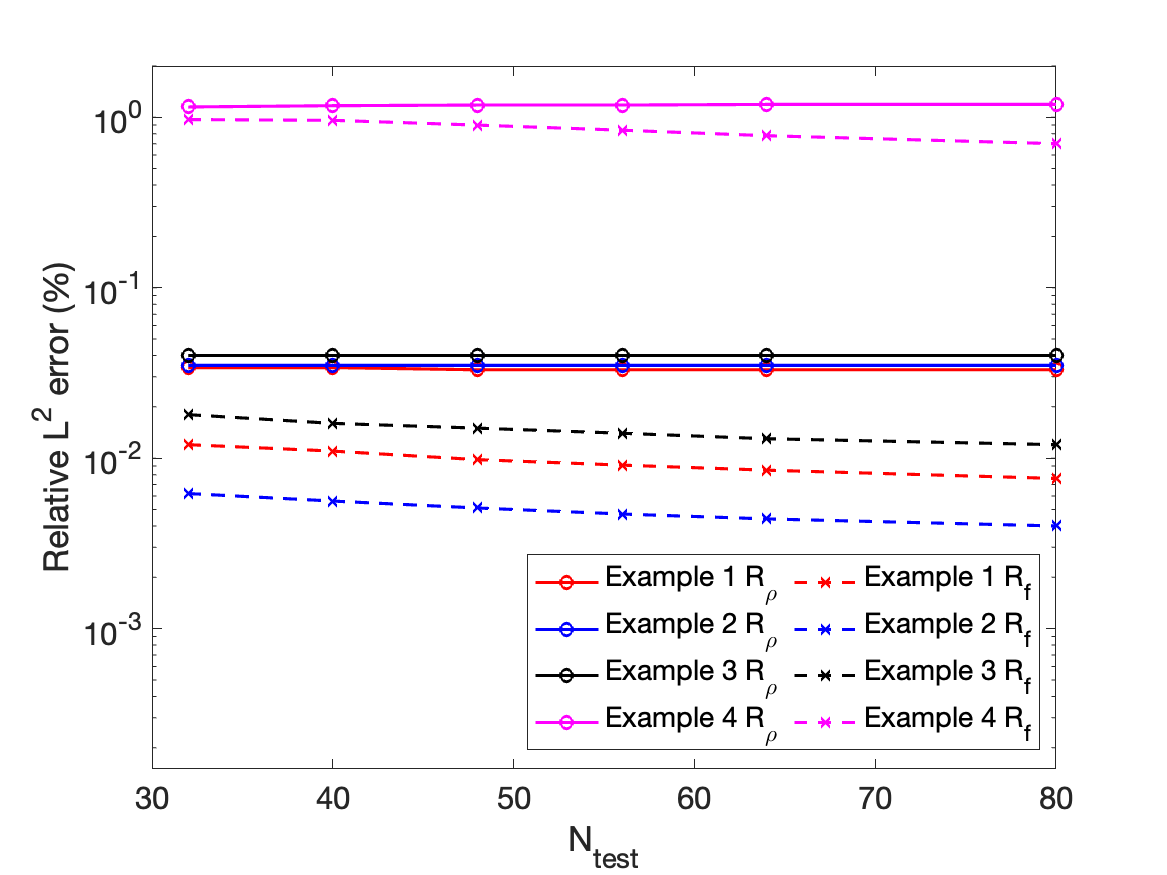}
}
\subfigure[\pzc{Relative online computational time} \label{fig:2d_relative_time}]{
 \includegraphics[width=0.45\textwidth]{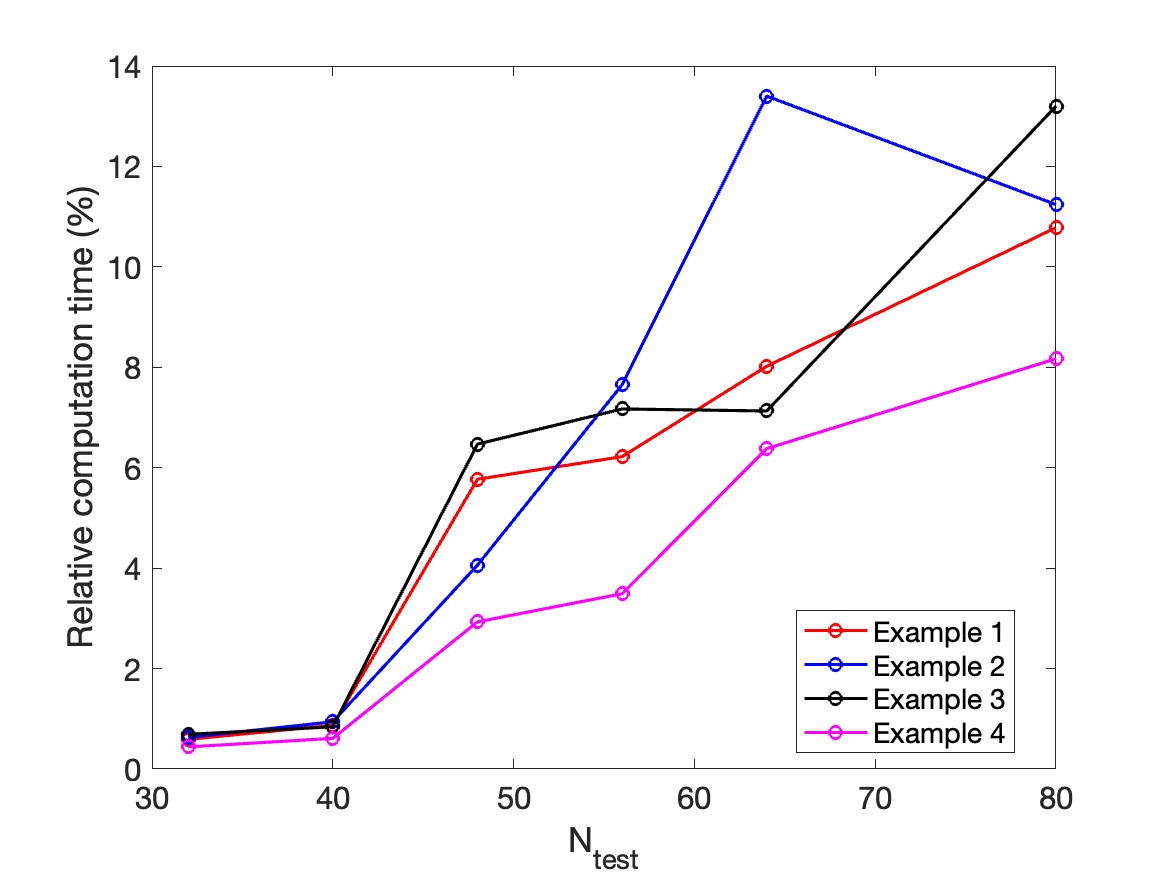}
}
  \caption{\fli{Relative testing errors and  relative online computational time for different test sets of size $N_\textrm{test}$, with 
 $r_{\rm tol}=10^{-3}$ for 2D Examples 1-3 and $r_{\rm tol}=10^{-2}$ for 2D Example 4.}
 \label{fig:2d_relative_time:a}}
\end{figure}

\medskip
\noindent {\bf \pzc{Online efficiency and the prediction accuracy:}} 
\fli{With the reduced basis functions obtained from a training set of $N_\BOmega=32$ quadrature points, we predict the solution $f$ at other angular samples. In Figure \ref{fig:2d_relative_error}, we present the relative testing errors for a group of  test sets. As the size $N_{test}$ of test sets varies,  the relative errors stay at almost the same level as that for the original training set. To illustrate the computational efficiency of online prediction, we further report in  Figure \ref{fig:2d_relative_time}  the relative online computational time with respect to that of the full order solve (with the same angular quadrature points as the test set), and they are always below $14\%$ as the size of test sizes varies.  This implies that the proposed RB algorithm can be used as a building block to construct ROMs when the model has essential parameters (e.g. scattering or absorption cross sections, boundary data etc), and substantial saving can be expected for online computation.}

\begin{figure}[htbp]
\subfigure[Example 1, $8$ reduced basis\label{fig:2d_test2}]{
 \includegraphics[width=0.33\textwidth]{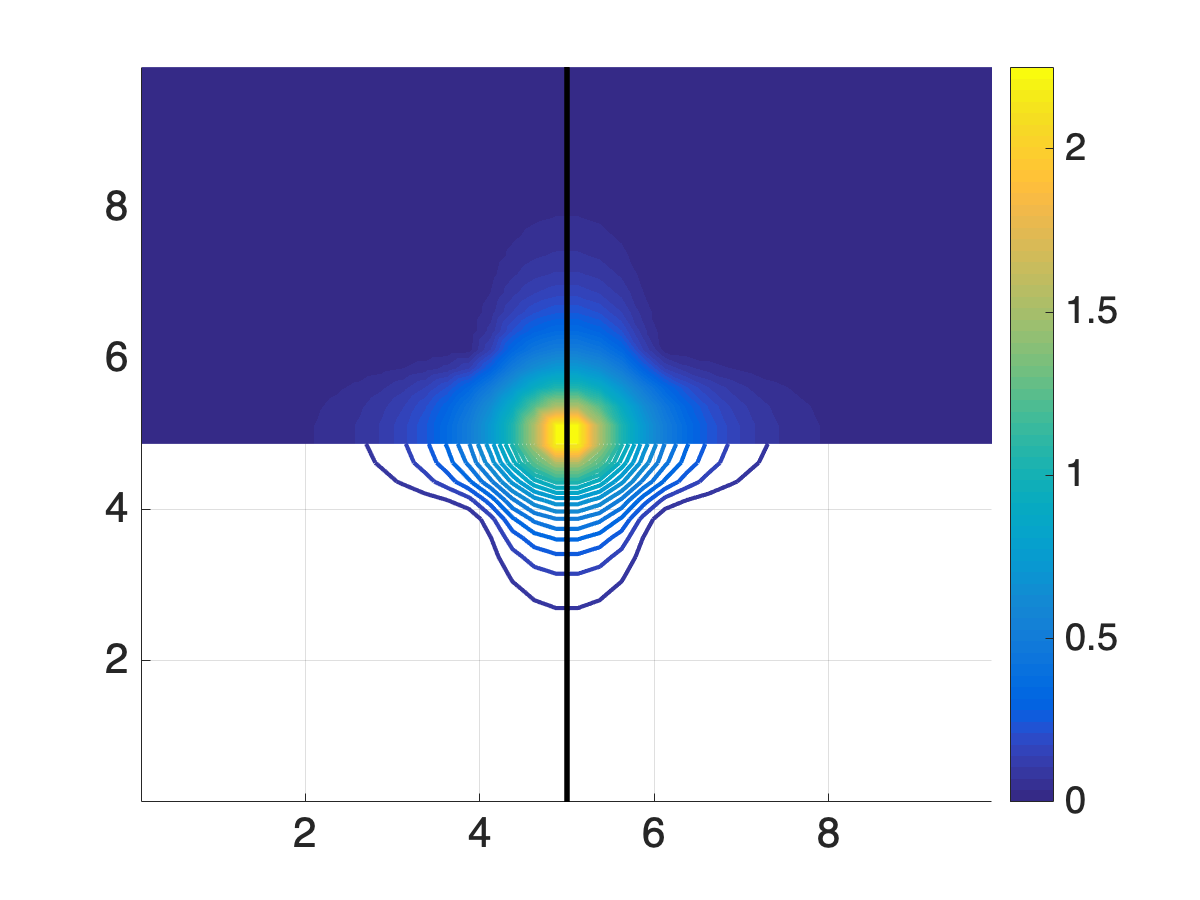}
}
\subfigure[Example 2, $4$ reduced basis\label{fig:2d_test1}]{
 \includegraphics[width=0.33\textwidth]{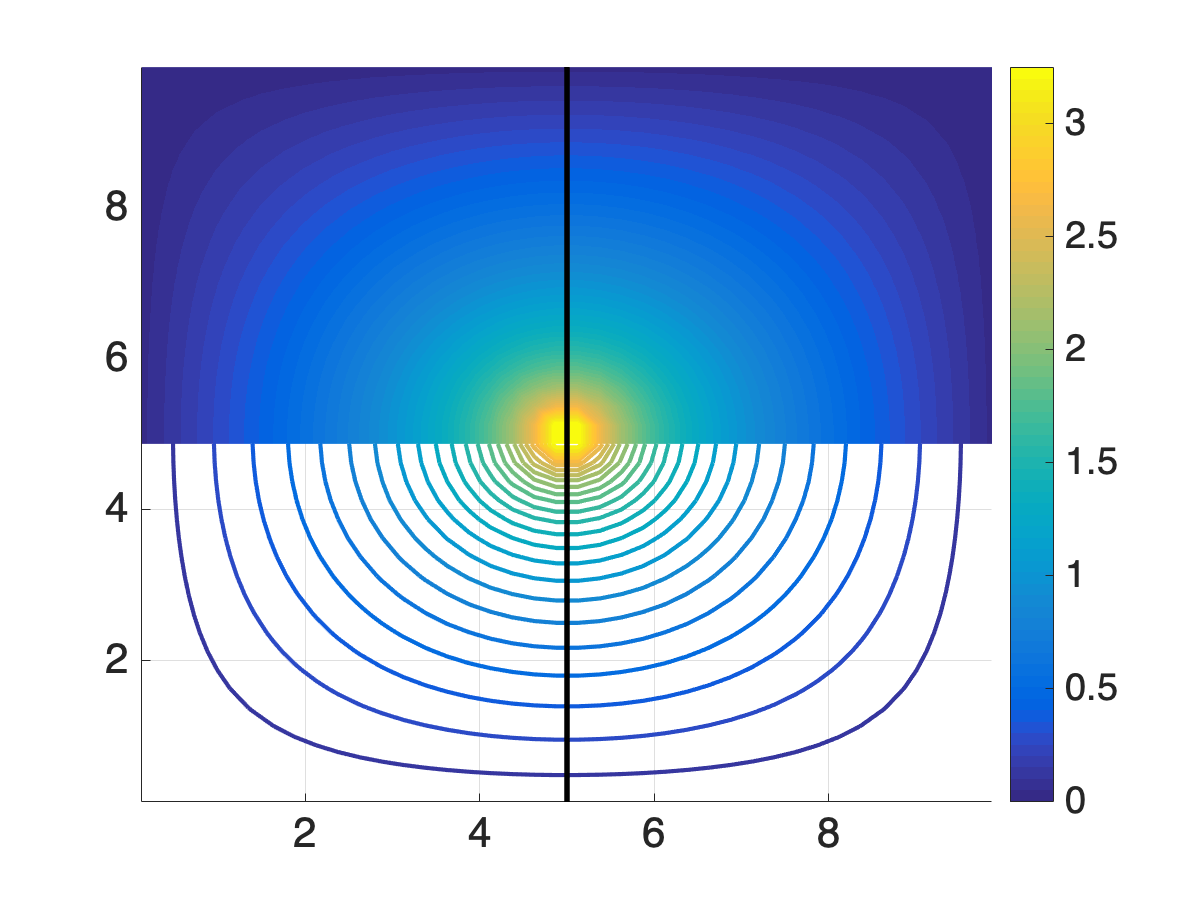}
}
\subfigure[Example 3, $10$ reduced basis\label{fig:2d_test3}]{
 \includegraphics[width=0.33\textwidth]{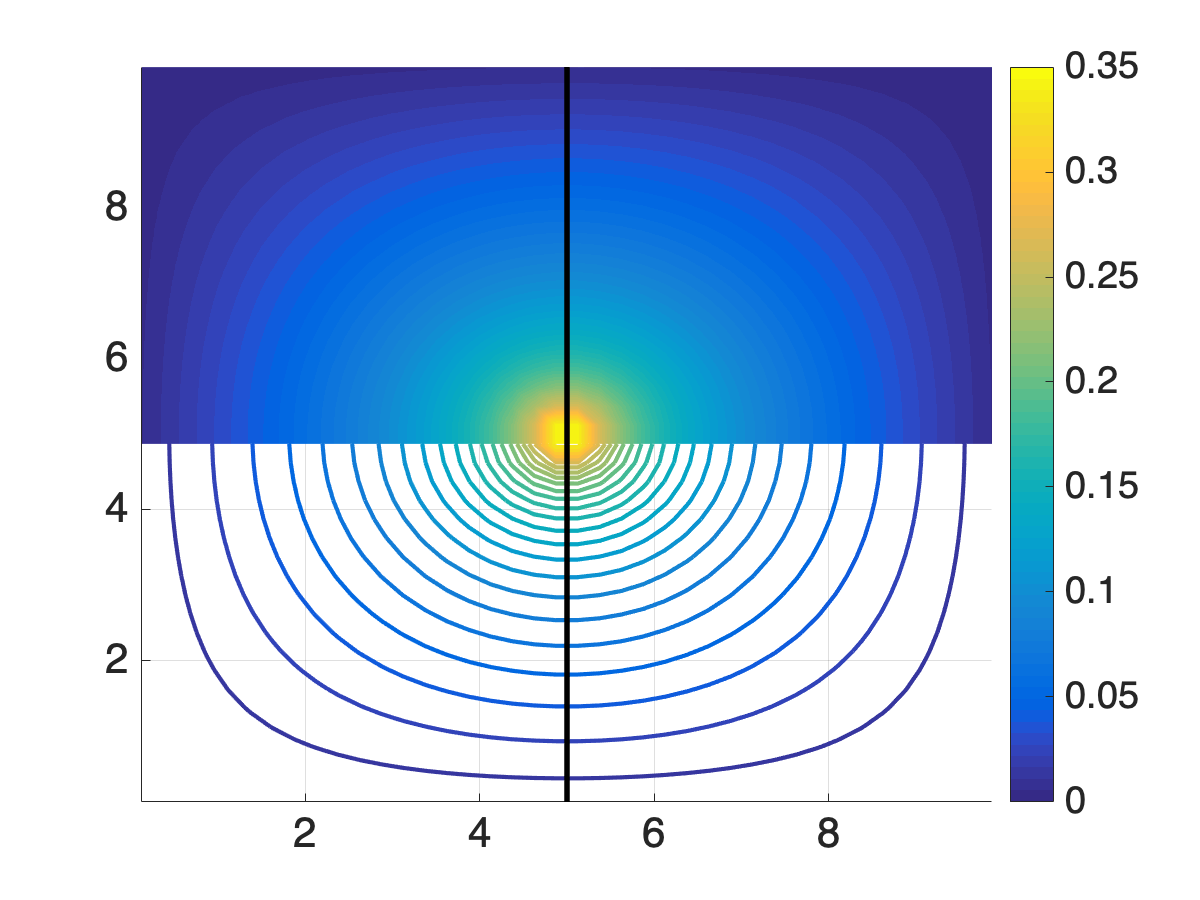}
}
\subfigure[Example 4, initial guess ($S_4$)\label{fig:2d_test4_rho_initial}]{
 \includegraphics[width=0.33\textwidth]{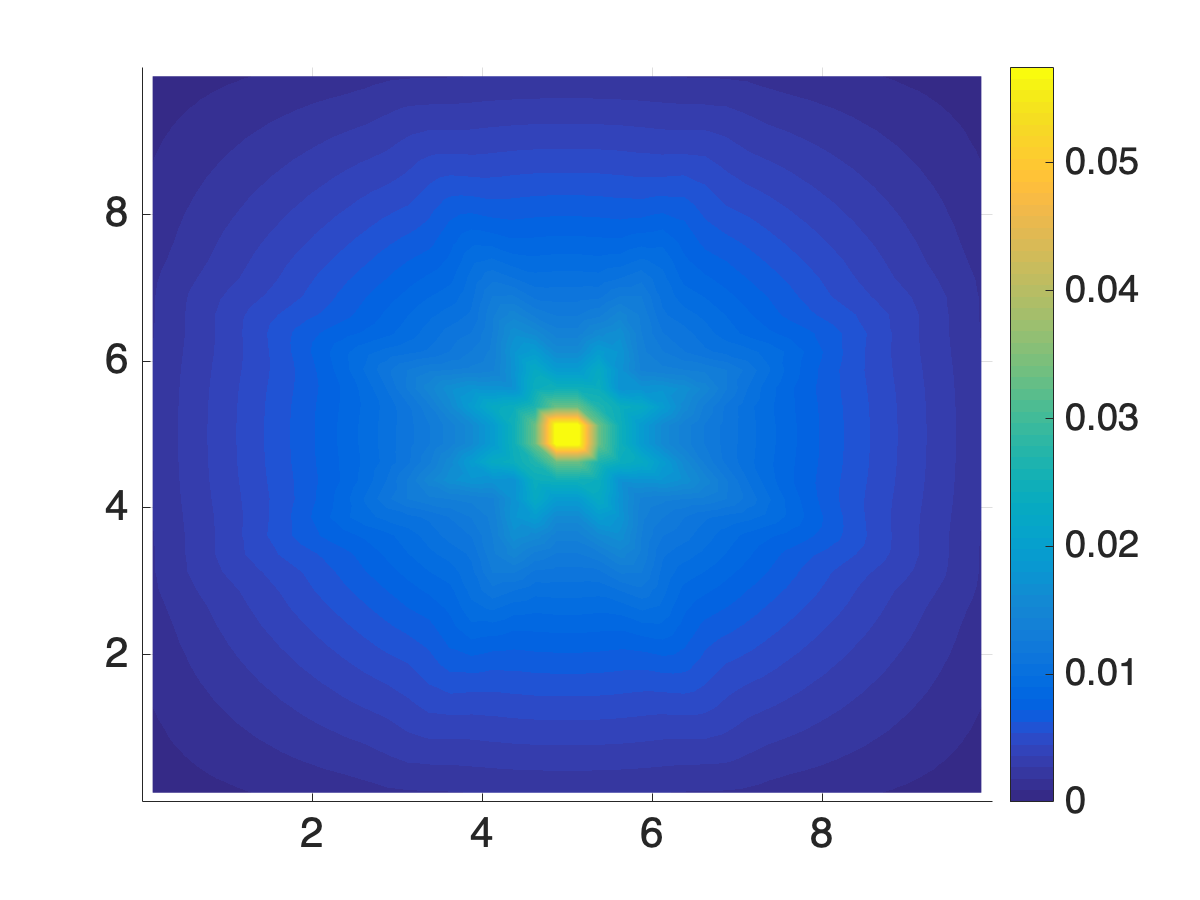}
}
\subfigure[Example 4, \pzc{$26$} reduced basis\label{fig:2d_test4_rho_rb}]{
 \includegraphics[width=0.33\textwidth]{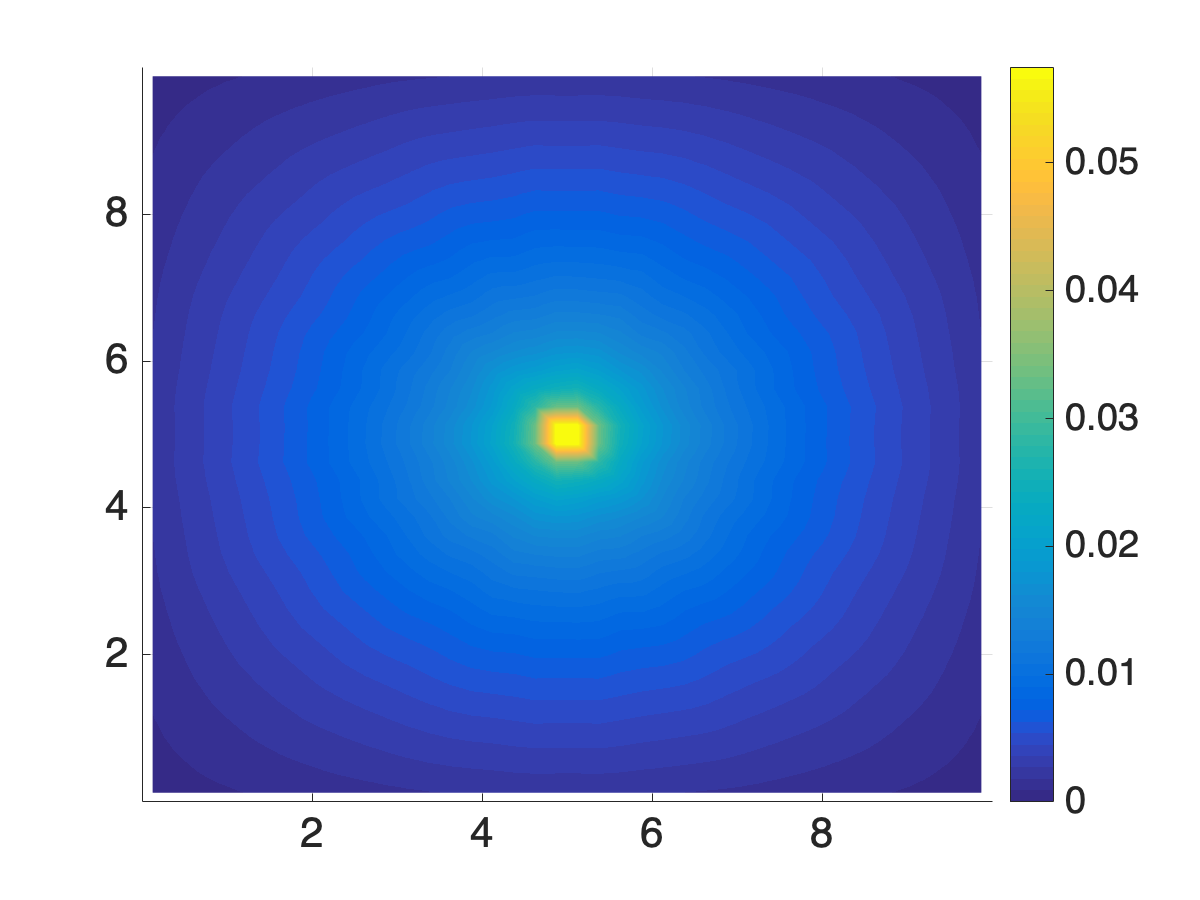}
}
\subfigure[Example 4, \pzc{$26$} reduced basis \label{fig:2d_test4}]{
 \includegraphics[width=0.33\textwidth]{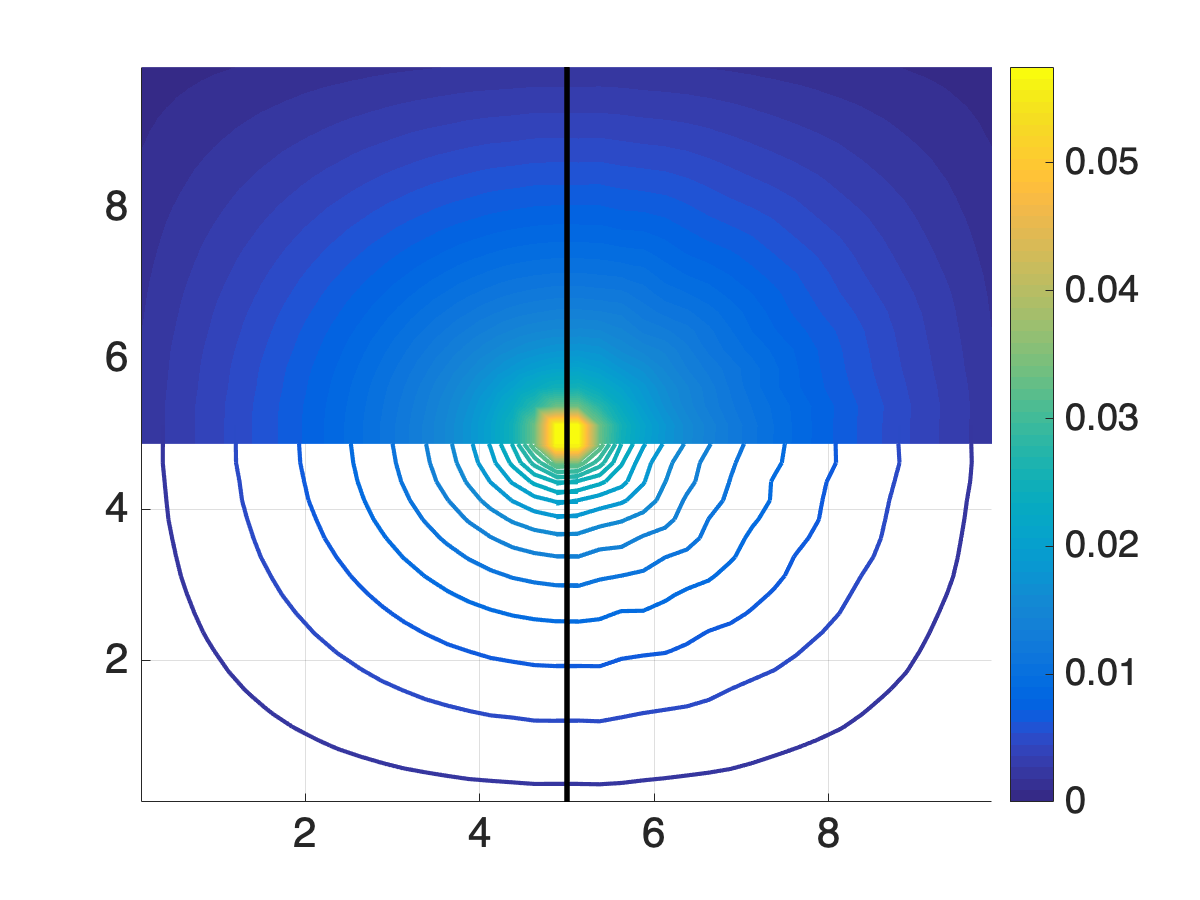}
}
\caption{Density for 2D examples. Top row: Examples 1-3, $r_{\rm tol}=10^{-3}$; Bottom row: Example 4, $r_{\rm tol}=10^{-2}$.  In each plot - left half: full order solution; right half: RB solution; Top half: full/RB solution; bottom half: contour of the full/RB solution. 
\label{fig:2d_test_rho}}
\end{figure}

\fli{The comparison between the densities obtained by the RB method (trained with $N_\BOmega=32$) and the full order method (with $40$ quadrature points) are presented in Figure \ref{fig:2d_test_rho}, with the former computed as the terminal density when  $r_{\rm tol}=10^{-3}$ is reached  for Examples 1-3 and $r_{\rm tol}=10^{-2}$ is reached for Example 4.}
 \pzc{ We see that the RB solution and the full order solution match each other well.  Moreover, the second row of Figure \ref{fig:2d_test_rho} demonstrates the effectiveness of our method in mitigating the ray effect \cite{lathrop1968ray}, which refers to the phenomenon that  the particles mainly propagate along the directions of sampled angular directions and the numerical solution has noticeable unphysical oscillations. As shown in Figure \ref{fig:2d_test4_rho_initial}, the initial guess for the problem in the transport regime suffers severely from the ray effect. \fli{Even with an initial guess of poor quality,  the ray effect in the reduced order solution is less pronounced}  as the dimension of the reduced space grows (see Figure \ref{fig:2d_test4_rho_rb} and Figure \ref{fig:2d_test4}).}

\pzc{
\begin{table}[htbp]
  \centering
 \medskip
    \begin{tabular}{|l|c|c|c|c|c|c|}
    \hline
 		&RB dimension &  	$\mE_f$&   $R_f$ & $\mE_\rho$  & $R_\rho$ \\ \hline
 Example 1&	8	&	1.50e-3&  1.21e-2\% & 7.69e-4 & 3.44e-2\% \\ \hline			
 Example 2&	4	&	2.40e-3&  6.23e-3\% & 2.40e-3 & 3.49e-2\%  \\ \hline	
 Example 3&	10	&	7.10e-4&  1.83e-2\% & 2.77e-4& 4.03e-2\%   \\ \hline
  Example 4&	26	&	2.95e-4&  9.97e-1\% & 9.72e-4& 1.15e-0\%  \\ \hline
 \end{tabular}
      \caption{\pzc{training errors for 2D Examples 1-3 with $r_{\rm tol}=10^{-3}$ and  2D Example 4  with $r_{\rm tol}=10^{-2}$ and $N_\BOmega=32$.} \label{tab:2d_test1_1e-3} }
\end{table}
}

\begin{figure}
\subfigure[Error history of $f$\label{fig:2d_test_error_history}]{
 \includegraphics[width=0.5\textwidth]{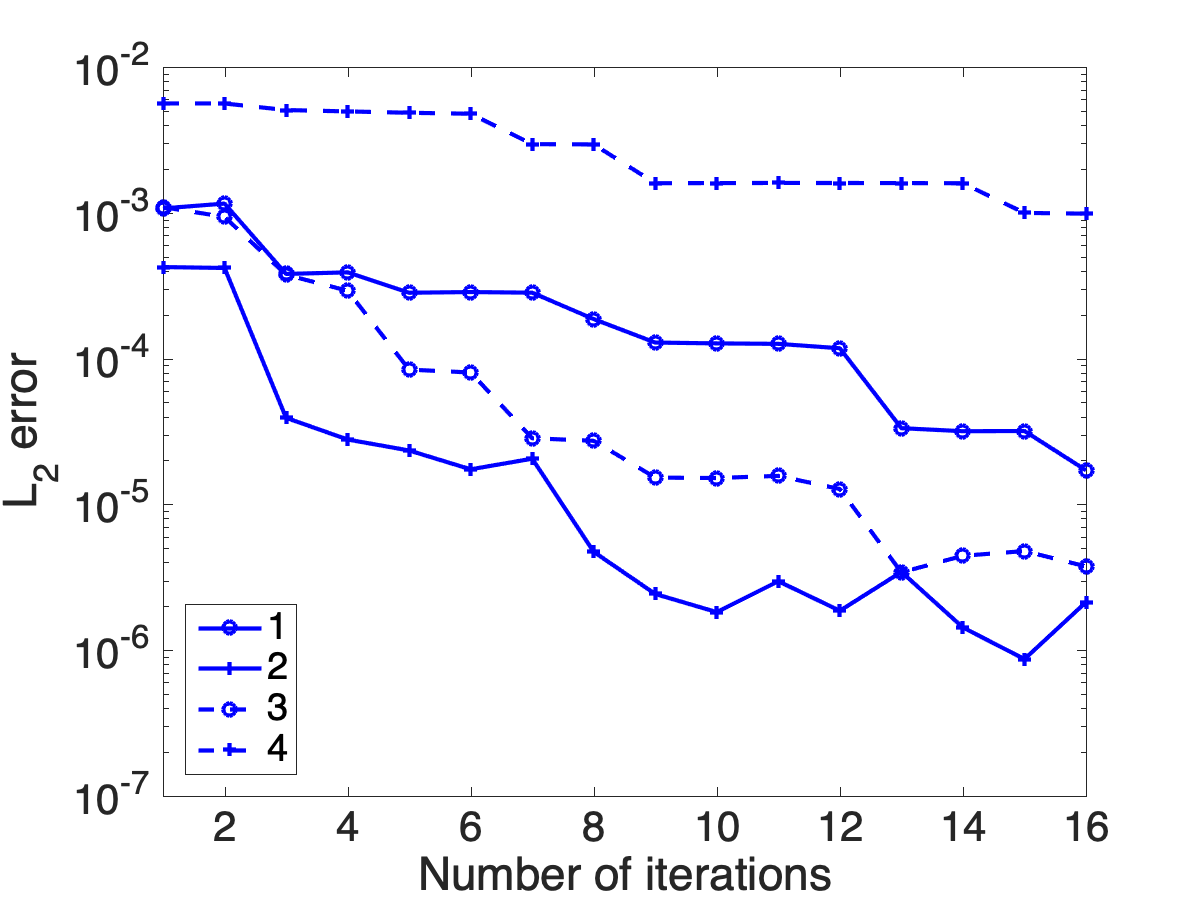}
 }
 \subfigure[History of spectral ratio 
 \label{fig:2d_test_svd_history}]{
 \includegraphics[width=0.5\textwidth]{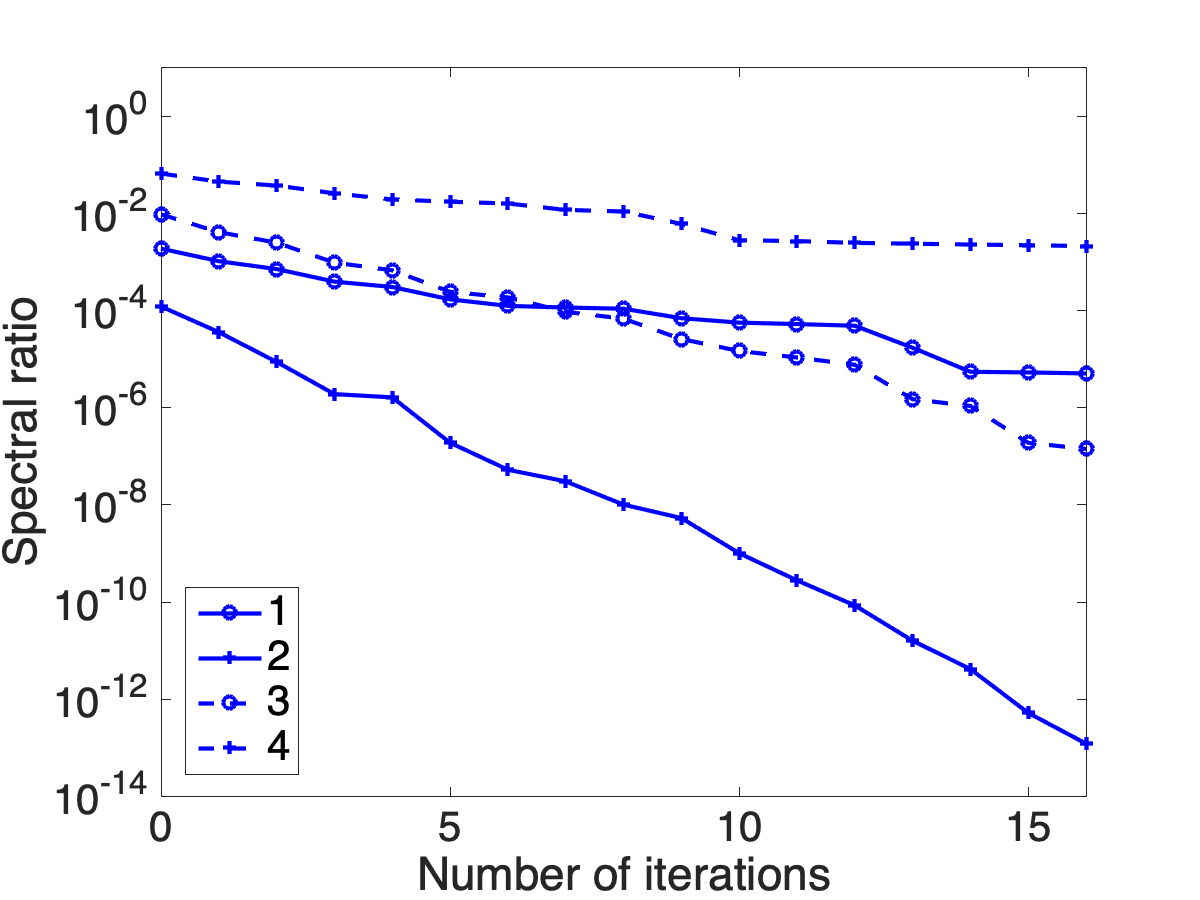}
 }
 \caption{Histories of convergence 
 of the \fli{$L^2$ error of $f$} and the monitored spectral ratio during training for 2D Examples 1-4  \label{fig:2d_test_history}}
\end{figure}

\noindent {\bf \pzc{Effectiveness of the error indicator:}} 
\pzc{In the last test,} 
we ask the algorithm to generate $M_{\rm tol}=16$ reduced bases for all problems and monitor the \fli{training} errors.  
The histories of convergence for {the $L^2$ errors of $f$}  and the spectral ratio $r^m$  are presented in Figure \ref{fig:2d_test_history}. Similar to the one-dimensional cases, we see that, as the dimension of the reduced space increases, the overall trend of the error and spectral ratio is decreasing for all examples, albeit slower than the one-dimensional cases. 
Moreover, the different speed of decay of the spectral ratio $r^m$ for different problems leads to the different terminal RB dimensions as shown in Tables \ref{tab:2d_test1_1e-3}. 
The fact that these quantities decay slower for transport dominant problems is consistent with the slow decay of the Kolmogorov $N$-width for transport problems, see e.g. \cite{greif2019decay,ohlberger2015reduced}. As a result, the performance of our method and many other linear ROMs suffer. 

\fli{Finally we want to mention that in Appendix \ref{appendix:b}, we plot and comment about the leading reduced basis functions generated by the proposed algorithm for selected 1D and 2D examples. }

\section{Conclusion}\label{sec:conclusion}
In this paper, we design a RB method to construct an angular-space reduced order model for the linear {radiative transfer equation.} Unlike the standard setting where RBM applies, the solutions for different parameter values are coupled through an integration operator. 
This coupling makes impossible the direct inquiry of the snapshot for any particular parameter value. 
An additional challenge is that parameter ensemble identified by traditional RBM is usually unstructured, and thus may not form a set of  quadrature points for a robust and accurate integration toward the macroscopic density which is however crucial for the efficiency and robustness of the full order SASI solver.

Via a careful iterative procedure where the macroscopic density is treated explicitly allowing a transport sweep and then updated afterwards, a least squares density reconstruction at each of the relevant physical locations, a $L^1$-based residual-free error indicator, and a symmetry enhancing greedy addition, we successfully designed the first RBM for the kinetic transport equation. Our numerical experiments indicate that the new method is highly effective for the scattering dominant problems, the intermediate regime problems, and the multiscale problems with scattering dominant subregions. Moreover, as the reduced space grows,  the ray effect can be mitigated. While less efficient for transport dominant problems, it is capable of decreasing the problem size by more than one order of magnitude and achieving one digit of accuracy. Designing a RBM that works equally well for transport dominant problems, likely by constructing nonlinear reduced manifolds, constitutes our future work. Other future work include the application to the full 5D model and the extension to time dependent problems and other kinetic models.

\begin{appendices}
\pzc{
\section{RB method  with the diffusion synthetic acceleration\label{sec:DSA}}
}
\pzc{
In this work,  we apply the S2SA 
to speed up the convergence of the source iteration  when algebraically solving the upwind DG discretization. One main advantage of the S2SA 
is that the same kind of kinetic solver as the full order one is applied. 
Another type
widely used synthetic acceleration is the diffusion synthetic acceleration (DSA) \cite{adams1992diffusion,wareing1993new,adams2002fast}. Instead of using a low order $S_N$ model for a kinetic problem to approximate the correction equation \eqref{eq:kinetic_correction} for $\delta f^{k+1}$ first and then to compute $\rho^{k,c}=\lgl \delta f^{k+1}\rgl$, the DSA method works with a discrete diffusion approximation, that is ``consistent'' (see  \cite{adams2002fast} for the definition of the consistency),  to approximate $\rho^{k,c}$ directly. It is known that the source iteration with ``inconsistent" DSA may converge slowly or even diverge in some regimes \cite{adams2002fast}. Compared with a straightforward S2SA method,  the DSA will involve fewer degrees of freedom. Next we will use the 1D model  on the slab geometry as an example to present a DSA method  that is consistent to the upwind DG discretization, and then demonstrate and compare the performance of the RB method with both synthetic acceleration strategies. 
}
\pzc{
\subsection{A consistent DSA method}
For the  1D slab geometry with $X=[0,1]$,  by using an ansatz $\delta f(x,v)=\rho^{k,c}(x)+3v g(x)$ and taking the zeroth, first moments of the correction equation \eqref{eq:kinetic_correction} in the angular variable, we obtain an approximated diffusion model 
\begin{subequations}
\label{eq:diffusion_correction}
\begin{align}
\underbrace{\lgl v^2\rgl}_{=\frac{1}{3}} &\partial_x\;({\sigma_t}^{-1}\partial_x \rho^{k,c})+\sigma_a \rho^{k,c}=\sigma_s(\rho^{k,*}-\rho^k), \quad \textrm{on}\; X
\end{align}
\end{subequations}
or, equivalently, in its first order form,
\begin{subequations}
\label{eq:p1_correction}
\begin{align}
\partial_x g + \sigma_a \rho^{k,c} = \sigma_s(\rho^{k,*}-\rho^k),\quad 
\underbrace{\lgl v^2\rgl}_{=\frac{1}{3}} \partial_x \rho^{k,c}+\sigma_t g = 0.
\end{align}
\end{subequations}
They will be complemented by the boundary conditions,
\begin{subequations}
\label{eq:diffusion_correctionBC}
\begin{align}
 & \lgl v(\rho^{k,c}+3vg) \rgl^{+}=\frac{1}{4}(\rho^{k,c}+2g)=\frac{1}{4}(\rho^{k,c}-\frac{2}{3\sigma_t}\partial_x \rho^{k,c}) =0,\quad\text{at}\; x=0,\\
 & \lgl v(\rho^{k,c}+3vg) \rgl^{-}=-\frac{1}{4}(\rho^{k,c}-2g)=- \frac{1}{4}(\rho^{k,c}+\frac{2}{3\sigma_t}\partial_x \rho^{k,c})=0,\quad\text{at}\; x=1.
\end{align}
\end{subequations}
Here $\lgl \eta\rgl^+=\frac{1}{2}\int_{v>0} \eta(v)dv$, $\lgl \eta \rgl^-=\frac{1}{2}\int_{v<0} \eta(v) dv$.
}

\pzc{
Let $\{T_i=[x_{i-\half},x_{i+\half}], \; i=1\dots, N_x\}$ be a partion of the domain $X=[0,1]$, we then discretize \eqref{eq:p1_correction} with a DG method: we seek $\rho^{k,c}_h, g_h\in U_h^K$, such that $\forall \phi_h, \psi_h\in U_h^K$, $\forall i=1,\dots, N_x$,  
\begin{subequations}
\label{eq:DSA:DG}
\begin{align}
&-\int_{I_i} g_h \partial_x \phi_h dx+ (\widehat{g_h} \phi^-_{h})_{i+\half}-(\widehat{g_h} \phi^+_{h})_{i-\half}+\int_{I_i} \sigma_a \rho_h^{k,c}  \phi_h dx = \int_{I_i} \sigma_s(\rho^{k,*}-\rho^k) \phi_h dx,\\
&\frac{1}{3} \left(-\int_{I_i} \rho_h^{k,c} \partial_x \psi_h dx + (\widehat{\rho_h^{k,c}} \psi^-_{h})_{i+\half}-(\widehat{\rho_h^{k,c}} \psi^+_{h})_{i-\half}\right)+\int_{I_i}  \sigma_t g_h \psi_h dx= 0.
\end{align}
\end{subequations}
The key to make the method \eqref{eq:DSA:DG} consistent to the upwind DG method for the transport sweep step lies in the numerical fluxes $(\widehat{g_h})_{i+\half}$ and $(\widehat{\rho_h^{k,c}})_{i+\half}$, that shall be based on the upwind flux for $\delta{f}$.\footnote{\pzc{In our RB method, we use the exact values of $\lgl v^2\rgl$ in \eqref{eq:DSA:DG} and $\lgl v^k\rgl^\pm, k=1, 2, 3$ in \eqref{eq:dsa_flux}. In literature, numerical integrations with certain property are often used.}}  With this in mind, we take 
\begin{subequations}
\label{eq:dsa_flux}
\begin{align}
(\widehat{g_h})_{i+\half}&= \lgl v\widehat{\delta f}^{\textrm{upwind}}\rgl_{i+\half} =
\lgl v (\rho_h^{k,c}+3v g_h)\rgl^+(x_{i+\half}^-) +  \lgl v (\rho_h^{k,c}+3v g_h)\rgl^-(x_{i+\half}^+) \notag\\
&=\{g_h\}_{i+\half}-\frac{1}{4}[\rho_h^{k,c}]_{i+\half}, \\ 
\frac{1}{3}(\widehat{\rho_h^{k,c}})_{i+\half} &= \lgl v^2\widehat{\delta f}^{\textrm{upwind}}\rgl_{i+\half} =\lgl v^2 (\rho_h^{k,c}+3v g_h)\rgl^+(x_{i+\half}^-) +  \lgl v^2 (\rho_h^{k,c}+3v g_h)\rgl^-(x_{i+\half}^+)\notag\\
&= \frac{1}{3}\left(\{\rho_h^{k,c}\}_{i+\half}-\frac{9}{8}[g_h]_{i+\half}\right).
\end{align}
\end{subequations}
Here
$u(x_{i+\half}^-)=u_{i+\half}^-$ (resp. $u(x_{i+\half}^+)=u_{i+\half}^+$) stands for the left (resp. right) limit of $u(x)$ at the cell interface $x_{i+\half}$. And $\{u\}_{i+\half}=\half(u_{i+\half}^++u_{i+\half}^-), \;[u]_{i+\half}=u_{i+\half}^+-u_{i+\half}^-$. At boundaries, the numerical fluxes are set as \eqref{eq:dsa_flux} with $\lgl v^k (\rho_h^{k,c}+3v g_h)\rgl^+(x_{\half}^-)=0$, $\lgl v^k (\rho_h^{k,c}+3v g_h)\rgl^-(x_{N_x+\half}^+)=0, \; k=1, 2$. 
In actual simulation, we eliminate $g_h$ in \eqref{eq:DSA:DG} at the algebraic level, leading to a smaller linear system for $\rho_h^{k,c}$ only, that is given  in its matrix-vector form as follows,
\begin{align}
\label{eq:consistent_dg_dsa}
\left(\BSIGMA_a -\frac{1}{4}\BD_{\rm jump} -\frac{1}{3}\BD_c(\BSIGMA_t-\frac{3}{8}\BD_{\textrm{jump}})^{-1}\BD_c \BRHO^{k,c}\right)\BRHO^{k,c}=\BSIGMA_S(\BRHO^{k,*}-\BRHO^{k}), 
\end{align}
with $\BD_c= \frac{1}{2}\left(\BD^++\BD^-\right)$, $\BD_\textrm{jump} =\BD^+-\BD^-$, where
\begin{subequations}
\begin{align}
(\BD^+)_{kl}&=-\sum_{i=1}^{N_x}\int_{T_i} \partial_x\phi_k(x) \phi_l(x) dx +\sum_{i=1}^{N_x-1}\phi_l(x^+_{i+\half}) \phi_k(x^-_{i+\half})-\sum_{i=1}^{N_x}\phi_l(x^+_{i-\half}) \phi_k(x^+_{i-\half}), \\
(\BD^-)_{kl}&=-\sum_{i=1}^{N_x}\int_{T_i} \partial_x\phi_k(x) \phi_l(x) dx +\sum_{j=1}^{N_x}\phi_l(x^-_{i+\half}) \phi_k(x^-_{i+\half})-\sum_{i=2}^{N_x}\phi_l(x^-_{i-\half}) \phi_k(x^+_{i-\half}).
\end{align}
\end{subequations}
One example of ``partially consistent'' DSA methods is to use central fluxes as in \cite{adams1992diffusion}. Partially consistent DSA methods may result in slower convergence. 
}

 \pzc{
 \subsection{Comparison of the RB method with the S2SA and the DSA}
 For the 1D examples, we replace the S2SA with the DSA in the source iteration, and compare the performance of the overall RB algorithm with the two different acceleration strategies.
The stopping criteria is $r_{\rm tol}\leq 10^{-4}$. For all  examples, the full order solvers with the S2SA and the DSA in the source iteration always converge to the same result.  The errors for the RB method with the S2SA and the DSA  are reported in Table \ref{tab:error_dsa_rb}, and the relative computational time, defined as (time with S2SA)/(time with DSA), is summarized in Table \ref{tab:s2sa_vs_dsa}. The full order model with the DSA is slightly more efficient in the diffusive regime, and it is comparable with the S2SA method in other regimes. For the diffusive and intermediate regimes \fli{(Examples 1-2)}, the RB method  with the DSA is more efficient in both offline, but the \fli{errors are relatively larger.} The angular samples picked by the greedy algorithm in the RB method are the same for both the DSA and S2SA, and \pzc{hence the online computational costs of the RB method with the DSA and the S2SA are close to each other. }
For Example 4 (two-material problem) and 
Example 5 (transport regime), the RB method with the DSA  fails to converge to the correct solution.  It is known that reduced order methods can be more sensitive to the choice of preconditioners compared with full order solvers \cite{ChenGottliebMaday, ZahmNouy2016, SantoDeparisManzoniQuarteroni2018}.  In summary, when combined with our RB method, the S2SA is more robust with various regimes and slightly more accurate. The RB method with the DSA is slightly more efficient if it converges, but it may fail to converge for problems with transport-dominant (sub)regions.
}

\begin{table}[!htb]
\centering
\pzc{
\begin{tabular}{|l|c|c|c|c|c|c|}
 \hline
 		&RB dimension &  	$\mE_f$-DSA&   $\mE_f$-S2SA& $\mE_\rho$-DSA &  $\mE_\rho$-S2SA\\ \hline
 Example 1&	4	&	6.26e-2&  8.04e-3 & 6.25e-2 & 8.00e-3\\ \hline	
 Example 2&	4	&	6.99e-2&  8.51e-3 & 6.98e-2 & 8.49e-3\\ \hline	
 Example 3&	4	&	3.28e-3&  9.35e-4 & 2.88e-3 & 8.25e-4\\ \hline	
 Example 4&	\multicolumn{5}{c|}{RB-DSA not convergent}  \\ \hline	
 Example 5&	\multicolumn{5}{c|}{RB-DSA not convergent}   \\ \hline	
 \end{tabular}
 }
\caption{\pzc{Errors of the RB method with the DSA and S2SA for 1D examples. \label{tab:error_dsa_rb}}}
\end{table}

 \begin{table}[!htb]
\centering
 \pzc{
\begin{tabular}{|c|c|c|c|c|c|c|}
\hline
 	 & Example 1 & Example 2 & Example 3 & Example 4  & Example 5\\
\hline 
RB: offline & $4.68$  & $4.67$ &  $3.37$ & \multicolumn{2}{c|}{RB-DSA not convergent}\\
\hline 
Full      & $1.57$ & $1.05$ &   $ 1.34$ & $1.12$ & 0.96\\
\hline 
\end{tabular}
}
\caption{\pzc{Relative computational time for 1D examples:  (time with S2SA)/(time with DSA).   \label{tab:s2sa_vs_dsa}}}
\end{table}

\section{\fli{Plots of selected RB functions}}
\label{appendix:b}

\begin{figure}
\centering
 \includegraphics[width=0.32\textwidth]{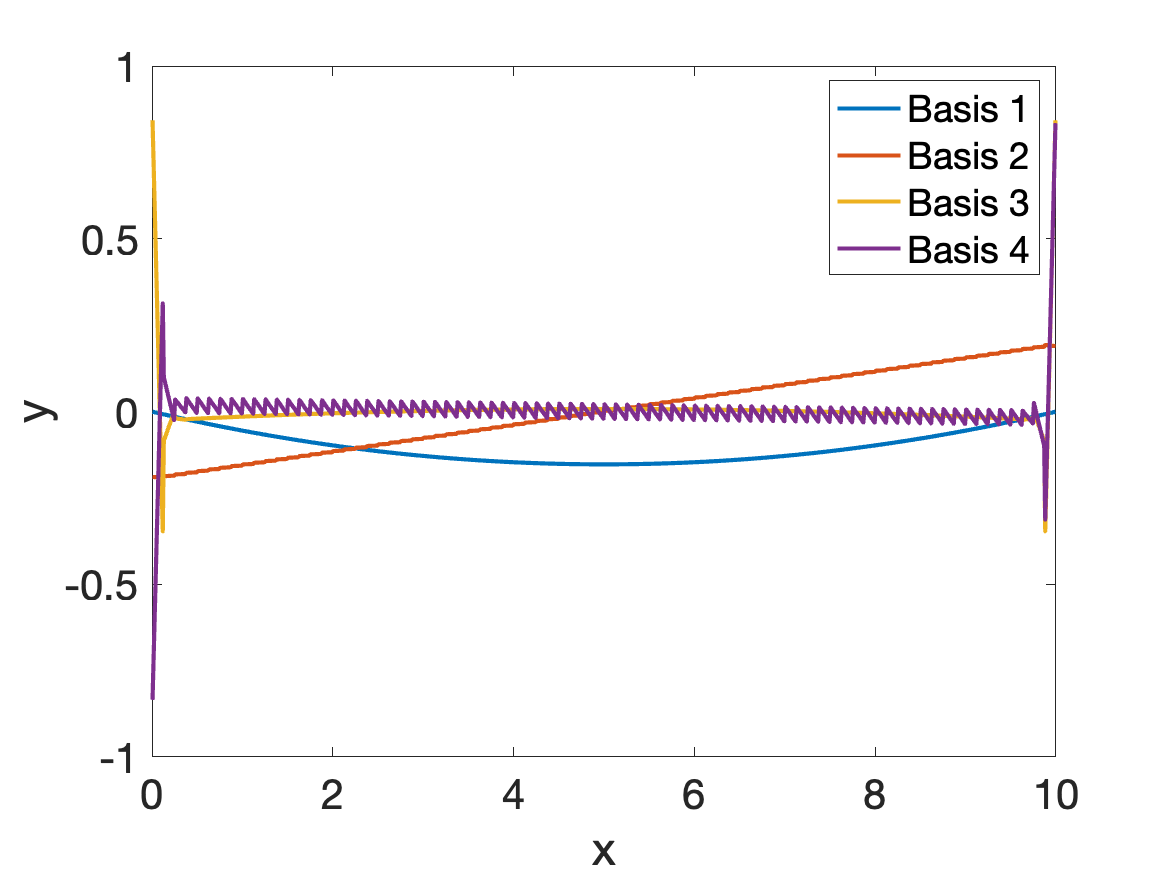}
 \includegraphics[width=0.32\textwidth]{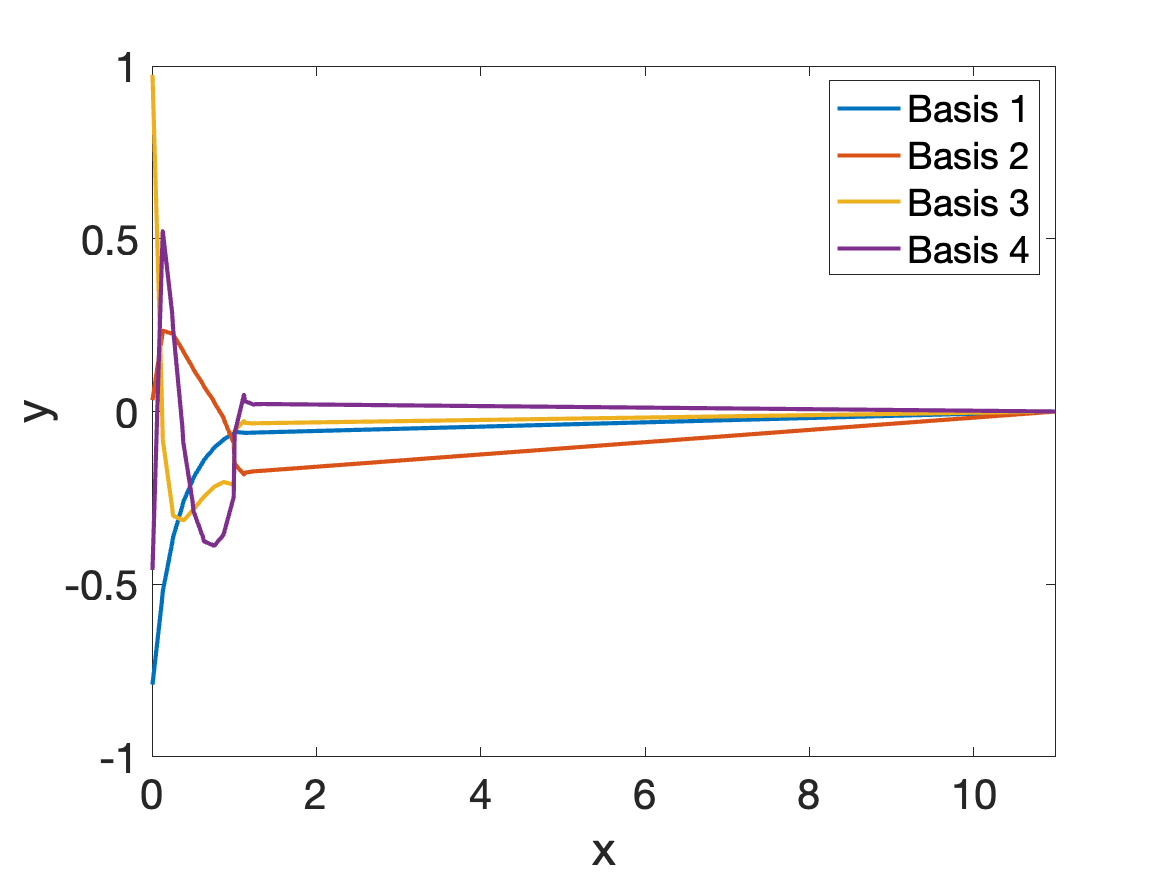}
 \includegraphics[width=0.32\textwidth]{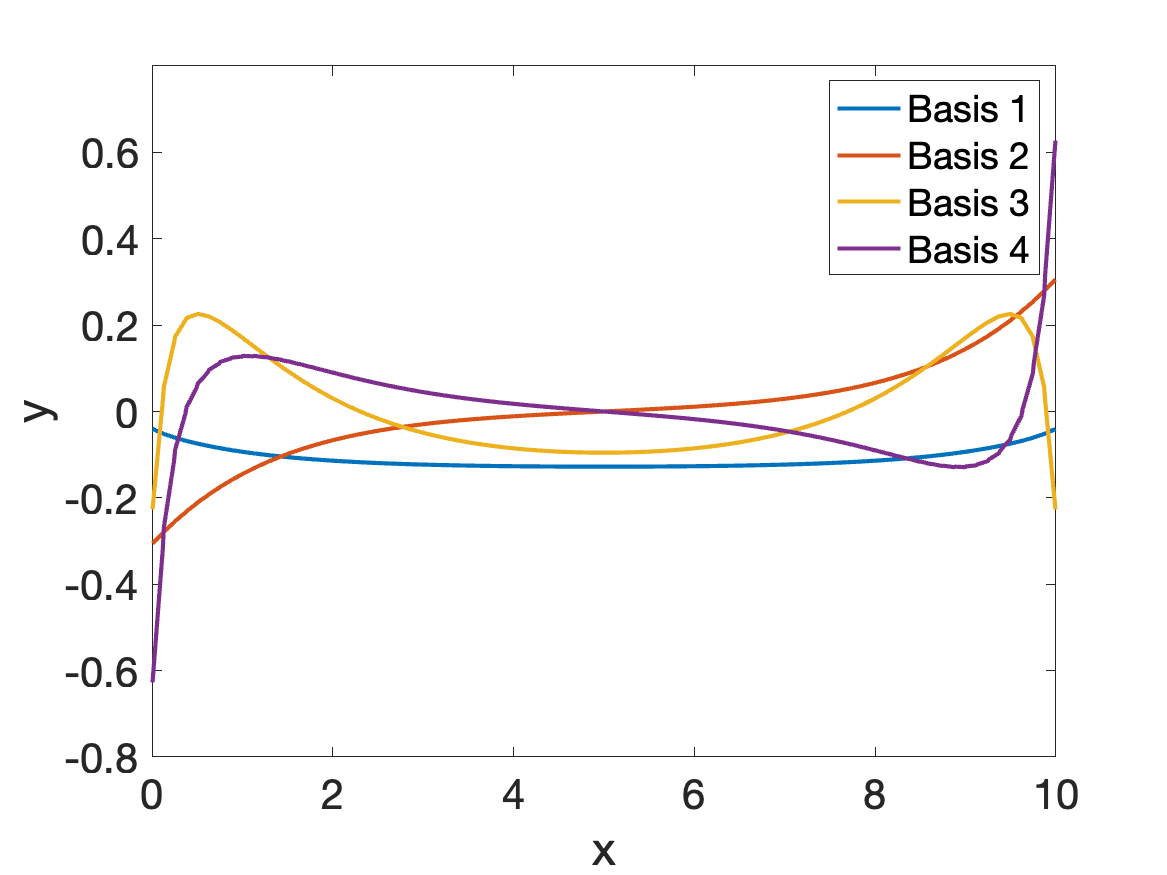}
 \caption{\fli{First four reduced basis functions for 1D Examples 1 (left), 4 (middle), 5 (right).
 \label{fig:1d_basis}}}
\end{figure}

\fli{
We here present some selected RB functions generated by the proposed algorithm.  In Figure \ref{fig:1d_basis}, we plot the first four RB functions (after the SVD orthogonalization step) for 1D  Examples 1, 4, 5 from Section \ref{sec:numericaL^1d}. The most interesting example is the two-material problem in Example 4. Particularly, the presence of a material interface is captured by most RB functions. Moreover, all four RB functions behave fairly  differently in the left subregion where the problem is transport dominant. 
}

\begin{figure}
\centering
 \includegraphics[width=0.4\textwidth]{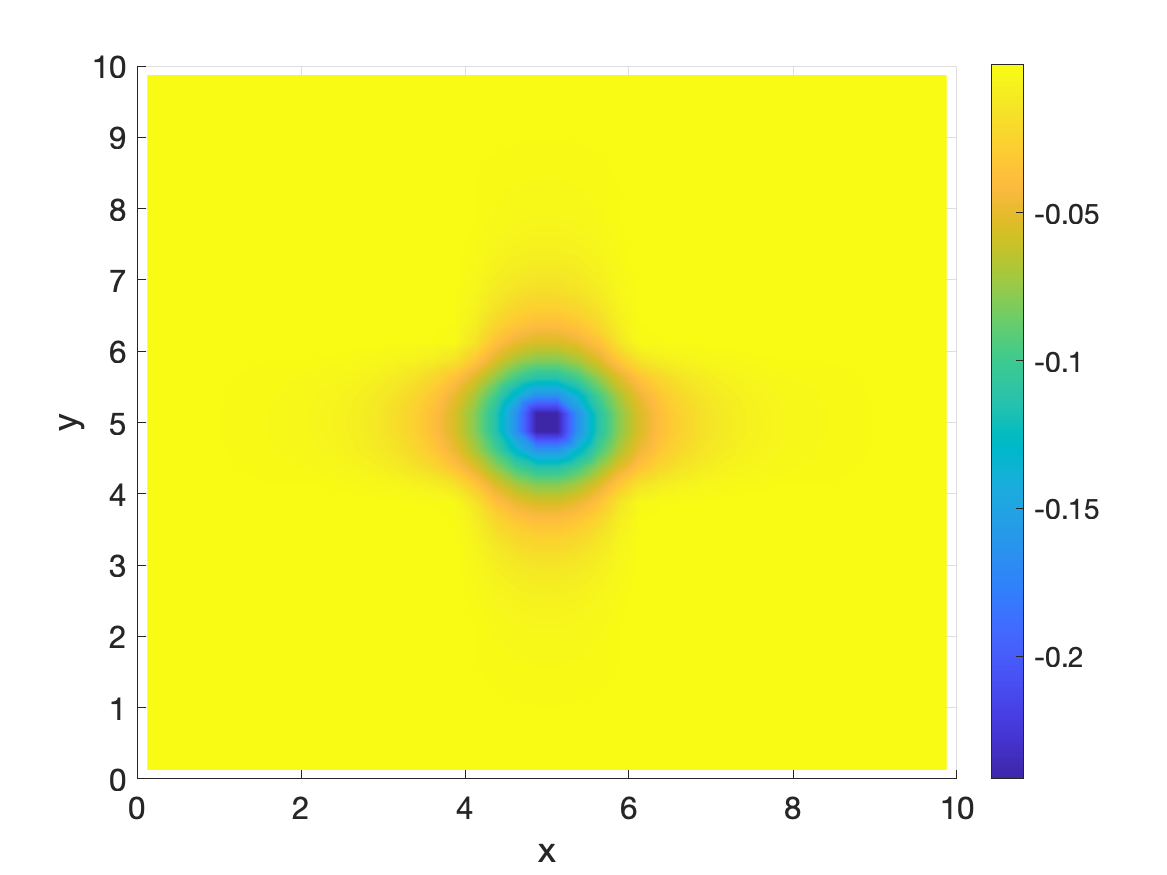}
 \includegraphics[width=0.4\textwidth]{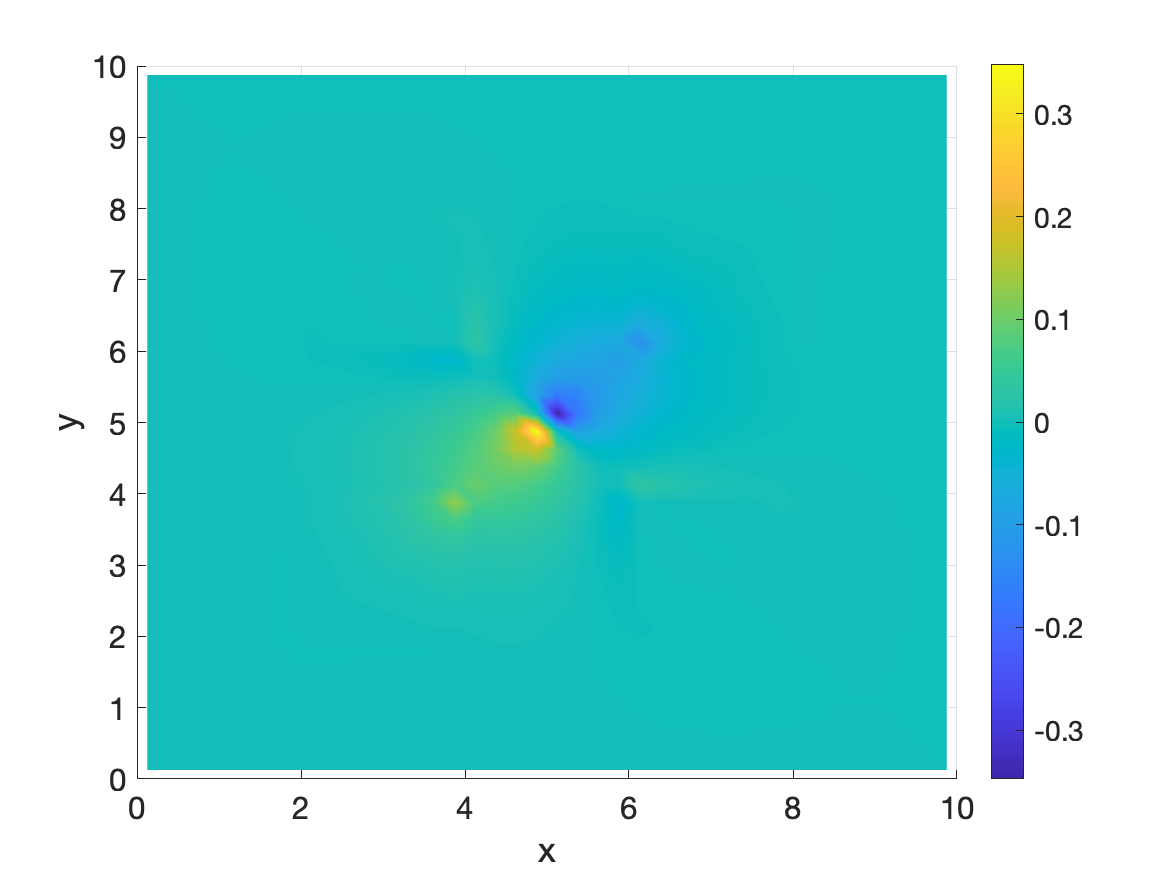}
 \includegraphics[width=0.4\textwidth]{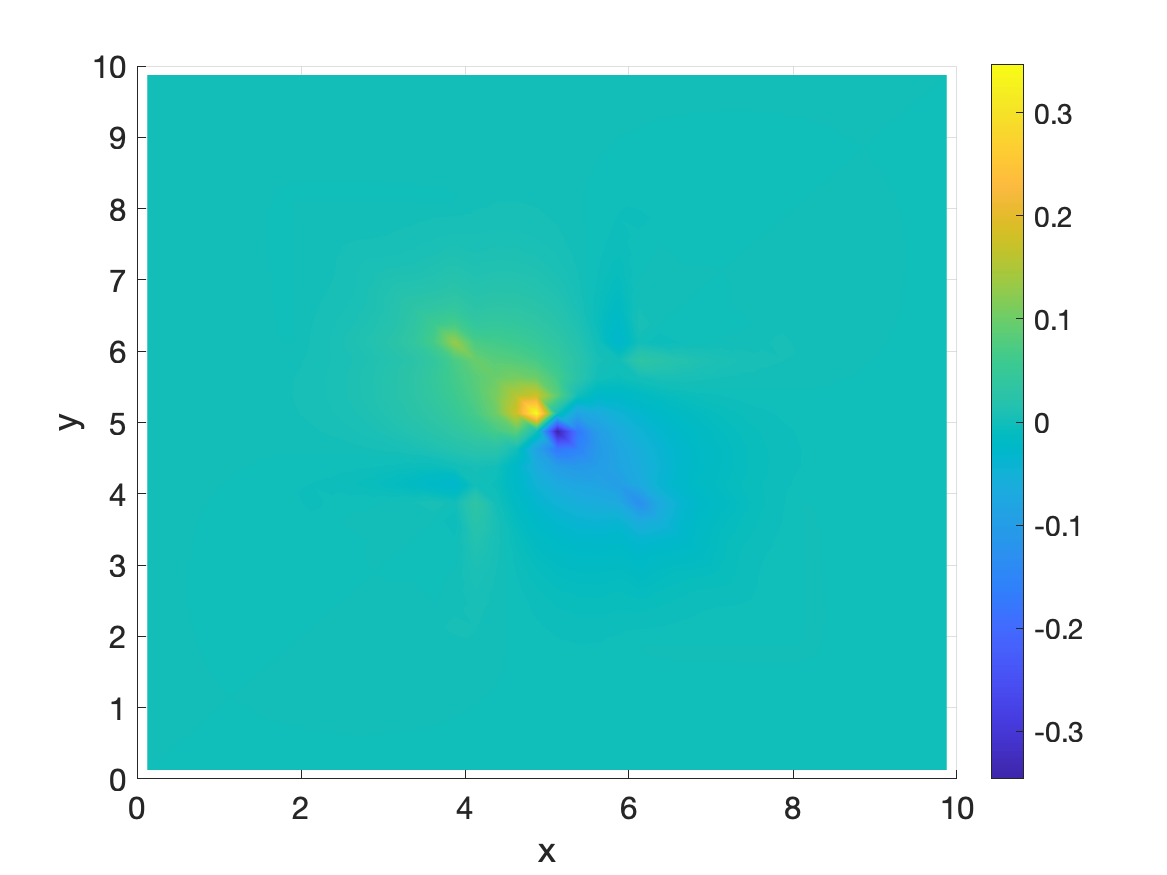}
 \includegraphics[width=0.4\textwidth]{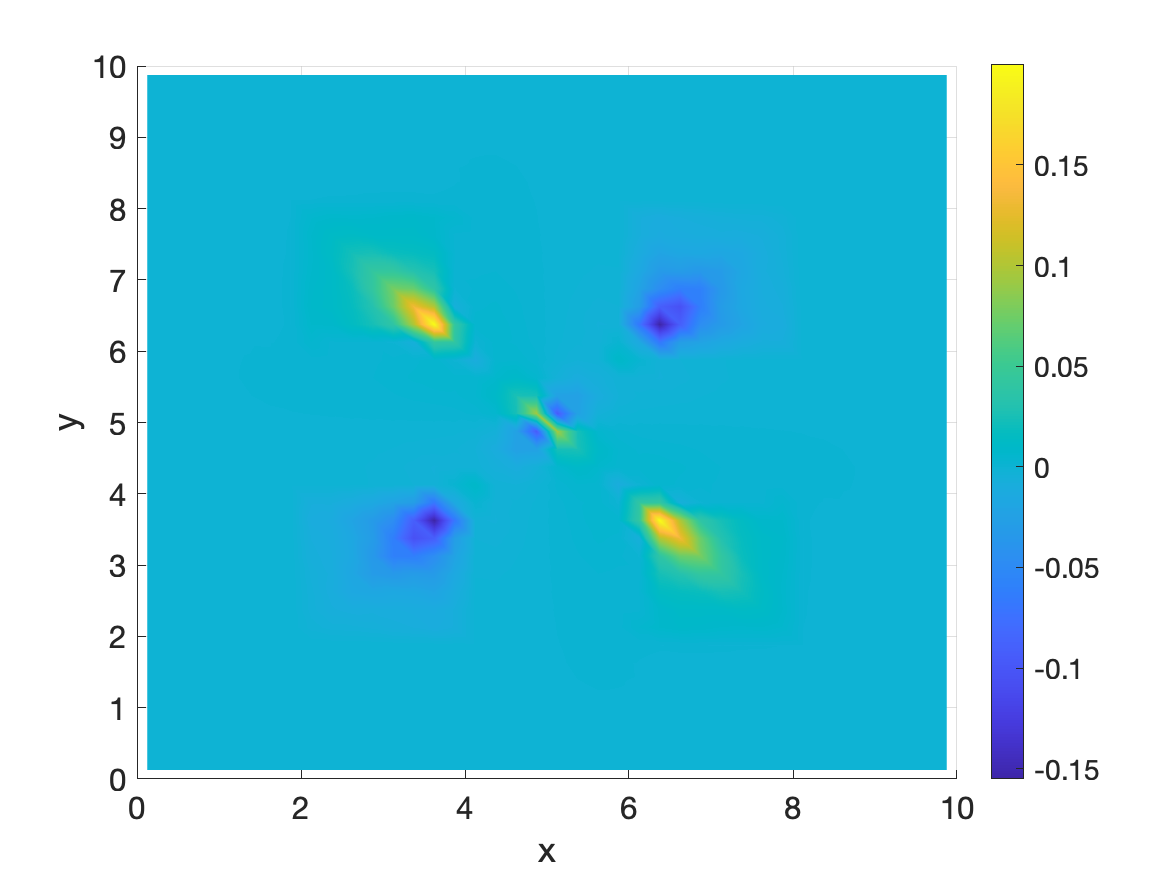}
 \caption{\pzc{First four reduced basis functions, \fli{ordered from left to right and from top to bottom}, for 2D Example 1 (checkerboard problem).\label{fig:2d_basis_test1}}}
\end{figure}

\begin{figure}
\centering
 \includegraphics[width=0.4\textwidth]{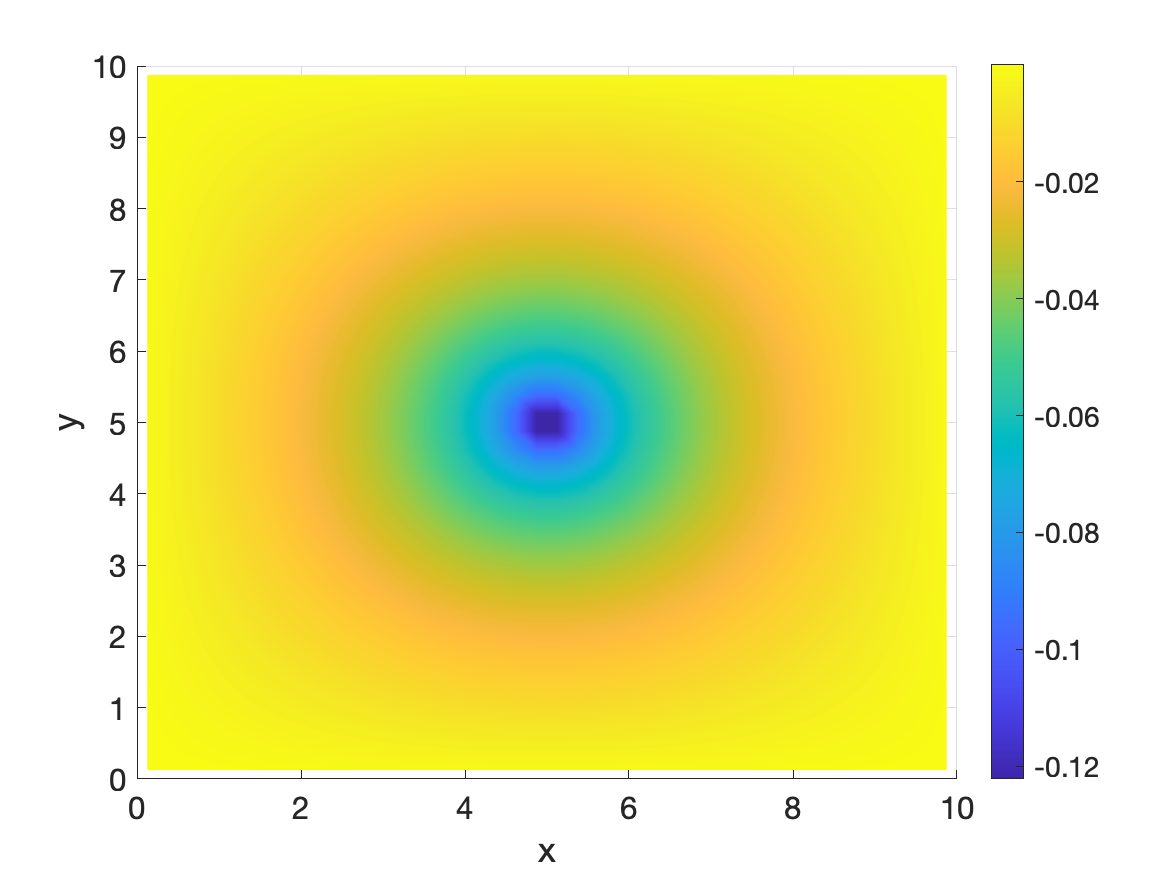}
 \includegraphics[width=0.4\textwidth]{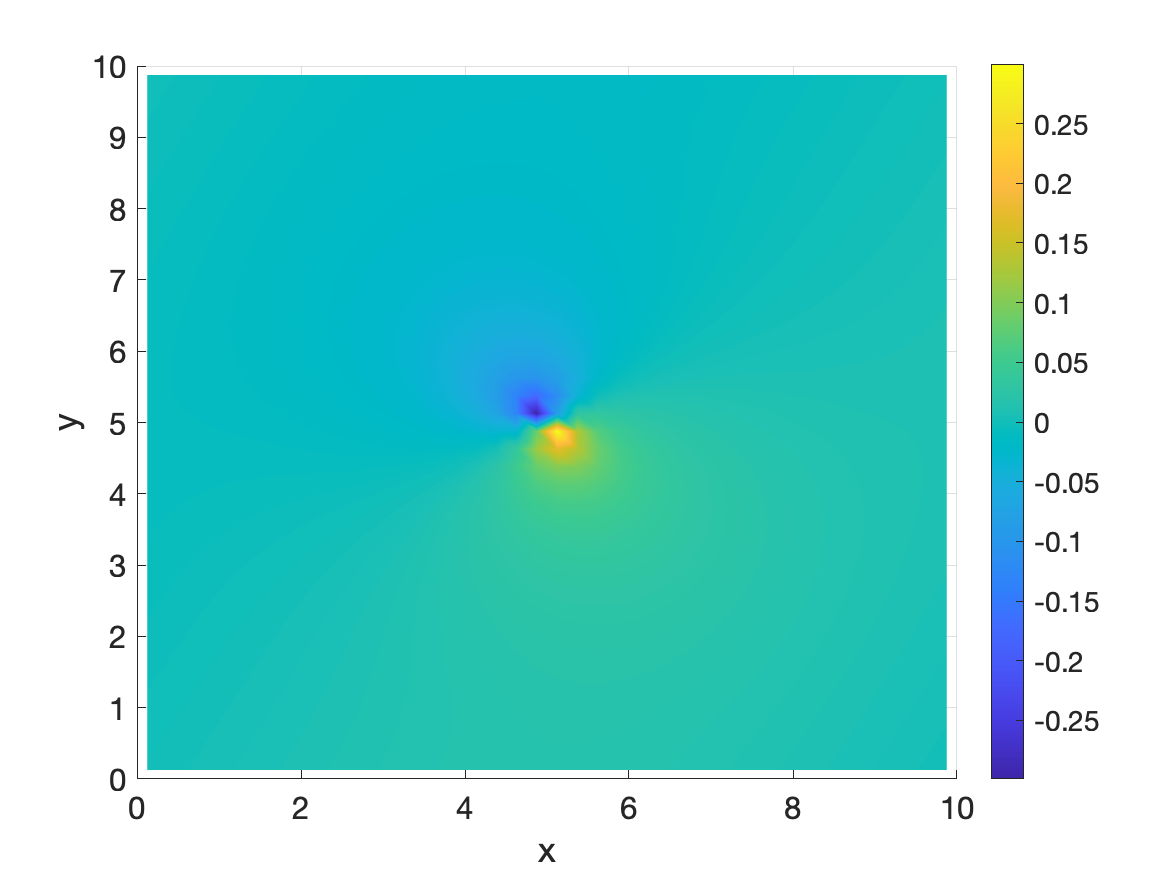}
 \includegraphics[width=0.4\textwidth]{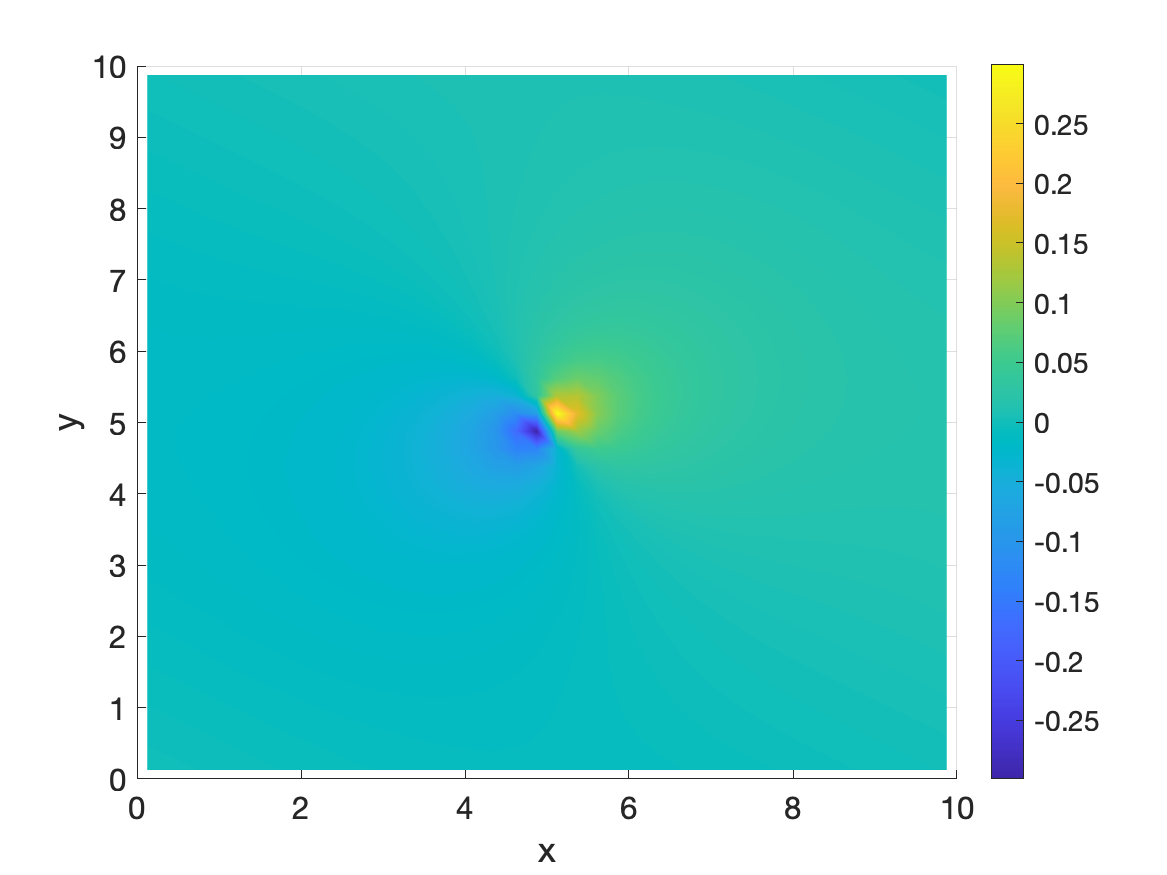}
 \includegraphics[width=0.4\textwidth]{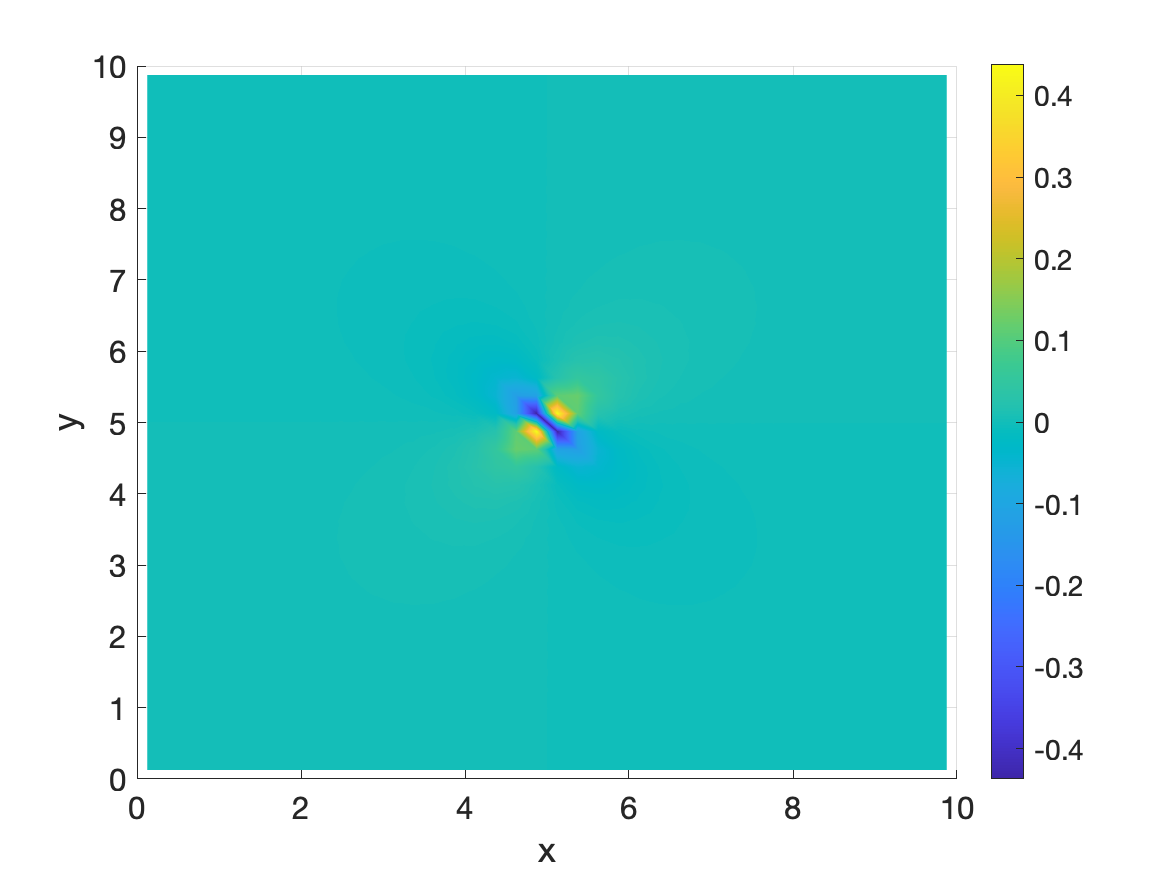}
 \caption{\pzc{First four reduced basis, \fli{ordered from left to right and from top to bottom}, for 2D Example 2 (scattering dominant).\label{fig:2d_basis_test2}}}
\end{figure}

\fli{
In Figures \ref{fig:2d_basis_test1}-\ref{fig:2d_basis_test2}, we present the first four reduced basis functions of 2D Examples 1 and 2 from Section \ref{sec:numericaL^2d}, again after the SVD orthogonalization step. The leading RB function captures the overall configuration of the density as in Figure \ref{fig:2d_test_rho}, while the remaining RB functions encode various multipole structures.}




\end{appendices}

\renewcommand\refname{References}
\bibliographystyle{siam}
\bibliography{bib_peng}

\begin{thebibliography}{10}

\bibitem{adams2001discontinuous}
{\sc M.~L. Adams}, {\em {Discontinuous finite element transport solutions in
  thick diffusive problems}}, Nuclear science and engineering, 137 (2001),
  pp.~298--333.

\bibitem{adams2002fast}
{\sc M.~L. Adams and E.~W. Larsen}, {\em {Fast iterative methods for
  discrete-ordinates particle transport calculations}}, Progress in nuclear
  energy, 40 (2002), pp.~3--159.

\bibitem{adams1992diffusion}
{\sc M.~L. Adams and W.~R. Martin}, {\em Diffusion synthetic acceleration of
  discontinuous finite element transport iterations}, Nuclear science and
  engineering, 111 (1992), pp.~145--167.

\bibitem{alberti2020reduced}
{\sc A.~L. Alberti and T.~S. Palmer}, {\em Reduced-order modeling of nuclear
  reactor kinetics using proper generalized decomposition}, Nuclear Science and
  Engineering, 194 (2020), pp.~837--858.

\bibitem{arridge2009optical}
{\sc S.~R. Arridge and J.~C. Schotland}, {\em Optical tomography: forward and
  inverse problems}, Inverse problems, 25 (2009), p.~123010.

\bibitem{bardos1984diffusion}
{\sc C.~Bardos, R.~Santos, and R.~Sentis}, {\em Diffusion approximation and
  computation of the critical size}, Transactions of the american mathematical
  society, 284 (1984), pp.~617--649.

\bibitem{behne2019model}
{\sc P.~Behne, J.~Ragusa, and J.~Morel}, {\em Model-order reduction for sn
  radiation transport}, in ANS International Conference on Mathematics and
  Computation (M\&C). Portland, OR, USA, 2019.

\bibitem{BinevCohenDahmenDevorePetrovaWojtaszczyk}
{\sc P.~Binev, A.~Cohen, W.~Dahmen, R.~Devore, G.~Petrova, and P.~Wojtaszczyk},
  {\em {Convergence rates for greedy algorithms in reduced basis methods}},
  SIAM Journal on Mathematical Analysis.

\bibitem{buchan2015pod}
{\sc A.~G. Buchan, A.~Calloo, M.~G. Goffin, S.~Dargaville, F.~Fang, C.~C. Pain,
  and I.~M. Navon}, {\em {A POD reduced order model for resolving angular
  direction in neutron/photon transport problems}}, Journal of Computational
  Physics, 296 (2015), pp.~138--157.

\bibitem{chen2020random}
{\sc K.~Chen, Q.~Li, J.~Lu, and S.~J. Wright}, {\em Random sampling and
  efficient algorithms for multiscale pdes}, SIAM Journal on Scientific
  Computing, 42 (2020), pp.~A2974--A3005.

\bibitem{ChenSigalJiMaday2021}
{\sc Y.~Chen, S.~Gottlieb, L.~Ji, and Y.~Maday}, {\em An {EIM}-degradation free
  reduced basis method via over collocation and residual hyper reduction-based
  error estimation}, arXiv preprint arXiv:2101.05902.

\bibitem{ChenGottliebMaday}
{\sc Y.~Chen, S.~Gottlieb, and Y.~Maday}, {\em Parametric analytical
  preconditioning and its applications to the reduced collocation methods.}, C.
  R. Acad. Sci. Paris, Ser. I, 352 (2014), pp.~661 -- 666.

\bibitem{ChenJiNarayanXu2021}
{\sc Y.~Chen, L.~Ji, A.~Narayan, and Z.~Xu}, {\em L1-based reduced over
  collocation and hyper reduction for steady state and time-dependent nonlinear
  equations}.

\bibitem{chen2019robust}
{\sc Y.~Chen, J.~Jiang, and A.~Narayan}, {\em A robust error estimator and a
  residual-free error indicator for reduced basis methods}, Computers \&
  Mathematics with Applications, 77 (2019), pp.~1963--1979.

\bibitem{choi2020space}
{\sc Y.~Choi, P.~Brown, W.~Arrighi, R.~Anderson, and K.~Huynh}, {\em
  Space--time reduced order model for large-scale linear dynamical systems with
  application to boltzmann transport problems}, Journal of Computational
  Physics, 424 (2020), p.~109845.

\bibitem{coale2019reduced}
{\sc J.~Coale and D.~Y. Anistratov}, {\em A reduced-order model for thermal
  radiative transfer problems based on multilevel quasidiffusion method}, in
  International Conference on Mathematics and Computational Methods Applied to
  Nuclear Science and Engineering, M and C 2019, 2019, pp.~278--287.

\bibitem{ding2019error}
{\sc Z.~Ding, L.~Einkemmer, and Q.~Li}, {\em Error analysis of an asymptotic
  preserving dynamical low-rank integrator for the multi-scale radiative
  transfer equation}, arXiv preprint arXiv:1907.04247,  (2019).

\bibitem{dominesey2019reduced}
{\sc K.~A. Dominesey and W.~Ji}, {\em Reduced-order modeling of neutron
  transport separated in space and angle via proper generalized decomposition},
  in ANS International Conference on Mathematics and Computation (M\&C).
  Portland, OR, USA, 2019.

\bibitem{dominesey2018reduced}
{\sc K.~A. Dominesey and J.~P. Senecal}, {\em A reduced-order neutron transport
  model separated in space and angle}, Transactions, 119 (2018), pp.~687--690.

\bibitem{einkemmer2020asymptotic}
{\sc L.~Einkemmer, J.~Hu, and Y.~Wang}, {\em An asymptotic-preserving dynamical
  low-rank method for the multi-scale multi-dimensional linear transport
  equation}, arXiv preprint arXiv:2005.06571,  (2020).

\bibitem{fuhrer1997posteriori}
{\sc C.~F{\"u}hrer and G.~Kanschat}, {\em A posteriori error control in
  radiative transfer}, Computing, 58 (1997), pp.~317--334.

\bibitem{greif2019decay}
{\sc C.~Greif and K.~Urban}, {\em {Decay of the Kolmogorov N-width for wave
  problems}}, Applied Mathematics Letters, 96 (2019), pp.~216--222.

\bibitem{guermond2010asymptotic}
{\sc J.-L. Guermond and G.~Kanschat}, {\em {Asymptotic analysis of upwind
  discontinuous Galerkin approximation of the radiative transport equation in
  the diffusive limit}}, SIAM Journal on Numerical Analysis, 48 (2010),
  pp.~53--78.

\bibitem{haasdonk2017reduced}
{\sc B.~Haasdonk}, {\em {Reduced basis methods for parametrized PDEs--a
  tutorial introduction for stationary and instationary problems}}, Model
  reduction and approximation: theory and algorithms, 15 (2017), p.~65.

\bibitem{hartmann2003adaptive}
{\sc R.~Hartmann and P.~Houston}, {\em Adaptive discontinuous galerkin finite
  element methods for nonlinear hyperbolic conservation laws}, SIAM Journal on
  Scientific Computing, 24 (2003), pp.~979--1004.

\bibitem{hesthaven2016certified}
{\sc J.~S. Hesthaven, G.~Rozza, B.~Stamm, et~al.}, {\em Certified reduced basis
  methods for parametrized partial differential equations}, vol.~590, Springer,
  2016.

\bibitem{jin2010asymptotic}
{\sc S.~Jin}, {\em {Asymptotic preserving (AP) schemes for multiscale kinetic
  and hyperbolic equations: a review}}, Lecture Notes for Summer School on
  “Methods and Models of Kinetic Theory (M\&MKT), Porto Ercole (Grosseto,
  Italy),  (2010), pp.~177--216.

\bibitem{larsen1989asymptotic}
{\sc E.~W. Larsen and J.~E. Morel}, {\em Asymptotic solutions of numerical
  transport problems in optically thick, diffusive regimes ii},  (1989).

\bibitem{larsen1987asymptotic}
{\sc E.~W. Larsen, J.~E. Morel, and W.~F. Miller~Jr}, {\em {Asymptotic
  solutions of numerical transport problems in optically thick, diffusive
  regimes}}, Journal of Computational Physics, 69 (1987), pp.~283--324.

\bibitem{lathrop1968ray}
{\sc K.~D. Lathrop}, {\em Ray effects in discrete ordinates equations}, Nuclear
  Science and Engineering, 32 (1968), pp.~357--369.

\bibitem{lewis1984computational}
{\sc E.~E. Lewis and W.~F. Miller}, {\em Computational methods of neutron
  transport},  (1984).

\bibitem{LiuChenChenShu2019}
{\sc Y.~Liu, T.~Chen, Y.~Chen, and C.-W. Shu}, {\em Certified offline-free
  reduced basis (cofrb) methods for stochastic differential equations driven by
  arbitrary types of noise}, Journal of Scientific Computing,  (2019).

\bibitem{lorence1989s}
{\sc L.~J. Lorence~Jr, J.~Morel, and E.~W. Larsen}, {\em {An $S_2$ synthetic
  acceleration scheme for the one-dimensional $S_n$ equations with linear
  discontinuous spatial differencing}}, Nuclear Science and Engineering, 101
  (1989), pp.~341--351.

\bibitem{mcclarren2019calculating}
{\sc R.~G. McClarren}, {\em Calculating time eigenvalues of the neutron
  transport equation with dynamic mode decomposition}, Nuclear Science and
  Engineering, 193 (2019), pp.~854--867.

\bibitem{mihalas2013foundations}
{\sc D.~Mihalas and B.~W. Mihalas}, {\em Foundations of radiation
  hydrodynamics}, Courier Corporation, 2013.

\bibitem{naldi1998numerical}
{\sc G.~Naldi and L.~Pareschi}, {\em Numerical schemes for kinetic equations in
  diffusive regimes}, Applied mathematics letters, 11 (1998), pp.~29--35.

\bibitem{ohlberger2015reduced}
{\sc M.~Ohlberger and S.~Rave}, {\em {Reduced basis methods: Success,
  limitations and future challenges}}, arXiv:1511.02021,  (2015).

\bibitem{peng2020low}
{\sc Z.~Peng, R.~G. McClarren, and M.~Frank}, {\em A low-rank method for
  two-dimensional time-dependent radiation transport calculations}, Journal of
  Computational Physics, 421 (2020), p.~109735.

\bibitem{pinkus_n-widths_1985}
{\sc A.~Pinkus}, {\em {N-widths in approximation theory}}, Springer, 1985.

\bibitem{pomraning1973equations}
{\sc G.~C. Pomraning}, {\em {The equations of radiation hydrodynamics}},
  International Series of Monographs in Natural Philosophy, Oxford: Pergamon
  Press,  (1973).

\bibitem{pomraning2005equations}
{\sc G.~C. Pomraning}, {\em The equations of radiation hydrodynamics}, Courier
  Corporation, 2005.

\bibitem{prince2019separated}
{\sc Z.~PRINCE and J.~RAGUSA}, {\em Separated representation of spatial
  dimensions in sn neutron transport using the proper generalized
  decomposition}, in ANS International Conference on Mathematics and
  Computation (M\&C). Portland, OR, USA, 2019.

\bibitem{quarteroni2015reduced}
{\sc A.~Quarteroni, A.~Manzoni, and F.~Negri}, {\em Reduced basis methods for
  partial differential equations: an introduction}, vol.~92, Springer, 2015.

\bibitem{quarteroni2011certified}
{\sc A.~Quarteroni, G.~Rozza, and A.~Manzoni}, {\em Certified reduced basis
  approximation for parametrized partial differential equations and
  applications}, Journal of Mathematics in Industry, 1 (2011), p.~3.

\bibitem{SantoDeparisManzoniQuarteroni2018}
{\sc N.~D. Santo, S.~Deparis, A.~Manzoni, and A.~Quarteroni}, {\em Multi space
  reduced basis preconditioners for large-scale parametrized pdes}, SIAM
  Journal on Scientific Computing, 40 (2018), pp.~A954--A983.

\bibitem{tencer2016reduced}
{\sc J.~Tencer, K.~Carlberg, R.~Hogan, and M.~Larsen}, {\em {Reduced order
  modeling applied to the discrete ordinates method for radiation heat transfer
  in participating media}}, in ASME 2016 Heat Transfer Summer Conference
  collocated with the ASME 2016 Fluids Engineering Division Summer Meeting and
  the ASME 2016 14th International Conference on Nanochannels, Microchannels,
  and Minichannels, American Society of Mechanical Engineers Digital
  Collection, 2016.

\bibitem{wareing1993new}
{\sc T.~Wareing}, {\em {New diffusion-sythetic acceleration methods for the
  $S_N$ equations with corner balance spatial differencing}},  (1993).

\bibitem{ZahmNouy2016}
{\sc O.~Zahm and A.~Nouy}, {\em Interpolation of inverse operators for
  preconditioning parameter-dependent equations}, SIAM Journal on Scientific
  Computing, 38 (2016), pp.~A1044--A1074.

\end{thebibliography}
\end{document}